\documentclass[12pt]{article}

\usepackage{amsmath}
\usepackage{amsfonts}
\usepackage{amssymb}
\usepackage{amsthm}
\usepackage{graphicx}

\usepackage{cite}
\usepackage{caption}
\usepackage{url}
\usepackage{color}\usepackage{placeins}
\usepackage{svg}

\usepackage{algorithm}
\usepackage{algpseudocode}
\usepackage{mathtools}
\usepackage{subcaption}
\usepackage{multirow}
\usepackage{bbm}
\usepackage{epstopdf}
\usepackage{graphicx}
\usepackage{bm}
\usepackage{enumerate}

\theoremstyle{definition}

\newtheorem{remark}{Remark}
\newcommand{\norm}[1]{\left\lVert#1\right\rVert}

\DeclareMathOperator*{\argmax}{arg\,max}

\title{FC-based shock-dynamics solver with neural-network localized artificial-viscosity assignment} 

\author{Oscar P. Bruno \footnote{Computing and Mathematical Sciences, Caltech,
    Pasadena, CA 91125, USA} \and Jan S. Hesthaven \footnote{Computational Mathematics and Simulation Science (MCSS), Ecole Polytechnique Federale de Lausanne, CH-1015 Lausanne, Switzerland} \and Daniel V. Leibovici $^*$}
    
\begin{document}
\date{}
\maketitle

\begin{abstract}
  This paper presents a spectral scheme for the numerical solution of
  nonlinear conservation laws in {\em non-periodic domains under
    arbitrary boundary conditions}. The approach relies on the use of
  the Fourier Continuation (FC) method for spectral representation of
  non-periodic functions in conjunction with smooth localized
  artificial viscosity assignments produced by means of a
  Shock-Detecting Neural Network (SDNN). Like previous shock capturing
  schemes and artificial viscosity techniques, the combined FC-SDNN
  strategy effectively controls spurious oscillations in the proximity
  of discontinuities. Thanks to its use of a {\em localized but smooth
    artificial viscosity term}, whose support is restricted to a
  vicinity of flow-discontinuity points, the algorithm enjoys spectral
  accuracy and low dissipation away from flow discontinuities, and, in
  such regions, it produces smooth numerical solutions---as evidenced
  by an essential absence of spurious oscillations in level set
  lines. The FC-SDNN viscosity assignment, which does not require use
  of problem-dependent algorithmic parameters, induces a significantly
  lower overall dissipation than other methods, including the
  Fourier-spectral versions of the previous entropy viscosity
  method. The character of the proposed algorithm is illustrated with
  a variety of numerical results for the linear advection, Burgers and
  Euler equations in one and two-dimensional non-periodic spatial
  domains.

\end{abstract}

\vspace{0.5 cm}
\maketitle
\newpage
\section{\label{sec:introduction}Introduction}
This paper presents a new ``FC-SDNN'' spectral scheme for the
numerical solution of nonlinear conservation laws under arbitrary
boundary conditions.  The proposed approaches relies on use of the
FC-Gram Fourier Continuation method~\cite{albin2011spectral,
  amlani2016fc,bruno2010high, lyon2010high} for spectral
representation of non-periodic functions in conjunction with localized
smooth artificial viscosity assignments prescribed by means of the
neural network-based shock-detection method proposed
in~\cite{schwander2021controlling}. The neural network
approach~\cite{schwander2021controlling} itself utilizes Fourier
series to discretize the gas dynamics and related equations, and it
eliminates Gibbs ringing at shock positions (which are determined by
means of an artificial neural network) by assigning artificial
viscosity over a small number of discrete points in a close vicinity
of shocks.  The use of the classical Fourier spectral method in that
contribution restricts the method's applicability to periodic problems
(so that, in particular, the outer computational boundaries cannot be
physical boundaries), and its highly localized viscosity assignments
give rise to a degree of non-smoothness, resulting in certain types of
unphysical oscillations manifested as serrated level-set lines in the
flow fields.  The FC-SDNN method presented in this paper addresses
these challenges. In particular, in view of its use of certain
newly-designed {\em smooth viscosity windows} introduced in
Section~\ref{subsec:artificial_viscosity_strategy}, the method avoids
the introduction of roughness in the viscosity assignments and thus it
yields smooth flows away from shocks and other flow discontinuities.
In addition, the underlying FC-Gram spectral representations enable
applicability to general non-periodic problems, and, in view of the
weak local viscosity assignments used, it gives rise to sharp
resolution of shocks.  As a result, and as demonstrated in this paper
via application to a range of well known 1D and 2D shock-wave test
configurations, the overall FC-SDNN approach yields accurate and
essentially oscillation-free solutions for general non-periodic
problems.

The computational solution of systems of conservation laws has been
tackled by means of a variety of numerical methods, including low-order
finite volume~\cite{leveque1992numerical, leveque2002finite} and finite difference methods equipped with slope
limiters~\cite{leveque1992numerical, leveque2002finite}, as well as
higher order shock-capturing methods such as the
ENO~\cite{harten1987uniformly} and WENO schemes~\cite{liu1994weighted,
  jiang1996efficient}. An efficient FC/WENO hybrid solver was proposed
in~\cite{shahbazi2011multi}.  The use of artificial viscosity as a
computational device for conservation laws, on the other hand, was
first proposed in~\cite{richtmyer1948proposed, vonneumann1950method}
and the subsequent contributions\cite{lapidus1967detached,
  lax1959systems, gentry1966eulerian}. The viscous terms proposed in
these papers, which incorporate derivatives of the square of the
velocity gradient, may induce oscillations in the vicinity of
shocks~\cite{lax1959systems} (since the velocity itself is not smooth
in such regions), and, as they do not completely vanish away from the
discontinuities, they may lead to significant approximation errors in
regions were the fluid velocity varies
rapidly. Reference~\cite{persson2006sub} proposes the use of a
shock-detecting sensor in order to localize the support of the
artificial viscosity, which is used in the context of a Discontinuous
Galerkin scheme.

The entropy viscosity method~\cite{guermond2011entropy} (EV)
incorporates a nonlinear viscous ``entropy-residual'' term which
essentially vanishes away from discontinuities---in view of the fact
that the flow is isentropic over smooth flow regions---and which is
thus used to limit non-zero viscosity assignments to regions near
flow discontinuities, including both shock waves and contact
discontinuities. This method, however, relies on several
problem-dependent algorithmic parameters that require tuning for every
application. Additionally, this approach gives rise to a significant
amount of dissipation even away from shocks, in particular in the
vicinity of contact discontinuities and regions containing fast
spatial variation in the flow-field variables. Considerable
improvements concerning this issue were obtained in
\cite{kornelus2018flux} (which additionally introduced a Hermite-based
method to discretize the hyperbolic systems) by modifying the EV
viscosity term. Like the viscous term introduced
by~\cite{vonneumann1950method}, the EV viscosity
assignments~\cite{guermond2011entropy,kornelus2018flux} are themselves
discontinuous in the vicinity of shocks, and, thus, their use may
introduce spurious oscillations. The
C-method~\cite{reisner2013space,ramani2019space1, ramani2019space2},
which augments the hyperbolic system with an additional equation used
to determine a spatio-temporally smooth viscous term, relies, like the
EV method, on use of several problem dependent parameters and
algorithmic variations.

Recently, significant progress was made by incorporating machine
learning-based techniques (ML) to enhance the performance of classical
shock capturing
schemes~\cite{ray2018artificial,discacciati2020controlling,stevens2020enhancement,schwander2021controlling}. The
approach~\cite{ray2018artificial,discacciati2020controlling,schwander2021controlling}
utilizes ML-based methods to detect discontinuities which are then
smeared by means of shock-localized artificial viscosity assignments
in the context of various discretization methods, including
Discontinuous Galerkin
schemes~\cite{ray2018artificial,discacciati2020controlling} and
Fourier spectral schemes~\cite{schwander2021controlling}. The ML-based
approach utilized in~\cite{stevens2020enhancement}, in turn, modifies
the finite volume coefficients utilized in the WENO5-JS scheme by
learning small perturbations of these coefficients leading to improved
accurate representations of functions at cell boundaries.

Like the strategy underlying the
contribution~\cite{schwander2021controlling}, the FC-SDNN method
proposed in this paper relies on the occurrence of Gibbs oscillations
for ML-based shock detection. Unlike the previous approach, however,
the present method assigns a {\em smooth} (albeit also
shock-localized) viscous term. In view of its smooth viscosity
assignments this procedure effectively eliminates Gibbs oscillations
while avoiding introduction of the flow-field roughness that is often
evidenced by the serrated level sets produced by other methods. In
view of its use of FC-based Fourier expansions, further, the proposed
algorithm enjoys spectral accuracy away from shocks (thus, delivering,
in particular, essentially vanishing dispersion in such regions; see
Section~\ref{sec:fc} and Figure~\ref{fig: LA dispersion}) while
enabling solution under general (non-periodic) boundary
conditions. Unlike other techniques, finally, the approach does not
require use of problem-dependent algorithmic parameters.

The capabilities of the proposed algorithm are illustrated by means of
a variety of numerical results, in one and two-dimensional contexts,
for the Linear Advection, Burgers and Euler equations.  In order to
provide a useful reference point, this paper also presents an FC-based
version, termed FC-EV, of the EV
algorithm~\cite{guermond2011entropy}. (The modified
version~\cite{kornelus2018flux} of the entropy viscosity approach,
which was also tested as a possible reference solver, was not found to
be completely reliable in our FC-based context, since it
occasionally led to spurious oscillations in shock regions as grids
were refined, and the corresponding results were therefore not
included in this paper.) We find that the FC-SDNN algorithm generally
provides significantly more accurate numerical approximations than the
FC-EV, as the localized artificial viscosity in the former approach
induces a much lower dissipation level than the latter.

This paper is organized as follows. After necessary preliminaries are
presented in Section~\ref{Preliminaries} (concerning the hyperbolic
problems under consideration, as well as the Fourier Continuation
method, and including basic background on the artificial-viscosity
strategies we consider), Section~\ref{sec:FC-SDNN} describes the
proposed FC-SDNN approach. Section~\ref{sec:numerical results} then
demonstrates the algorithm's performance for a variety of non-periodic
linear and nonlinear hyperbolic problems. In particular, cases are
considered for the linear advection equation, Burgers equation and
Euler's equations in one-dimensional and two-dimensional rectangular
and non-rectangular spatial domains, including cases in which shock
waves meet smooth and non-smooth physical boundaries.

\section{\label{Preliminaries} Preliminaries}

\subsection{Conservation laws\label{laws}}
This paper proposes novel spectral methodologies, applicable in
general non-periodic contexts and with general boundary conditions,
for the numerical solution of conservation-law equations of the form
\begin{equation} \label{eq: nonlinear conservative eqn}
    \dfrac{\partial}{\partial t} \textit{\textbf{e}}(\textit{\textbf{x}}, t) + \nabla \cdot \big( f(\textit{\textbf{e}}(\textit{\textbf{x}}, t))\big) = 0
\end{equation}
on a bounded domain $\Omega\subset \mathbb{R}^q $, where
$\textit{\textbf{e}} : \Omega \times [0, T] \rightarrow \mathbb{R}^r$
and $f : \mathbb{R}^r \rightarrow \mathbb{R}^r \times \mathbb{R}^q$
denote the unknown solution vector and a (smooth) convective flux,
respectively.

The proposed spectral approaches are demonstrated for several
equations of the form~\eqref{eq: nonlinear conservative eqn},
including the Linear Advection equation
\begin{equation} \label{eq: linear advection}
    \frac{\partial u}{\partial t} + a\frac{\partial u}{\partial x} = 0
\end{equation}
with a constant propagation velocity $a$, where we have
$\textit{\textbf{e}} = u$, and $f(u) = a u$; the one- and
two-dimensional scalar Burgers equations
\begin{equation} \label{eq: burgers 1d}
     \frac{\partial u}{\partial t} + \frac{1}{2}\Big(\frac{\partial u^2}{\partial x}\Big) = 0
\end{equation}
and
\begin{equation}\label{eq: burgers 2d} \frac{\partial u}{\partial t}
  + \frac{1}{2}\Big(\frac{\partial u^2}{\partial x}\Big) +
  \frac{1}{2}\Big(\frac{\partial u^2}{\partial y}\Big) = 0,
\end{equation}
for each of which we have $\textit{\textbf{e}} = u$ and
$f(u) = \frac{u^2}{2}$; as well as the one- and two-dimensional Euler
equations
\begin{equation} \label{eq: euler 1d equation}
\dfrac{\partial}{\partial t}
\begin{pmatrix}
\ \rho \\[\jot]
\ \rho u\\[\jot]
\ E
\end{pmatrix} +  \frac{\partial}{\partial x}
\begin{pmatrix}
\ \rho u\\[\jot]
\ \rho u^2 + p\\[\jot]
\ u (E + p)
\end{pmatrix} = 0
\end{equation}
and
\begin{equation} \label{eq: euler 2d equation}
\dfrac{\partial}{\partial t}
\begin{pmatrix}
\ \rho \\[\jot]
\ \rho u\\[\jot]
\ \rho v\\[\jot]
\ E
\end{pmatrix} + \frac{\partial}{\partial x}
\begin{pmatrix}
\ \rho u\\[\jot]
\ \rho u^2 + p\\[\jot]
\ \rho u v\\[\jot]
\ u (E + p)
\end{pmatrix}
+ \frac{\partial}{\partial y}
\begin{pmatrix}
\ \rho v\\[\jot]
\ \rho u v\\[\jot]
\ \rho v^2 + p\\[\jot]
\ v (E + p)
\end{pmatrix} = 0
\end{equation}
with
\begin{equation}\label{eq: tp}
    E = \frac{p}{\gamma - 1} + \frac{1}{2} \rho \lvert \textit{\textbf{u}} \rvert ^2,
\end{equation}
for each of which we have
\begin{equation}\label{flow_flux}
  \textit{\textbf{e}} = (\rho, \rho \textit{\textbf{u}}, E)^T,\quad 
  f(\textit{\textbf{e}}) = (\rho \textit{\textbf{u}}, \rho
  \textit{\textbf{u}} \otimes \textit{\textbf{u}} + p\mathbb{I},
  \textit{\textbf{u}}(E + p))^T.
\end{equation}
Here $\mathbb{I}$ denotes the identity tensor,
$\textbf{a}\otimes \textbf{b} = (a_ib_j)$ denotes the tensor product
of the vectors $\textbf{a} = (a_i)$ and $\textbf{b} = (b_j)$, and
$\rho$, \textit{\textbf{u}}, $E$ and $p$ denote the density, velocity
vector, total energy and pressure, respectively. The speed of
sound~\cite{leveque1992numerical}
\begin{equation}\label{snd_sp}
  a = \sqrt{\frac{\gamma p}{\rho}}
\end{equation}
for the Euler equations plays important roles in the various
artificial viscosity assignments considered in this paper for Euler
problems in both 1D and 2D.
\begin{remark} 
  As an example concerning notational conventions, note that in the
  case of the 2D Euler equations, for which $f$ is given
  by~\eqref{flow_flux},
  $\nabla \cdot \big(f(\textit{\textbf{e}})\big)$ can be viewed as a
  three coordinate vector whose first, second and third coordinates
  are a scalar, a vector and a scalar, respectively. Using the
  Einstein summation convention, these three components are
  respectively given by
  $\nabla\cdot(\rho\textbf{u}) = \partial_i (\rho u_i)$,
  $\big( \nabla\cdot (\rho \textit{\textbf{u}} \otimes
  \textit{\textbf{u}} + p\mathbb{I}) \big)_j =\partial_i (\rho u_j u_i
  + p)$ and
  $\nabla\cdot((E + p)\textbf{u}) = \partial_i ((E + p) u_i)$.
\end{remark}

\subsection{Artificial viscosity\label{art_vis}}

As is well known, the shocks and other flow discontinuities that arise
in the context of nonlinear conservation laws of the form~(\ref{eq:
  nonlinear conservative eqn}) give rise to a number of challenges
from the point of view of computational simulation. In particular, in
the framework of classical finite difference methods as well as
Fourier spectral methods, such discontinuities are associated with the
appearance of spurious ``Gibbs oscillations''. Artificial viscosity
methods aim at tackling this difficulty by considering, instead of the
inviscid equations~(\ref{eq: nonlinear conservative eqn}), certain
closely related equations which include viscous terms containing
second order spatial derivatives. Provided the viscous terms are
adequately chosen and sufficiently small, the resulting solutions,
which are smooth functions on account of viscosity, approximate well
the desired (discontinuous) inviscid solutions.  In general terms, the
viscous equations are obtained by adding a viscous term of the form
$\nabla \cdot \big( f_\textit{visc}[\textit{\textbf{e}}]\big)$ to the
right hand side of (\ref{eq: nonlinear conservative eqn}), where the
``viscous flux'' operator
\begin{equation}\label{f_visc}
f_\textit{visc}[\textit{\textbf{e}}] = \mu[\textit{\textbf{e}}] \mathbf{D}[\textit{\textbf{e}}],
\end{equation}
(which, for a given vector-valued function $\textbf{e}(x,t)$, produces
a vector-valued function $f_\textit{visc}[\textit{\textbf{e}}](x,t)$
defined in the complete computational domain), is given in terms of a
certain ``viscosity'' operator $\mu[\textit{\textbf{e}}](x,t)$ (which
may or may not include derivatives of the flow variables
$\textbf{e}$), and a certain matrix-valued first order differential
operator $\mathbf{D}$. Once such a viscous term is included, the
viscous equation
\begin{equation} \label{eq: convection diffusion eqn}
  \frac{\partial \textit{\textbf{e}}(\textit{\textbf{x}}, t)}{\partial t} + \nabla \cdot
  \big( f(\textit{\textbf{e}}(\textit{\textbf{x}}, t))\big) = \nabla \cdot
  \big( f_\textit{visc}[\textit{\textbf{e}}](\textit{\textbf{x}}, t)\big)
\end{equation}
results.

Per the discussion in Section~\ref{sec:introduction}, this paper
exploits and extends, in the context of the Fourier-Continuation
discretizations, two different approaches to
viscosity-regularization---each one resulting from a corresponding
selection of the operators $\mu$ and $\mathbf{D}$.  One of these
approaches, the EV method, produces a viscosity assignment
$\mu[\textit{\textbf{e}}](x,t)$ on the basis of certain differential
and algebraic operations together with a number of tunable
problem-dependent parameters that are specifically designed for each
particular equation considered, as described in Section~\ref{EV}. The
resulting viscosity values $\mu[\textit{\textbf{e}}](x,t)$ are highest
in a vicinity of discontinuity regions and decrease rapidly away from
such regions. The neural-network approach introduced in
Section~\ref{SDNN}, in turn, uses machine learning methods to pinpoint
the location of discontinuities, and then produces a viscosity
function whose support is restricted to a vicinity of such
discontinuity locations. As a significant advantage, the
neural-network method, which does not require use of tunable
parameters, is essentially problem independent, and it can use a
single pre-trained neural network for all the equations
considered. Details concerning these two viscosity-assignment methods
considered are provided in what follows.

\subsubsection{Artificial viscosity via shock-detecting neural network
  (SDNN)\label{SDNN}} The SDNN approach proposed in this paper is
based on the neural-network strategy introduced
in~\cite{schwander2021controlling} for detection of discontinuities on
the basis of Gibbs oscillations in Fourier series, together with a novel
selection of the operator $\mu$ in~\eqref{f_visc} that yields
spatially localized but smooth viscosity assignments: per the
description in Section~\ref{subsec:artificial_viscosity_strategy}, the
FC-SDNN viscosity $\mu[\textit{\textbf{e}}](\textit{\textbf{x}}, t)$
is a smooth function that vanishes except in narrow regions around
flow discontinuities.  The differential operator $\mathbf{D}$, in
turn, is simply given by
\begin{equation} \label{eq: viscosity}
  \mathbf{D}[\textit{\textbf{e}}](\textit{\textbf{x}},
  t) = \nabla(
  \textit{\textbf{e}}(\textit{\textbf{x}}, t)), 
\end{equation}
where the gradient is computed component-wise. As indicated in
Section~\ref{sec:introduction}, the smoothness of the proposed
viscosity assignments is inherited by the resulting flows away from
flow discontinuities, thus helping eliminate the serrated level-set
lines that are ubiquitous in the flow patterns produced by other
methods.

\subsubsection{Entropy viscosity methodology (EV)\label{EV}}
The operators $\mu$ and $\mathbf{D}$ employed by the EV
approach~\cite{guermond2011entropy} are defined in terms of a number
of problem dependent functions, vectors and operators. Indeed,
starting with an equation dependent \textit{entropy pair} $(\eta,\nu)$
where $\eta$ is a scalar function and $\nu$ is a vector of the same
dimensionality as the velocity vector, the EV approach utilizes an
associated scalar \textit{entropy residual} operator
\begin{equation} \label{eq: ev residual}
  R_{EV}[\textit{\textbf{e}}](\textit{\textbf{x}}, t) =
  \frac{\partial \eta (\textit{\textbf{e}}(\textit{\textbf{x}},
    t))}{\partial t} + \nabla \cdot \nu
  (\textit{\textbf{e}}(\textit{\textbf{x}}, t))
\end{equation}
together with a function $C = C(\textit{\textbf{e}})$ related to the
local wave speed, and a normalization operator $N = N[e](x,t)$
obtained from the function $\eta$.

In practice, reference~\cite{guermond2011entropy} proposes
$\eta(\textit{\textbf{e}}) = \frac{u^2}{2}$,
$\nu(\textit{\textbf{e}}) = a \frac{u^2}{2}$ and
$C(\textit{\textbf{e}})= a$ for the Linear Advection
equation~\eqref{eq: linear advection},
$\eta(\textit{\textbf{e}}) = \frac{u^2}{2}$,
$\nu(\textit{\textbf{e}}) = \frac{u^3}{3}$ and
$C(\textit{\textbf{e}})= u$ for the 1D and 2D Burgers
equations~\eqref{eq: burgers 1d} and~\eqref{eq: burgers 2d}, and
$\eta(\textit{\textbf{e}}) = \frac{\rho}{\gamma -
  1}\log(p/\rho^{\gamma})$,
$\nu(\textit{\textbf{e}}) = \textit{\textbf{u}} \frac{\rho}{\gamma -
  1}\log(p/\rho^{\gamma})$ and
$C(\textit{\textbf{e}}) = \|\textit{\textbf{u}}\| + a$ (where $a$
denotes the speed of sound~\eqref{snd_sp}) for the 1D and 2D Euler
equations~\eqref{eq: euler 1d equation} and~\eqref{eq: euler 2d
  equation}. As for the normalization operator,
reference~\cite{guermond2011entropy} proposes $N=1$ for the Euler
equations and $N[e](x, t) =|\eta(e)(x,t) - \overline{\eta}(e)(t)|$ for
the Linear advection and Burgers equations, where
$\overline{\eta}(e)(t)$ denotes the spatial average of $\eta(e)$ at
time~$t$.

For a numerical discretization with maximum spatial mesh size
$h$, the EV viscosity function is defined by\looseness = -1
\begin{equation}\label{eq: mu}
  \mu[\textit{\textbf{e}}](\textit{\textbf{x}}, t) = \min(\mu_\textit{max}[\textit{\textbf{e}}](t), \mu_{E}[\textit{\textbf{e}}](\textit{\textbf{x}}, t))
\end{equation}
where the maximum viscosity $\mu_\textit{max}$ is given by
\begin{equation} \label{eq: mu max} \mu_\textit{max}[\textit{\textbf{e}}](t) = c_\textit{max} h \max_{x \in \Omega} |C(\textit{\textbf{e}}(\textit{\textbf{x}}, t))| 
\end{equation}
and where 
\begin{equation} \label{eq: mu ev}
  \mu_{E}[\textit{\textbf{e}}](\textit{\textbf{x}}, t) = c_E h^2
  \frac{|R_{EV}[\textit{\textbf{e}}](\textit{\textbf{x}},
    t)|}{N[\textit{\textbf{e}}](\textit{\textbf{x}}, t)}.
\end{equation}
In particular, the EV viscosity function depends on two parameters,
$c_\textit{max}$ and $c_E$, both of size $\mathcal{O}(1)$, that,
following~\cite{guermond2011entropy}, are to be tuned to each
particular problem.

Finally, the EV differential operator $\mathbf{D}$ for the Linear
Advection and Burgers equations is defined by
\begin{equation}
  \mathbf{D}[\textit{\textbf{e}}](\textit{\textbf{x}}, t) = \nabla(
  \textit{\textbf{e}}(\textit{\textbf{x}}, t) ),
\end{equation}
while for the Euler equations it is given by
\begin{equation}
    \mathbf{D}[\textit{\textbf{e}}](\textit{\textbf{x}}, t) = 
    \begin{pmatrix}
    \ 0 \\[\jot]
    \ \frac{1}{2}(\nabla \textit{\textbf{u}} + (\nabla \textit{\textbf{u}})^T) \\[\jot]
    \ \frac{1}{2}(\nabla \textit{\textbf{u}} + (\nabla \textit{\textbf{u}})^T)  \textit{\textbf{u}}  + \kappa \nabla (p/\rho)
\end{pmatrix}
\end{equation}
where, using the Einstein notation
$\left\{(\nabla \textit{\textbf{u}} + (\nabla \textit{\textbf{u}})^T)
  \textit{\textbf{u}}\right\}_i = (\partial_iu_j +\partial_ju_i)u_j$,
and where $\kappa = \frac{\mathcal{P}}{\gamma - 1} \mu$, with the
Prandtl number $\mathcal{P}$ taken to equal 1.

\subsection{\label{sec:fc} Fourier Continuation spatial approximation}

The straightforward Fourier-based discretization of nonlinear
conservation laws generally suffers from crippling Gibbs oscillations
resulting from two different sources: the physical flow
discontinuities, on one hand, and the overall generic non-periodicity
of the flow variables, on the other. Unlike the Gibbs ringing in
flow-discontinuity regions, the ringing induced by lack of periodicity
is not susceptible to treatment via artificial viscosity assignments
of the type discussed in~\ref{art_vis}. In order to tackle this
difficulty we resort to use of the Fourier Continuation (FC) method
for equispaced-grid spectral approximation of non-periodic functions.
\begin{figure}[H]
\centering
    \includegraphics[width=0.6\linewidth,]{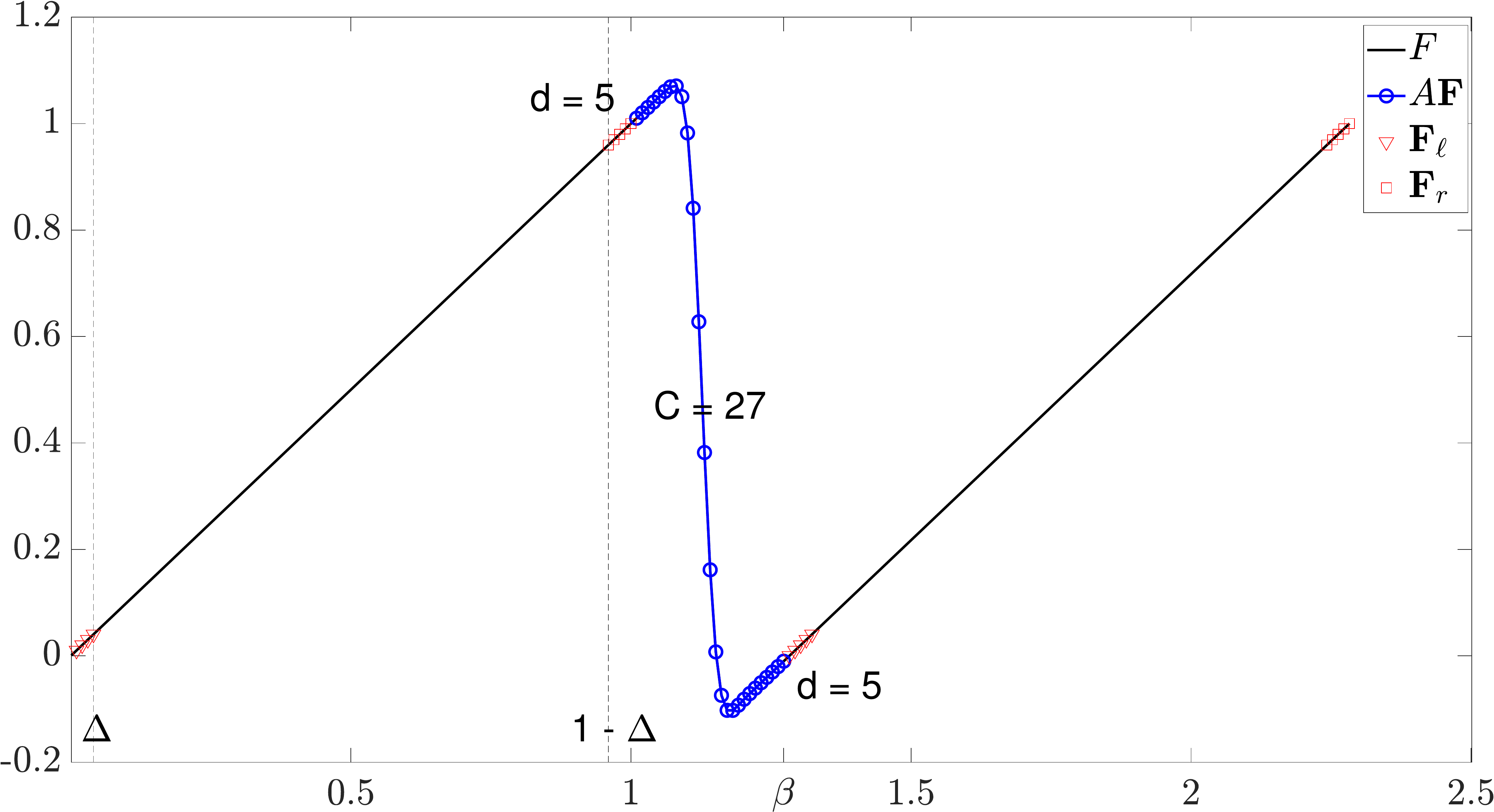}
    \caption{Fourier Continuation of the non-periodic function
      $F(x) = x$ on the interval $[0, 1]$. With reference to the text,
      the red triangles (resp. squares), represent the $d = 5$ left
      (resp. right) matching points, while the blue circles represent
      the $C = 27$ continuation points.}
    \label{fig:fcexample}
\end{figure}

The basic FC algorithm~\cite{albin2011spectral}, called FC-Gram in
view of its reliance on Gram polynomials for near-boundary
approximation, constructs an accurate Fourier approximation of a given
generally non-periodic function $F$ defined on a given one-dimensional
interval---which, for definiteness, is assumed in this section to
equal the unit interval $[0,1]$.  The Fourier continuation algorithm
relies on use of the function values $F_j =F(x_j)$ of the function
$F:[0,1]\to \mathbb{R}$ at $N$ points $x_j=jh\in [0,1]$ ($h=1/(N-1)$)
to produce a function
\begin{equation}\label{FC_exp}
 F^c(x) = \sum^{M}_{k = -M} \hat{F}^c_k \exp (2 \pi i k x /\beta)
\end{equation}
which is defined (and periodic) in an interval $[0,\beta]$ that
strictly contains $[0,1]$, where $\hat{F}^c_k$ denote the FC
coefficients of $F$ and where, as detailed below, $M$ is an integer
that, for $N$ large, is close to (but different from) the integer
$\lfloor N/2\rfloor$.

In order to produce the FC function $F^c$, the FC-Gram algorithm first
uses two subsets of the function values in the vector
$\mathbf{F} = (F_0,\dots,F_{N-1})^T$ (namely the function values at
the``matching points'' $\{x_0,..,x_{d-1}\}$ and
$\{x_{N-d},...,x_{N-1}\}$ located in small matching subintervals
$[0,\Delta]$ and $[1-\Delta,1]$ of length $\Delta = (d-1)h$ near the
left and right ends of the interval $[0,1]$, where $d$ is a small
integer independent of $N$), to produce, at first, a discrete (but
``smooth'') periodic extension vector $\mathbf{F}^c$ of the vector
$\mathbf{F}$. Indeed, using the matching point data, the FC-Gram
algorithm produces and appends a number $C$ of continuation function
values in the interval $[1,\beta]$ to the data vector $\mathbf{F}$, so
that the extension $\mathbf{F}^c$ transitions smoothly from $F_{N-1}$
back to $F_0$, as depicted in Figure~\ref{fig:fcexample}. (The FC
method can also be applied on the basis certain combinations of
function values and derivatives by constructing the continuation
vector $\mathbf{F}^c$, as described below in this section, on the
basis of e.g. the vector
$\mathbf{F} = (F_0,\dots,F_{N-2},F'_{N-1})^T$, where
$F_j\approx F(x_j)$ for $1\leq j\leq N-2$ and where
$F'_{N-1} \approx F'(x_{N-1})$.  Such a procedure enables imposition
of Neumann boundary conditions in the context of the FC method.) The
resulting vector $\mathbf{F}^c$ can be viewed as a discrete set of
values of a smooth and periodic function which can be used to produce
the Fourier continuation function $F^c$ via an application of the FFT
algorithm. The function $F^c$ provides a spectral approximation of $F$
throughout the interval $[0,1]$ which does not suffer from either
Gibbs-ringing or the associated interval-wide accuracy
degradation. Throughout this paper we assume, for simplicity, that
$N+C$ is an odd integer, and, thus, the resulting series has bandwidth
$M = \frac{N + C-1}{2}$; consideration of even values of $N+C$ would
require a slight modification of the index range in~\eqref{FC_exp}.

To obtain the necessary discrete periodic extension $\mathbf{F}^c$,
the FC-Gram algorithm first produces two polynomial interpolants, one
per matching subinterval, using, as indicated above, a small number
$d$ of function values or a combination of function values and a
derivative near each one of the endpoints of the interval $[0,1]$.
This approach gives rise to high-order interpolation of the function
$F$ over the matching intervals $[0,\Delta]$ and $[1-\Delta,1]$. The
method for evaluation of the discrete periodic extension is based on a
representation of these two polynomials in a particular orthogonal
polynomials basis (the Gram polynomials), for each element of which
the algorithm utilizes a precomputed smooth function which blends the
basis polynomial to the zero function over the distance
$\beta-1$~\cite{albin2011spectral,amlani2016fc}. Certain simple
operations involving these ``blending to zero'' functions are then
used, as indicated in these references and as illustrated in
Figure~\ref{fig:fcexample}, to obtain smooth transitions-to-zero from
the left-most and right-most function values to the extension interval
$[1,\beta]$.  The values of this transition function at the points
$N/(N-1),(N+1)/(N-1),\ldots,(N+C-1)/(N-1)$ provide the necessary $C$
additional point values from which, as mentioned above, the discrete
extension $\mathbf{F}^c$ is obtained. The continuation function $F^c$
then easily results via an application of the FFT algorithm to the
function values $F^c$ in the interval $[0,\beta]$.

The discrete continuation procedure can be expressed in the matricial form
\begin{equation}\label{eq: continuation}
    \mathbf{F}^c = 
    \begin{pmatrix}
        \mathbf{F} \\[\jot]
        A_{\ell} Q^T \mathbf{F}_{\ell} + A_r Q^T\mathbf{F}_r 
    \end{pmatrix}
\end{equation}
where the $d$-dimensional vectors $\mathbf{F}_{\ell}$ and
$\mathbf{F}_r$ contain the point values of $F$ at the first and last
$d$ discretization points in the interval $[0,1]$, respectively; where
$Q$ is a $d \times d$ matrix, whose columns contain the point values
of the elements of the Gram polynomial bases on the left matching
intervals; and where $A_{\ell}$ and $A_r$ are $C \times d$, matrices
containing the $C$ values of the blended-to-zero Gram polynomials in
the left and right Gram bases, respectively. These small matrices can
be computed once and stored on disc, and then read for use to produce
FC expansions for functions $G:[a,b]\to \mathbb{R}$ defined on a given
1D interval $[a,b]$, via re-scaling to the interval $[0,1]$.

A minor modification of the procedure presented above suffices to
produce a Fourier continuation function on the basis of data points at
the domain interior and a derivative at interval endpoints. For
example, given the vector
$\mathbf{F} = (F_0,\dots,F_{N-2},F'_{N-1})^T$, using an adequately
modified version $\widetilde Q$ of the matrix $Q$, an FC series $F^c(x)$
can be produced which matches the function values $F_0,\dots,F_{N-2}$
at $x=x_0,\dots,x_{N-2}$, and whose derivative equals $F'_{N-1}$ at
$x_{N-1}$. The matrix $\widetilde Q$ is obtained by using the matrix $Q$
to obtain a value $F_{N-1}$ such that the derivative $F'(x_{N-1})$
equals the given value $F'_{N-1}$. Full details in this regard can be
found in~\cite[Sec. 3.2]{amlani2016fc}.

Clearly, the approximation order of the Fourier Continuation method
(whether derivative values or function values are prescribed at
endpoints) is restricted by the corresponding order $d$ of the Gram
polynomial expansion, which, as detailed in various cases in
Section~\ref{sec:numerical results}, is selected as a small integer:
$d=2$ or $d = 5$.  The relatively low order of accuracy afforded by
the $d=2$ selection, which must be used in some cases to ensure
stability, is not a matter of consequence in the context of the
problems considered in the present paper, where high orders of
accuracy are not expected from any numerical method on account of
shocks and other flow discontinuities. Importantly, even in this
context the FC method preserves one of the most significant numerical
properties of Fourier series, namely its extremely small numerical
dispersion. In fact, with exception of the cyclic advection example
presented in Section~\ref{dispersionless}, for which errors can
accumulate on account of the spatio-temporal periodicity, for all
cases in Section~\ref{sec:numerical results} for which both the $d= 2$
and $d=5$ simulations were performed (which include those presented in
Sections~\ref{bound_1D} (1D linear advection), \ref{subsubsec:2D
  burgers results} (2D Burgers equation) and~\ref{Euler 1D results}
(1D Euler equations), the lower and higher order results obtained were
visually indistinguishable.

\begin{figure}
  \centering
  \includegraphics[width=0.45\textwidth]{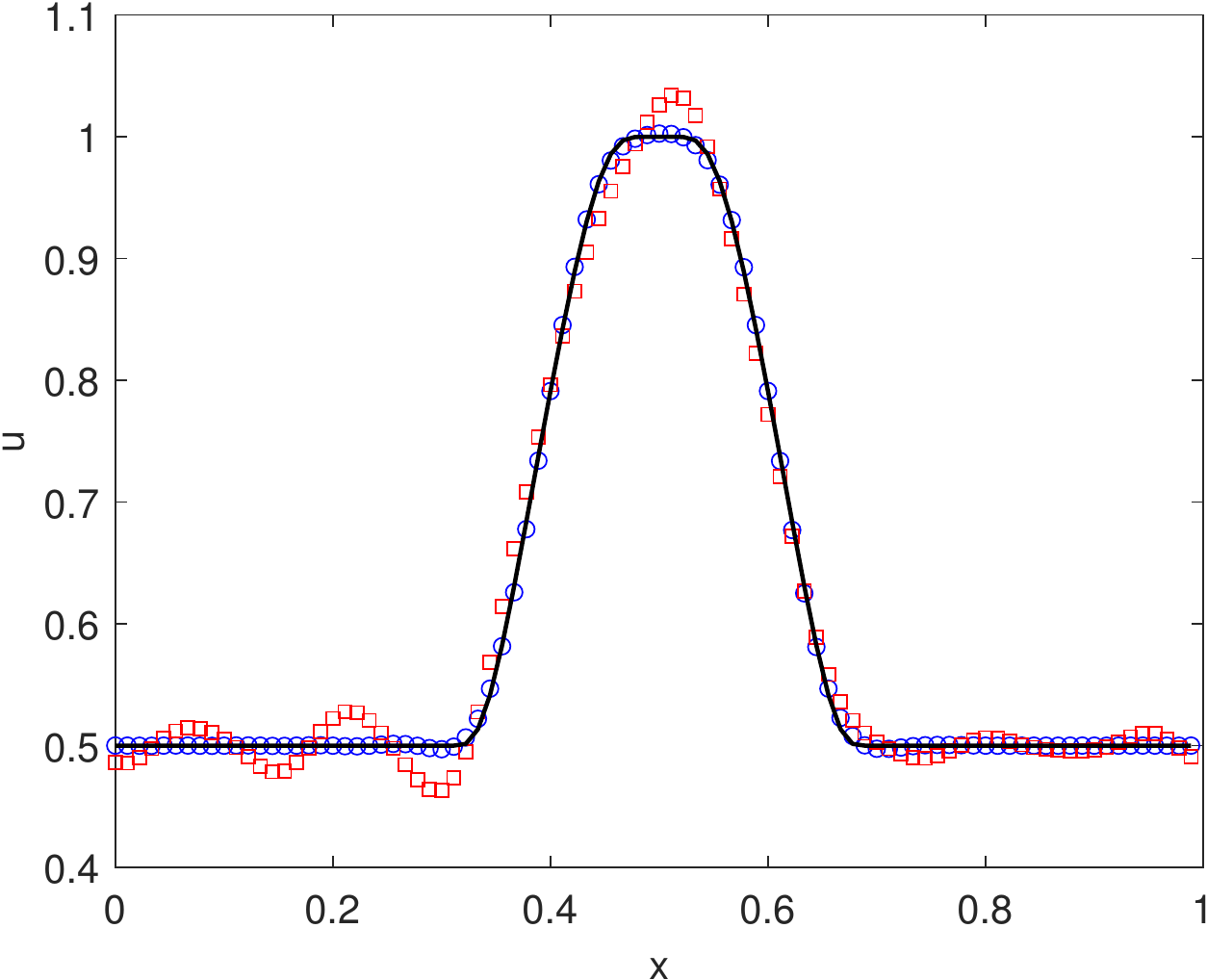}
  \caption{Numerical solutions to the periodic one-dimensional linear
    advection problem with $a=1$ up to time $T = 500$: Exact solution
    (solid black line), fifth order FC-SDNN method (blue circles) and
    sixth order centered finite-difference scheme (red squares). Both
    numerical solutions were obtained using $N = 90$ discretization
    points. In view of its nearly dispersionless character, the
    FC-based solution remains significantly more accurate than its
    higher order finite-difference counterpart.}
  \label{fig: LA dispersion}
\end{figure}

The low dispersion character resulting from use of the FC method is
demonstrated in Figure~\ref{fig: LA dispersion}, which displays
solutions produced by means of two different methods, namely, the
FC-based order-5 FC-SDNN algorithm (Section~\ref{sec:FC-SDNN}) and the
order-6 centered finite-difference scheme (both of which use the
SSPRK-4 time discretization scheme), for a linear advection problem.
The FC-SDNN solution presented in the figure does not deteriorate even
for long propagation times, thus illustrating the essentially
dispersion-free character of the FC-based approach. The
finite-difference solution included in the Figure, in turn, does
exhibit clear degradation with time, owing to the dispersion and
diffusion effects associated with the underlying finite difference
discretization.

\section{\label{sec:FC-SDNN}FC-based time marching under neural network-controlled artificial viscosity}

\subsection{\label{FC-tstep}Spatial grid functions and spatio-temporal
  FC-based differentiation}

We consider in this work 1D problems on intervals
$I = [\xi_{\ell}, \xi_r]$ as well as 2D problems on open domains
$\Omega$ contained in rectangular regions $I\times J$, where
$I=[\xi_{\ell}, \xi_r]$ and $J=[\xi_{d}, \xi_u]$
($\xi_{\ell} < \xi_r$, $\xi_d < \xi_u$). Using a spatial meshsize $h$,
the spatially discrete vectors of unknowns and certain related flow
quantities will be represented by means of scalar and vector grid
functions defined on 1D or 2D discretization grids of the form
\[
  G = \{x_i\ :\ x_i = x_0 +ih,\quad i=0,\dots,N-1 \}\quad (x_0 = \xi_{\ell},\quad x_{N-1} = \xi_r),
\]
and
\[
  G = \overline{\Omega}\cap \big\{(x_i, y_j)\ :\ x_i = x_0 +ih, y_j =
  y_0 +jh, 0 \leq i \leq N_1-1,\ 0 \leq j \leq N_2 - 1\big\},
\]
respectively. Here $\overline{\Omega}$ denotes the closure of
$\Omega$, $x_0 = \xi_{\ell}$, $x_{N_1-1} = \xi_r$, $y_0 = \xi_d $ and
$y_{N_2-1} = \xi_u $.  In either case a function
\[
  b:G\to \mathbb{R}^q
\]
will be called a ``$q$-dimensional vector grid function''. Letting
\[
  \mathcal{I} = \big\{(i, j) \in \{0, \dots, N_1 - 1\} \times \{0, \dots, N_2-1\}\ :\ (x_i, y_j) \in G \big\},
\]
we will also write $b(x_i) = b_i$
($0\leq i\leq N-1$) and $b(x_i,y_j) = b_{ij}$
($(i, j) \in \mathcal{I}$). The set of $q$-dimensional vector grid
functions defined on $G$ will be denoted by $\mathcal{G}^q$.

It is important to mention that, although the two-dimensional setting
described above does not impose any restrictions on the character of
the domain $\Omega$, for simplicity, the FC-SDNN solver presented in
this paper assumes that the boundary of $\overline{\Omega}$ is given
by a union of horizontal and vertical straight segments, each one of
which runs along a Cartesian discretization line; see e.g. the Mach 3
forward-facing step case considered in Figure~\ref{M3
  solutions_a}. Extensions to general domains $\Omega$, which could
rely on either an
embedded-boundary~\cite{bruno2010high,lyon2010high,bruno2020two}
approach, or an overlapping patch boundary-conforming curvilinear
discretization
strategy~\cite{albin2011spectral,amlani2016fc,bruno2019higher}, is
left for future work.

A spatially-discrete but time-continuous version of the solution
vector $\textbf{\textit{e}}(\mathbf{x},t)$ considered in
Section~\ref{Preliminaries} for 1D problems (resp.  2D problems)
can be viewed as a time-dependent $q$-dimensional vector grid function
$\textbf{\textit{e}}_i=\textbf{\textit{e}}_i(t)$
(resp. $\textbf{\textit{e}}_{ij}=\textbf{\textit{e}}_{ij}(t)$).  Using
$\textbf{\textit{e}}_h=\textbf{\textit{e}}_h(t)$ to refer generically
to the $1D$ and $2D$ time-dependent grid functions
$\textbf{\textit{e}}_i$ and $\textbf{\textit{e}}_{ij}$, the
semidiscrete scheme for equation~\eqref{eq: convection diffusion eqn}
becomes
\begin{equation}
  \frac{d \textbf{\textit{e}}_h(t)}{dt} = L[\textbf{\textit{e}}_h(t)],
\end{equation} 
where $L$ denotes a consistent discrete approximation of the spatial
operator in~\eqref{eq: convection diffusion eqn}.

The discrete time evolution of the problem, on the other hand, is
produced, throughout this paper, by means of the 4-th order strong
stability preserving Runge-Kutta scheme
(SSPRK-4)~\cite{gottlieb2005high}---which, while not providing high
convergence orders for the non-smooth solutions considered in this
paper, does lead to low temporal dispersion and diffusion over smooth
space-time regions of the computational domain.  The corresponding
time step is selected adaptively at each time-step $t=t_n$ according to
the expression
\begin{equation} \label{eq: CFL} \Delta t = \frac{\textrm{CFL}}{\pi
    (\frac{ \max_{\textit{\textbf{x}} \in \Omega}\lvert
      S[\textit{\textbf{e}}](\mathbf{x}, t))) \rvert}{h} + \frac{
      \max_{\textit{\textbf{x}} \in \Omega} \mu[\textit{\textbf{e}}](\textit{\textbf{x}},
      t)}{h^2})}.
\end{equation}  
Here $\textrm{CFL}$ is a constant parameter that must be selected for
each problem considered (as illustrated by the various selections
utilized in Section~\ref{sec:numerical results}), and
$\mu[\textit{\textbf{e}}](\mathbf{x},t)$ and
$S = S[\textit{\textbf{e}}](\mathbf{x}, t)$ denote the artificial
viscosity (equations~\eqref{eq: V} and~\eqref{eq: V2D}) and a {\em
  maximum wave speed bound} (MWSB) operator (which must be
appropriately selected for each equation; see
Section~\ref{subsec:artificial_viscosity_strategy}) at the
spatio-temporal point $(\mathbf{x},t)$. (To avoid confusion we
emphasize that equation~\eqref{eq: CFL} utilizes the {\em maximum}
value for all $\mathbf{x}\in \Omega$ of the selected bound
$S[\textit{\textbf{e}}](\mathbf{x}, t)$ on the {\em maximum} wave
speed.)

To obtain FC-based approximate derivatives of a function
$F:K\to\mathbb{R}$ defined on a one-dimensional interval
$K=[x_0,x_{N-1}]$, whose values $(F_0,F_1,\dots,F_{N-1})^T$ are given
on the uniform mesh $\{x_0, x_1,\dots,x_{N-1}\}$, the interval $K$ is
re-scaled to $[0, 1]$ and the corresponding continuation function
$F^c$ is obtained by means of the FC-Gram procedure described in
Section~\ref{sec:fc}. The approximate derivatives at all mesh points
are then obtained by applying the IFFT algorithm to the Fourier
coefficients
\begin{equation}\label{fc_der}
    (\hat{F}^{c})'_k= \frac{2 \pi i k}{\beta}\hat{F}^c_k.
\end{equation}
of the derivative of the series~\eqref{FC_exp} and re-scaling back to
the interval $K$.

All of the numerical derivatives needed to evaluate the spatial
operator $L[\textbf{\textit{e}}_h(t)]$ are obtained via repeated
application of the 1D FC differentiation procedure described
above. For a function $F=F(x,y)$ defined on a two-dimensional domain
$\Omega$ and whose values $F_{ij}$ ($(i, j) \in \mathcal{I}$) are
given on a grid $G$ of the type described above in this section, for
example, partial derivatives with respect to $x$ along the line
$y=y_{j_0}$ for a relevant value of $j_0$ are obtained by
differentiation of the FC expansion obtained for the function values
$(F(x_i,y_{j_0}))_i$ for integers $i$ such that
$(i, j_0) \in \mathcal{I}$. The $y$ differentiation process proceeds
similarly. Mixed derivatives, finally, are produced by successive
application of the $x$ and $y$ differentiation processes. Details
concerning the filtered derivatives used in the proposed scheme are
provided in Section~\ref{filtering}.

The boundary conditions of Dirichlet and Neumann considered in this
paper are imposed as part of the differentiation process described
above. Dirichlet boundary conditions at time
$t_{n,\nu}$ ($t_n<t_{n,\nu}\leq t_{n+1}$) corresponding to the $\nu$-th
SSPRK-4 stage ($\nu = 1, \dots, 4$) for the time-step starting at
$t=t_n$, are simply imposed by overwriting the boundary values of the
unknown solution vector $\textit{\textbf{e}}_h$ obtained at time
$t=t_{n,\nu}$ with the given boundary values at that time,  prior to the
evaluation of the spatial derivatives needed for the subsequent
SSPRK-4 stage.  Neumann boundary conditions are similarly enforced by
constructing appropriate continuation vectors (Section~\ref{sec:fc})
after each stage of the SSPRK-4 scheme on the basis of
the modified pre-computed matrix $\widetilde Q$ mentioned
in Section~\ref{sec:fc}.

It is known that enforcement of the given physical boundary conditions
at intermediate Runge-Kutta stages, which is referred to as the
``conventional method'' in~\cite{carpenter1995theoretical}, may lead
to a reduced temporal order of accuracy at spatial points in a
neighborhood of the boundary of the domain boundary. This is not a
significant concern in the context of this paper, where the global
order of accuracy is limited in view of the discontinuous character of
the solutions considered. Alternative approaches that preserve the
order of accuracy for smooth solutions, such as those introduced
in~\cite{carpenter1995theoretical,pathria1997correct}, could also be
used in conjunction with the proposed approach. Another alternative,
under which no boundary conditions are enforced at intermediate
Runge-Kutta stages~\cite{kopriva2009implementing}, can also be
utilized in our context, but we have found the conventional method
leads to smoother solutions near boundaries.

\subsection{Neural network-induced
  smoothness-classification}\label{subsec: nn smoothness
  classification}

\subsubsection{Smoothness-classification operator and data
  pre-processing}\label{subsec: smoothness classifier}

The method described in the forthcoming
Section~\ref{subsec:artificial_viscosity_strategy} for determination
of the artificial viscosity values
$\mu[\textit{\textbf{e}}](\mathbf{x},t)$ (cf. also
Section~\ref{art_vis}) relies on the ``degree of smoothness'' 
of a certain function $\Phi(\textit{\textbf{e}})(\mathbf{x},t)$
(called the ``proxy variable'') of the unknown solution vector
$\textit{\textbf{e}}$.  In detail,
following~\cite{schwander2021controlling}, in this paper a proxy
variable $\Phi(\textit{\textbf{e}})$ is used, which equals the
velocity $u$, $\Phi(\textit{\textbf{e}}) = u$, (resp. the Mach number,
$\Phi(\textit{\textbf{e}}) = \norm {\mathbf{u}}
\sqrt{\frac{\rho}{\gamma p}}$) for equations~\eqref{eq: linear
  advection} through~\eqref{eq: burgers 2d}
(resp. equations~\eqref{eq: euler 1d equation} and~\eqref{eq: euler 2d
  equation}). The degree of smoothness of the function
$\Phi(\textit{\textbf{e}})$ at a certain time $t$ is characterized by
a smoothness-classification operator
$\tau = \tau[\Phi(\textit{\textbf{e}})]$ that analyzes the
oscillations in an FC expansion of $\Phi(\textit{\textbf{e}})$---which
is itself obtained from the discrete numerical values
$\bm{\phi}=\Phi(\textbf{\textit{e}}_h)$, so that, in particular,
$\tau[\Phi(\textit{\textbf{e}})] = \tilde\tau[\bm{\phi}]$ for some
discrete operator $\tilde\tau$. The determination (or, rather,
estimation) of the degree of smoothness by the operator $\tilde\tau$
is effected on the basis of an Artificial Neural Network (ANN). (We
introduce the operators $\tau$ and $\tilde\tau$ for the specific
function $\Phi(\textbf{\textit{e}})$, but, clearly, the algorithm
applies to arbitrary scalar or vector functions, as can be seen e.g. in
the application of these operators, in the context of network
training, in Section~\ref{subsubsec: nn archirecture and training}.)

We first describe the operator $\tilde\tau = \tilde\tau[\bm{\phi}]$
for for a conservation law over a one-dimensional interval
$I=[\xi_\ell,\xi_r]$ discretized by an $N$-point equispaced mesh
$(x_0,\dots,x_{N-1})$ of mesh-size $h$, and for which FC expansions
are obtained on the basis of the extended equispaced mesh
$\{x_0,\dots,x_{N+C-1}\}$. (Note that, in accordance with
Section~\ref{sec:fc}, this extended mesh includes the discrete points
$\{x_0,\dots,x_{N-1}\}$ in the interval $I$ as well as the discrete
points $\{x_{N},\dots,x_{N+C-1}\}$ in the FC extension region.) In
this case, the evaluation of the operator $\tilde \tau$ proceeds as
follows. 
\begin{enumerate}[(i)]
\item\label{pt1} Obtain the FC expansion coefficients
  $(\hat{\bm{\phi}}^c_{-M}, \dots, \hat{\bm{\phi}}^c_{M})^T$ of
  $\Phi(\textbf{\textit{e}})$ by applying the FC procedure described
  in Section~\ref{sec:fc} to the column vector
  $(\bm{\phi}_0, \dots, \bm{\phi}_{N-1})^T$. (Note that in the present
  1D case we have $\bm{\phi}_j = \Phi(\textbf{\textit{e}}_j)$.)
\item\label{pt2} For a suitable selected non-negative number
  $\delta < h$, evaluate the values $\bm{\phi}^{(\delta)}_{j}$
  ($0\leq j\leq N+C-1$ ) of the FC expansion obtained in
  point~\eqref{pt1} at the shifted grid points
  $x_0+\delta ,x_1+\delta ,\dots,x_{N+C-1}+\delta$. This is achieved
  by applying the FFT algorithm to the ``shifted'' Fourier
  coefficients
  $\hat{\bm{\phi}}^\delta = (\hat{\bm{\phi}}^\delta_{-M}, \dots,
  \hat{\bm{\phi}}^\delta_{M})$ where
  $\hat{\bm{\phi}}^\delta_j = \hat{\bm{\phi}}^c_{j} \exp(\frac{2 \pi i
    j \delta}{\beta})$. Here, as in Section~\ref{sec:fc} and
  equation~\eqref{FC_exp}, $\beta$ denotes the length of the FC
  periodicity interval. Throughout this work, the value
  $\delta = \frac{h}{10}$ is used for classification of flow
  discontinuities. As indicated in Section~\ref{subsubsec: nn
    archirecture and training}, different values of $\delta$ are used
  in the training process.
\item\label{pt3} For each $j \in \{0, \dots, N-1 \}$, form the seven-point
  stencil
  \[
    \bm{\phi}^{(\delta,j)}= \big(\bm{\phi}^{(\delta)}_{m(j-3,N+C)},
  \dots, \bm{\phi}^{(\delta)}_{m(j,N+C)}, \dots,
  \bm{\phi}^{(\delta)}_{m(j+3,N+C)} \big)^T
\]
of values of the shifted grid function obtained per
point~\eqref{pt2}. (Here, for an integer $0 \leq j \leq P $, $m(j,P)$
denotes the remainder of $j$ modulo $P$, that is to say, $m(j,P)$ is
the only integer between $0$ and $P-1$ such that $j - m(j,P)$ is an
integer multiple of $P$. In view of the extended domain inherent in
the continuation method, use of the remainder function $m$ allows for
the smoothness classification algorithm to continue to operate
correctly even at points $x_j$ near physical boundaries---for which
the seven-point subgrid $(x_{j-3},\dots,x_j,\dots,x_{j+3})$ may not be
fully contained within the physical domain.)
\item\label{pt4} Obtain the modified stencils
  $\widetilde{\bm{\phi}}^{(\delta,j)} =
  \big(\widetilde{\bm{\phi}}^{(\delta)}_{m(j-3, N+C)}, \ldots,
  \widetilde{\bm{\phi}}^{(\delta)}_{m(j, N+C)}, \ldots
  \widetilde{\bm{\phi}}^{(\delta)}_{m(j + 3, N+C)}\big)^T$ given by
\begin{equation}
    \widetilde{\bm{\phi}}^{(\delta)}_{m(j + r, N+C)} = \bm{\phi}^{(\delta)}_{m(j + r, N+C)} - \ell_{j + r}\quad (-3 \leq r \leq 3),
\end{equation}
  that result by subtracting the ``straight line''
\begin{equation}
  \ell_{j + r} = \bm{\phi}^{(\delta)}_{m(j - 3, N+C)}
  + \frac{r+3}{6}(\bm{\phi}^{(\delta)}_{m(j + 3, N+C)} -
  \bm{\phi}^{(\delta)}_{m(j - 3,  N+C)})\quad (-3 \leq r \leq 3)
\end{equation}
passing through the first and last stencil points.
\item\label{pt5} Rescale each stencil $\widetilde{\bm{\phi}}^{(\delta,j)}$
  so as to obtain the ANN input stencils 
  $$ \check{\bm{\phi}}^{(\delta,j)} =
  \big(\check{\bm{\phi}}^{(\delta)}_{m(j-3, N+C)}, \ldots,
  \check{\bm{\phi}}^{(\delta)}_{m(j, N+C)}, \ldots
  \check{\bm{\phi}}^{(\delta)}_{m(j+3, N+C)}\big)^T,$$ given by
\begin{equation}\label{norm_st_dat}
  \check{\bm{\phi}}^{(\delta)}_{m(j + r, N+C)} = \frac{2 \widetilde{\bm{\phi}}^{(\delta)}_{m(j + r, N+C)} - M^{(+)}_j - M^{(-)}_j}{M^{(+)}_j - M^{(-)}_j} \quad (-3 \leq r \leq 3)
\end{equation}
where
\begin{equation}
  M^{(+)}_j  = \max_{-3 \leq r \leq 3} \widetilde{\bm{\phi}}^{(\delta)}_{m(j + r, N+C)}
  \quad\mbox{and}\quad    M^{(-)}_j  = \min_{-3 \leq r \leq 3} \widetilde{\bm{\phi}}^{(\delta)}_{m(j + r, N+C)}.
  \end{equation}
  Clearly, the new stencil entries satisfy
  satisfy
  $-1\leq \check{\bm{\phi}}^{(\delta)}_{m(j + r, N+C)}\leq 1$.

\item \label{pt6} Apply the ANN algorithm described in
  Section~\ref{subsubsec: nn archirecture and training} to each one of
  the stencils
  $ \check{\bm{\phi}}^{(\delta,j)} =
  \big\{\check{\bm{\phi}}^{(\delta)}_{m(j-3, N+C)}, \ldots,
  \check{\bm{\phi}}^{(\delta)}_{m(j, N+C)}, \ldots
  \check{\bm{\phi}}^{(\delta)}_{m(j+3, N+C)}\big\}$, to produce a
  four-dimensional vector $w^j$ of estimated probabilities (EP)
  for each $j \in [0, \dots, N-1]$,
  where $w^j_1$ is the EP that $\Phi(\textit{\textbf{e}})$ is
  discontinuous on the subinterval $I_j = [x_j - 3h, x_j + 3h]$,
  where, for $i = 2,3$, $w^j_i$ equals the EP that
  $\Phi(\textit{\textbf{e}})\in\mathcal{C}^{i-2}\setminus\mathcal{C}^{i-1}
  $ on $I_j$, and where, for $i=4$, $w^j_i$ equals the EP that
  $\Phi(\textit{\textbf{e}})\in\mathcal{C}^{2}$ on $I_j$. Define
  $\tau[\bm{\phi}]_j$ as the index $i$ corresponding to the maximum
  entry of $w^j_i$ ($i=1,\dots 4$):
  \begin{equation}\label{tau}
    \tilde
    \tau[\bm{\phi}]_j = \argmax_{1 \leq i \leq 4}(w^j_i)\quad (0 \leq j \leq N-1).
  \end{equation}
  (Note that, for points $x_j$ close to physical boundaries, the
  interval $I_j = [x_j - 3h, x_j + 3h]$, within which the smoothness
  of the function $\Phi(\textit{\textbf{e}})$ is estimated, can extend
  beyond the physical domain and into the extended Fourier
  Continuation region; cf. also point~\ref{pt3} above.)
\end{enumerate}
This completes the definition of the 1D smoothness classification
operator $\tilde\tau$.

For 2D configurations, in turn, we define a two-dimensional smoothness
classification operator
$\tau_{xy}[\Phi(\textit{\textbf{e}})] = \tilde\tau_{xy}[\bm{\phi}]$,
similar to the 1D operator, which classifies the smoothness of the
proxy variable $\Phi(\textit{\textbf{e}})$ on the basis of its
discrete values $\bm{\phi} = \bm{\phi}_{ij}$. Note the $xy$ subindex
which indicates 2D classification operators $\tau_{xy}$ and
$\tilde\tau_{xy}$; certain associated 1D ``partial'' discrete
classification operators in the $x$ and $y$ variables, which are used
in the definition of $\tilde\tau_{xy}$, will be denoted by
$\tilde{\tau}_x$ and $\tilde{\tau}_y$, respectively.

In order to introduce the operator $\tilde\tau_{xy}$ we utilize
certain 1D sections of both the set $\mathcal{I}$ and the grid
function $\bm{\phi} = \bm{\phi}_{ij}$ (see
Section~\ref{FC-tstep}). Thus, the $i$-th horizontal section
(resp. the $j$ vertical section) of $\mathcal{I}$ is defined by
$\mathcal{I}_{i:} = \{j\in \mathbb{Z}\ :\ (i,j)\in\mathcal{I}\}$
(resp.
$\mathcal{I}_{:j} = \{i\in \mathbb{Z}\ :\
(i,j)\in\mathcal{I}\}$). Similarly, for a given 2D grid function
$\bm{\phi} = \bm{\phi}_{ij}$, the $i$-th horizontal section
$\bm{\phi}_{i:}$ (resp. the $j$ vertical section $\bm{\phi}_{:j}$) of
$\bm{\phi}$ is defined by
$\left(\bm{\phi}_{i:}\right)_j = \bm{\phi}_{ij}$,
$j\in \mathcal{I}_{i:}$ (resp.
$\left(\bm{\phi}_{:j}\right)_i = \bm{\phi}_{ij}$,
$i\in \mathcal{I}_{:j}$). Utilizing these notations we define
\begin{equation}\label{tauxy}
  \tilde\tau_{xy}[\mathbf{\bm{\phi}}]_{ij} =
  \min\left\{\tilde\tau_x[\bm{\phi}_{i:}]_{j},
    \tilde\tau_y[\bm{\phi}_{:j}]_{i}\right\},
\end{equation}
where, as suggested above, $\tilde\tau_x$ (resp. $\tilde\tau_y$)
denotes the discrete one-dimensional classification operator along the
$x$ direction (resp. the $y$ direction), given by~\eqref{tau} but with
$j\in \mathcal{I}_{i:}$ (resp. $i\in \mathcal{I}_{:j}$). In other
words, the 2D smoothness operator $\tilde\tau_{xy}$ equals the lowest
degree of smoothness between the classifications given by the two
partial classification operators.

\begin{remark}\label{noise}
  Small amplitude noise in the proxy variable can affect ANN analysis,
  leading to misclassification of stencils and under-prediction of the
  smoothness of the proxy variable. In order to eliminate the effect
  of noise, stencils $\check{\bm{\phi}}^{(\delta,j)}_j$ for which
  $M^{(+)}_j - M^{(-)}_j \leq \varepsilon$, for a prescribed value of
  $\varepsilon$, are assigned regularity
  $\tilde\tau[\bm{\phi}]_j = 4$. Throughout this paper we have used
  the value $\varepsilon = 0.01$.
\end{remark}

\subsubsection{Neural network architecture and training}\label{subsubsec:
nn archirecture and training}
The proposed strategy relies on standard neural-network techniques and
nomenclature~\cite[Sec. 6]{goodfellow2016}: it utilizes an ANN with a
depth of four layers, including three fully-connected hidden layers of
sixteen neurons each. The ANN takes as input a seven-point
``preprocessed stencil''
$z = (z_1,z_2,z_3,z_4,z_5, z_6, z_7)^T$---namely, a stencil $z$ that
results from an application of points~\eqref{pt1} to~\eqref{pt5} in
the previous section to the 401-coordinate vector $\mathbf{F}$ of grid
values obtained for a given function $F$ on a 401-point equispaced
grid in the interval $[0,2\pi]$---in place of the grid values of the
proxy variable $\Phi(\textbf{\textit{e}})$---resulting in a total of
$401$ stencils, one centered around each one of the $401$ grid points
considered; cf. points~\eqref{pt1} to~\eqref{pt3} and note that, on the
basis of the FC-extended function, the stencils near endpoints draw
values at grid points outside the interval $[0,2\pi]$. (A variation of
point~\eqref{pt2} is used in the training process: shift values
$\delta = \frac{h}{10},\frac{2h}{10}\dots\frac{10h}{10}$ are used to
produce a variety of seven-point stencils for {\em training} purposes
instead of the single value $\delta = \frac{h}{10}$ used while
employing the ANN in the {\em classification} process.) The output of
the final layer of the ANN is a four-dimensional vector
$\widetilde w = (\widetilde w_1,\widetilde w_2,\widetilde
w_3,\widetilde w_4)^T$, from which, via an application of the softmax
activation function~\cite[Sec. 4.1]{goodfellow2016}, the EP mentioned
in point~\eqref{pt6} of the previous section, are obtained:
\begin{equation}\label{sft_mx}
  w_{i}=\frac{e^{\widetilde w_{i}}}{\sum_{\ell=1}^{4} e^{\widetilde w_{\ell}}},\quad 1\leq i\leq 4.
\end{equation}
(The values $w_i^j$ ($1\leq i\leq 4$) mentioned in point~\eqref{pt6}
result from the expression~\eqref{sft_mx} when the overall scheme
described above in the present Section~\ref{subsubsec: nn archirecture
  and training} is applied to $ z = \check{\bm{\phi}}^{(\delta,j)}$.)
The \textrm{ELU} activation function
\begin{equation}
   \textrm{ELU}(x; \alpha_0) = \left\lbrace
    \begin{array}{ccc}
        x & \mbox{if } x > 0\\
        \alpha_0 (e^x - 1) & \mbox{if } x \leq 0,
    \end{array}\right. \\   
\end{equation}
with $\alpha_0 = 1$, is used in all of the hidden layers.

In what follows we consider, for both the ANN training and validation
processes, the data set $\mathcal{D}_\mathcal{F}$ of preprocessed
stencils resulting from the set
$\mathcal{F}=\big\{\big(F_k,D_k\big), k=1,2,\dots \big\}$ of all pairs
$\big(F_k,D_k\big)$, where $F_k$ is a function defined on the interval
$[0, 2\pi]$ and where $D_k$ is a certain ``restriction domain'', as
described in what follows. The functions $F_k$ are all the functions
obtained on the basis of one of the five different parameter-dependent
analytic expressions
\begin{equation}\label{subeq: data set}
\begin{split}
    f_1(x) &= \sin(2 a x) \\
    f_2(x) &= a |x - \pi| \\
    f_3(x) &= \left\lbrace
        \begin{array}{ccc}
            a_1 & \mbox{if} \quad |x - \pi| \leq a_3\\
            a_2 & \mbox{if} \quad |x - \pi| > a_3
        \end{array}\right. \\
    f_4(x) &= \left\lbrace
        \begin{array}{ccc}
            a_1|x - \pi| - a_1 a_3  & \mbox{if} \quad  |x - \pi|\leq a_3\\
            a_2|x - \pi| - a_2 a_3  & \mbox{if} \quad |x - \pi| > a_3
        \end{array}\right. \\
    f_5(x) &= \left\lbrace
        \begin{array}{ccc}
            0.5a_1|x - \pi|^2 - a_1 a_3  & \mbox{if} |x - \pi|\leq a_3\\
            a_2|x - \pi|^2 - a_2  - 0.5a_3^2(a_1 - a_2)  & \mbox{if} |x - \pi| > a_3
        \end{array}\right. \\
\end{split}
\end{equation}
proposed in~\cite{schwander2021controlling}, for each one of the
possible selections of the parameters $a,a_1,a_2,a_3$, as prescribed
in Table~\ref{training_set_parameters}. The corresponding parameter
dependent restriction domains $D_k$ are also prescribed in
Table~\ref{training_set_parameters}; in all cases $D_k$ is a
subinterval of $[0, 2\pi]$.

The restriction domains $D_k$ are used to constrain the choice of
stencils to be used among all of the $401$ stencils available for each
function $F_k$---so that, for a given function $F_k$, the preprocessed
stencils associated with gridpoints contained within $D_k$, but not
others, are included within the set $\mathcal{D}_{\mathcal{F}}$.  The
set $\mathcal{D}_{\mathcal{F}}$ is randomly partitioned into a
training set $\mathcal{D}_{\mathcal{F}}^{\mathcal{T}}$ containing
$80\%$ of the elements in $\mathcal{D}_{\mathcal{F}}$ (which is used
for optimization of the ANN weights and biases), and a validation set
$\mathcal{D}_{\mathcal{F}}^{\mathcal{V}}$ containing the remaining
$20\%$---which is used to evaluate the accuracy of the ANN after each
epoch~\cite[Sec. 7.7]{goodfellow2016}.

The network training and validation processes rely on use of a
``label'' function $C$ defined on $\mathcal{D}_{\mathcal{F}}$ which
takes one of four possible values. Thus each stencil
$z \in \mathcal{D}_{\mathcal{F}}$ is labeled by a class vector
$C(z) = (C_1(z),C_2(z),C_3(z),C_4(z))^T$, where for each $z$,
$C(z) = (1, 0, 0, 0)^T$, $(0, 1, 0, 0)^T$, $(0, 0, 1, 0)^T$ or
$(0, 0, 0, 1)^T$ depending on whether $z$ was obtained from a function
$F_k$ that is $\mathcal{C}^{2}$,
$\mathcal{C}^{1}\setminus\mathcal{C}^{2}$,
$\mathcal{C}^{0}\setminus\mathcal{C}^{1}$, or discontinuous over the
subinterval $I_z\cap [0,2\pi]$, where $I_z$ denotes the interval
spanned by the set of seven consecutive grid points associated with
the prepocessed stencil $z$.

The ANN is characterized by a relatively large number of parameters
contained in four weight matrices of various dimensions (a
$16\times 7$ matrix, a $4\times 16$, and two $16\times 16$ matrices),
as well as four bias vectors (one 4-dimensional vector and three
16-dimensional vectors). In what follows a single parameter vector $X$
is utilized which contains all of the elements in these matrices and
vectors in some arbitrarily prescribed order. Utilizing the parameter
vector $X$, for each stencil $z$ the ANN produces the estimates $w_i$,
given by~\eqref{sft_mx}, of the actual classification vector
$C(z)$. In order to account for the dependence of $w_i$ on the
parameter vector $X$ for each stencil $z$, in what follows we write
$w_i = A_i(X,z)$ ($1\leq i \leq 4$).

The parameter vector $X$ itself is obtained by training the network on
the basis of existing data, which is accomplished in the present
context by selecting $X$ as an approximate minimizer of the ``cross
entropy'' loss function~\cite[Sec. 6.2]{goodfellow2016}
\begin{equation}\label{entropy}
  \mathcal{L}(X) = -\frac{1}{N_{\mathcal{T}}} \sum_{z \in \mathcal{D}_{\mathcal{F}}^{\mathcal{T}}}\sum_{i = 1}^4C_i(z)\log(A_i(X, z))
\end{equation}
over all $z$ in the training set
$\mathcal{D}_{\mathcal{F}}^{\mathcal{T}}$, where $N_{\mathcal{T}}$
denotes the number of elements in the training set. The loss function
$\mathcal{L}$ provides an indicator of the discrepancy between the EP
$A(X, z) = (A_1(X, z),A_2(X, z),A_3(X, z),A_4(X, z) )$ produced by the
ANN and the corresponding classification vector
$C(z) = (C_1(z),C_2(z),C_3(z),C_4(z))^T$ introduced above, over all
the preprocessed stencils $z \in \mathcal{D}_{\mathcal{F}}$.  The
minimizing vector $X$ of weights and biases define the network, which
can subsequently be used to produce $A(X, z)$ for any given
preprocessed stencil $z$.

The Neural Network is trained (that is, the loss function $L$ is
minimized with respect to $X$) by exploiting the stochastic gradient
descent algorithm without momentum~\cite[Secs. 8.1,
8.4]{goodfellow2016}, with mini batches of size 128 and with a
constant learning rate of $10^{-6}$. The weight matrices and bias
vectors are initialized using the Glorot
initialization~\cite{glorot2010understanding}. The training set is
randomly re-shuffled after every epoch, and the validation data is
re-shuffled before each network validation. The best performing
network obtained, which is used for all the illustrations presented in
this paper, has a training accuracy of $99.61\%$ and validation
accuracy of $99.58\%$.

\begin{table}
\centering
\begin{tabular}{|c|c|c|c|}
  \hline
  $f(x)$     & Parameters   & restriction domains    & $\tau$  \\ \hline
  $f_1$    & $a = \big\{ -20, -19.5, \ldots, 19,5 \big\}$      & $\big[ 0, 2 \pi\big]$            & 4     \\ \hline
  $f_2$    & $a = \big\{ -10, -9, \ldots, 10 \big\}$   & $\big[3.53, 5.89\big]$            & 4     \\ \hline
  $f_3$    & \begin{tabular}{c} $a_1 =  \big\{ -10, -9, \ldots, 9 \big\}$ \\ $a_2 =  \big\{ -10, -9, \ldots, 9 \big\}$ \\ $a_3 =  \big\{0.25, 0.5, \ldots, 2.5 \big\}$ \\ s.t.  $a_1 \neq a_2$ \end{tabular}    & $\big[\pi + a_3 - 0.05, \pi + a_3 + 0.05\big]$            & 1     \\ \hline
  $f_4$   & \begin{tabular}{c} $a_1 =  \big\{ -10, -9, \ldots, 9 \big\}$ \\ $a_2 =  \big\{ -10, -9, \ldots, 9 \big\}$ \\ $a_3 =  \big\{0.25, 0.5, \ldots, 2.5 \big\}$ \\ s.t. $a_1 > 2a_2$ or $a_1 < 0.5a_2$ \end{tabular}   & $ \big[\pi + a_3 - 0.05, \pi + a_3 + 0.05\big]$             & 2     \\ \hline
  $f_5$    &  \begin{tabular}{c} $a_1 =  \big\{ -10, -9, \ldots, 9 \big\}$ \\ $a_2 =  \big\{ -10, -9, \ldots, 9 \big\}$ \\ $a_3 =  \big\{0.25, 0.5, \ldots, 2.5 \big\}$ \\ s.t. $a_1 > 5a_2$ or $a_1 < 0.2a_2$ \end{tabular}     & $ \big[\pi + a_3 - 0.05, \pi + a_3 + 0.05\big]$            & 3     \\ \hline
\end{tabular}
\caption{Data set.}
\label{training_set_parameters}
\end{table}

\subsection{SDNN-localized artificial viscosity
  algorithm} \label{subsec:artificial_viscosity_strategy} As indicated
in Section~\ref{sec:introduction}, in order to avoid introduction of
spurious irregularities in the flow field, the algorithm proposed in
this paper relies on use of smoothly varying artificial viscosity
assignments. For a given discrete solution vector
$\textit{\textbf{e}}_h$, the necessary grid values of the artificial
viscosity, which correspond to discrete values of the continuous
operator $\mu = \mu[\textit{\textbf{e}}]$ in~\eqref{f_visc}, are
provided by a certain discrete viscosity operator
$\tilde\mu = \tilde\mu[\textit{\textbf{e}}_h]$. The discrete operator
$\tilde\mu$ is defined in terms of a number of flow- and
geometry-related concepts, namely the proxy variable $\bm{\phi}$
defined in Section~\ref{subsec: smoothness classifier} and the
smoothness-classification operator given by equations~\eqref{tau}
and~\eqref{tauxy} for the 1D and 2D cases, respectively, as well as
certain additional functions and operators, namely a ``weight
function'' $R$ and ``weight operator'' $\widetilde R$, an MWSB operator
$S$ (see Section~\ref{FC-tstep}) and its discrete version $\widetilde S$,
a sequence of ``localization stencils'' (denoted by $L^i$ with
$0\leq i\leq N-1$ in the 1D case, and by $L^{i, j}$ with
$(i, j) \in \mathcal{I}$ in the 2D case) , and a
''windowed-localization'' operator $\Lambda$. A detailed description
of the 1D and 2D discrete artificial viscosity operators
$\tilde\mu = \tilde\mu[\textit{\textbf{e}}_h]$ is provided in
Sections~\ref{subsubsec:1d_case} and~\ref{subsubsec: 2d case},
respectively.

\subsubsection{One-dimensional case}
\label{subsubsec:1d_case}

The proxy variable $\bm{\phi}$ and 1D smoothness-classification
operator $\tau$ that are used in the definition of the 1D artificial
viscosity operator have been described earlier in this paper; in what
follows we introduce the additional necessary functions and operators
mentioned above.

The weight function $R$ assigns a viscosity weight according to the
smoothness classification; throughout this paper we use the weight
function given by $R(1) = 2$, $R(2)=1$, $R(3) = 0$, and $R(4)=0$; the
corresponding grid-function operator $\widetilde R$, which acts over the
set of grid functions $\eta$ with grid values $1$, $2$, $3$ and $4$,
is defined by $\widetilde R[\eta]_i = R(\eta_i)$.

The MWSB operator $S : \mathcal{G}^{q} \to \mathcal{G}$ maps the
$q$-dimensional vector grid function $\textit{\textbf{e}}_h$ onto a
grid function corresponding to a bound on the maximum eigenvalue of
the flux Jacobian
$\left(J_\textbf{\textit{e}}f\right)_{k\ell}
=\left(\partial_{\textbf{\textit{e}}^\ell}f_k\right)$ at
$\textbf{\textit{e}}= \textit{\textbf{e}}_h$, where
$\textbf{\textit{e}}^\ell$, (resp. $f_k$) denotes the $\ell$-th
(resp. $k$-th) component of the unknowns solution vector
$\textbf{\textit{e}}$ (resp. of the convective flux
$f(\textbf{\textit{e}})$). For the one-dimensional problems, the MWSB
operator $S(\textit{\textbf{e}})$ (resp. the discretized operator
$\widetilde S[\textit{\textbf{e}}_h]$ on the grid $\{x_i\}$) is taken to
{\em equal} the maximum characteristic speed (since the maximum
characteristic speed is easily computable from the velocity in the 1D
case), so that $S(\textit{\textbf{e}}) = a$ (resp
$\widetilde S[\textit{\textbf{e}}_h]_i = a_i$) for the 1D Linear Advection
equation~\eqref{eq: linear advection}, $S(\textit{\textbf{e}}) = |u|$
(resp $\widetilde S[\textit{\textbf{e}}_h]_i = |u_i|$) for the 1D Burgers
equation~\eqref{eq: burgers 1d}, and
$S(\textit{\textbf{e}}) = |u| +a$~\cite{leveque1992numerical} (resp
$\widetilde S[\textit{\textbf{e}}_h]_i = |u_i| + a_i$) in the case of the
1D Euler problem~\eqref{eq: euler 1d equation} (in terms of the sound
speed~\eqref{snd_sp}).

The localization stencil $L^i$ ($0 \leq i \leq N-1$) is a set of seven
points that surround $x_i$:
$L^i = \big\{x_{i-3}, \ldots, x_{i}, \ldots x_{i+3}\big\}$ for
$4 \leq i \leq N-4$, $L^i = \big\{x_{0}, \ldots x_{6}\big\}$ for
$i \leq 3$, and $L^i = \big\{x_{N-7}, \ldots x_{N-1}\big\}$ for
$i \geq N-3$.

\begin{figure}
  \centering
  \includegraphics[width=0.3\textwidth]{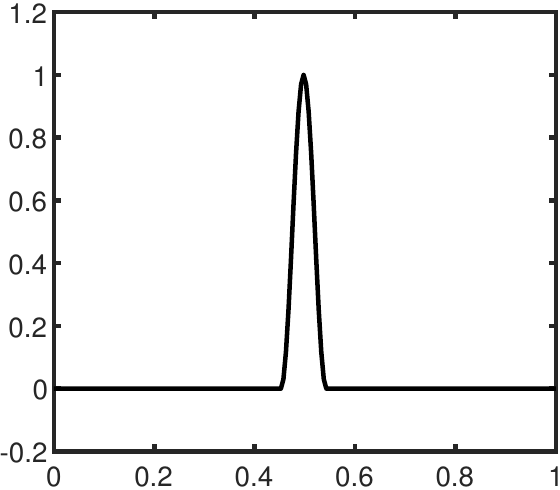}
  \hspace{0.5 cm}
  \includegraphics[width=0.3\textwidth]{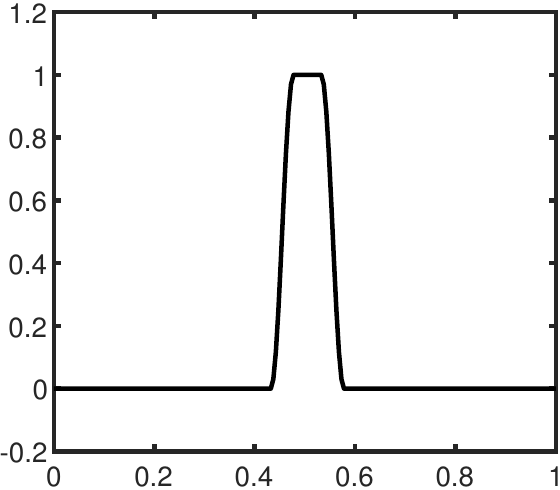}
  \caption{Left: Windowing function $W^{j}(x) $ (with $j=100$, on the
    domain $[0, 1]$, and using $N = 200$ discretization points)
    utilized in the definition of the windowed-localization operator
    $\Lambda$. Right: Windowing function $q_{18, 9}(x - z)$ utilized
    for the localized filtering of the initial condition on the domain
    $[0, 1]$ (using $N = 200$ discretization points) for a
    discontinuity located at $z = 0.5$.}
         \label{fig:visc_filt_window}    
\end{figure}

The windowed-localization operator $\Lambda$ is
constructed on the basis of the window function
\begin{equation} \label{eq: qf}
q_{c, r}(x)=
\left\lbrace
  \begin{array}{ccl}
    1 & \mbox{if} &  |x| <  ch/2\\
    \cos^2\Big(\frac{\pi (|x| - ch/2)}{rh}\Big) & \mbox{if} &    ch/2 \leq |x| \leq  (c/2 + r)h \\
    0 & \mbox{if} & |x| > (c/2 + r)h, \\
  \end{array}\right.
\end{equation}
depicted in Figure~\ref{fig:visc_filt_window} left, where $c$ and $r$
denote small positive integer values, with $c$ even. (Note that the
$q_{c, r}$ notation does not explicitly display the $h$-dependence of
this function.) Using the window function $q_{c, r}$, two sequences of
windowing functions, denoted by $W^j$ and $\check{W}^j$
($0 \leq j \leq N-1$), are defined, where the second sequence is a
normalized version of the former. In detail $W^j(x)$ is obtained by
translation of the function $q_{c, r}$ with $c=0$ and $r = 9$:
$W^j(x) = q_{0, 9}(x - x_j)$; the corresponding grid values of this
function on the grid $\{x_i\}$ are denoted by $ W^j_i =
W^j(x_i)$.
The normalized windowing functions $\check{W}^j$ and the
windowed-localization operator $\Lambda$, finally, are given by
\begin{equation} \label{eq: W_tilde}
\check{W}^j_i = \frac{W^j_i}{\sum_{k=0}^{N-1}W^j_k},
\end{equation}
and
\begin{equation} \label{eq:Lambda} \Lambda[b]_i =
  \sum_{k=0}^{N-1}\check{W}^k_i b_k,
\end{equation}
respectively.  Using these operators and functions, we define the 1D
artificial viscosity operator
\begin{equation} \label{eq: V}
\tilde\mu[\textit{\textbf{e}}_h]_i = \Lambda[\widetilde{R}(\tilde\tau[\bm{\phi}])]_i \cdot \max_{j\in L^i}(\widetilde S[\textit{\textbf{e}}_h]_j) h;
\end{equation}
as mentioned in Section~\ref{sec:introduction} and demonstrated in
Section~\ref{sec:numerical results}, use of the smooth artificial
viscosity assignments produced by this expression yield smooth flows
away from shocks and other discontinuities.

\begin{figure}
  \centering
         \includegraphics[width=0.4\textwidth]{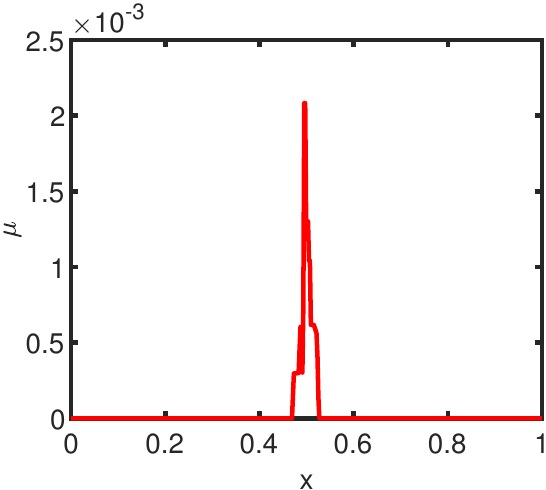}
    \hspace{0.5 cm}
    \includegraphics[width=0.4\textwidth]{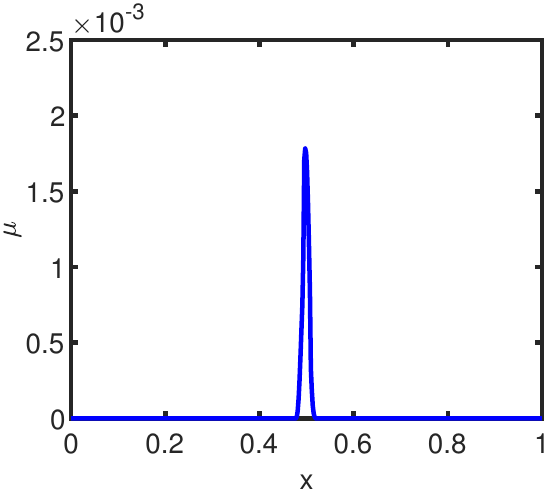}
    \caption{Comparison of the viscosity functions arising at the
      first viscous time-step for the Sod problem, using $N=500$
      discretization points. Left: viscosity used
      in~\cite{schwander2021controlling}. Right: viscosity used in the
      present FC-SDNN method (equation~\eqref{eq: V}).}
\label{fig:viscosity_comparisons}       
\end{figure}

It is important to note the essential role of the
windowed-localization operator in the assignment of smooth viscosity
profiles. The smooth character of the resulting viscosity functions is
illustrated in Figure~\ref{fig:viscosity_comparisons}, which showcases
the viscosity assignments corresponding to the second time-step in the
solution process. (This run corresponds to the Sod problem described
in Section~\ref{Euler 1D results}.) The left image displays the
viscosity profiles used in~\cite{schwander2021controlling} (which do
not utilize the smoothing windows~\eqref{eq:Lambda}) and the right
image presents the window-based viscosity profile~\eqref{eq: V}. The
right-hand profile, which is comparable in size but, in fact, more
sharply focused around the shock than the non-smooth profile on the
left-hand image, helps eliminate spurious oscillations that otherwise
arise from viscosity non-smoothness, and allows the FC-SDNN method to
produce smooth flow fields, as demonstrated in Section~\ref{Euler 1D
  results}.

\begin{remark}\label{periodic}
  For the case of a 1D periodic problem, such as those considered in
  Section~\ref{dispersionless}, the localization stencils and the
  windowing functions are defined by
  $L^i = \big\{x_{m(i-3,N)}, \ldots, x_{m(i, N)}, \ldots
  x_{m(i+3,N)}\big\}$, where the $\textit{modulo}$ function $m$ is
  defined in Section~\ref{subsec: smoothness classifier}, and where
  $W^j_i = q_{c, r}(|x_j - x_{\widetilde{m}(i, j, N)}|)$. Here, for
  $s = \frac{c}{2} + r$ we have set
\begin{equation} \label{eq: W non-periodic}
    \widetilde{m}(i, j, N) =
    \left\lbrace
        \begin{array}{ccc}
            j + N - i & \mbox{if} & j < 2s  \mbox{ and } N - 2s + j \leq i \leq N - 1\\
            j - 1 - i & \mbox{if} & N - 2s  < j \leq N - 1 \mbox{ and } 0 \leq i \leq 2s - (N - j) \\
            i & \mbox{else.} \\
        \end{array}\right.
\end{equation}
The values $c=0$ and $r=9$ considered previously are once again used
in the periodic context.
\end{remark}

\begin{figure}
    \centering
    \includegraphics[width=0.7\textwidth]{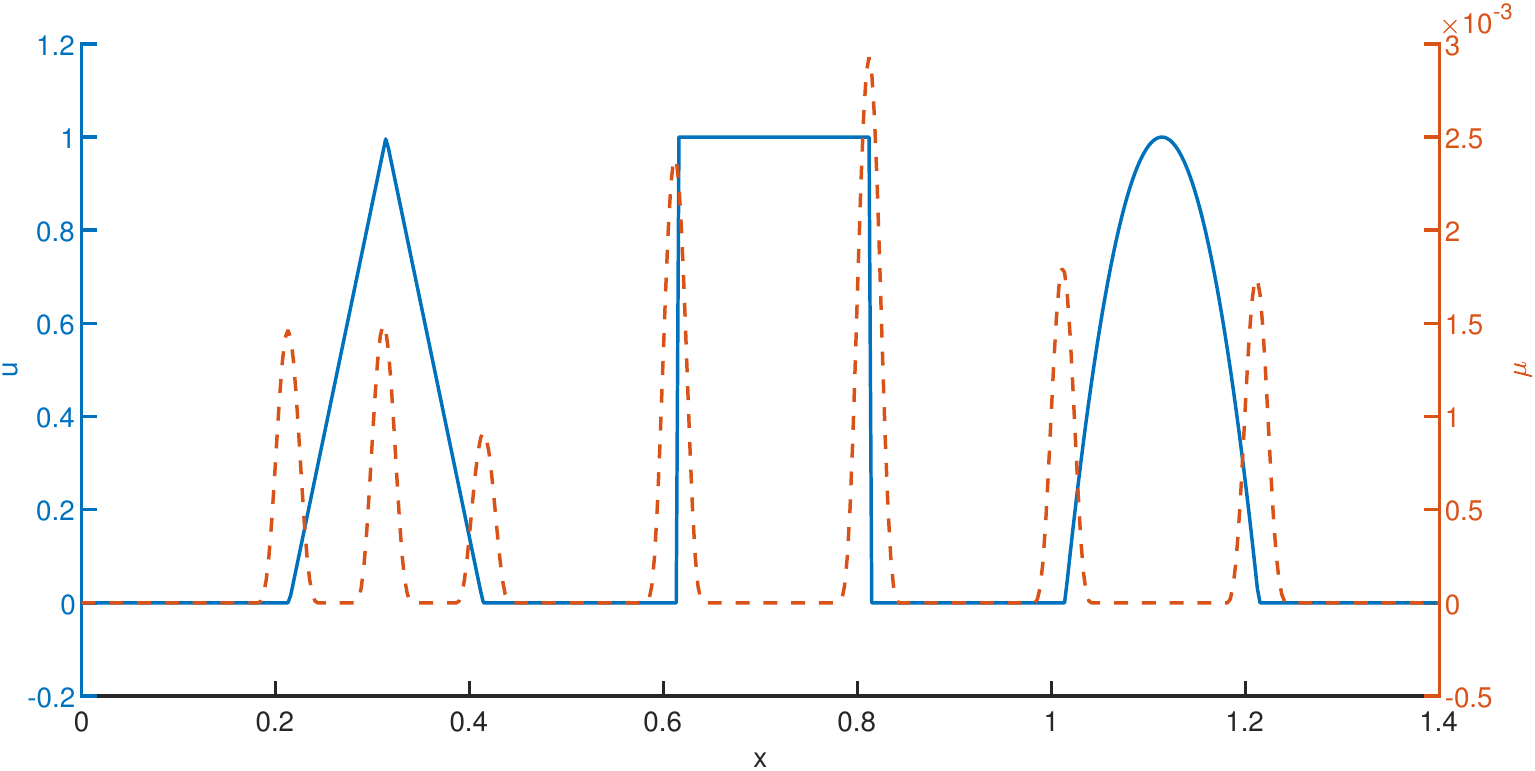} 
    \caption{Viscosity assignment (orange dashed line) resulting from
      application of the SDNN-localized artificial viscosity algorithm
      to the initial condition~\eqref{eq: linear advection IC} (blue
      solid line), using $N = 500$ discretization points.}
\label{fig: viscosity_analysis}
\end{figure}

As an example, Figure~\ref{fig: viscosity_analysis} displays the
viscosity assignments produced, by the method described in this
section for the function
\begin{equation} \label{eq: linear advection IC}
    u(x)=
    \left\lbrace
        \begin{array}{cl}
            10(x - 0.2) & \mbox{if }  0.2 < x \leq 0.3\\
            10(0.4 - x) & \mbox{if } 0.3 < x \leq 0.4\\
            1 & \mbox{if } 0.6 < x \leq 0.8\\
            100(x - 1)(1.2 - x) & \mbox{if } 1 < x \leq 1.2\\
            0 & \mbox{otherwise.}
        \end{array}\right.
\end{equation}
in the interval $[0, 1.4]$. Clearly, the viscosity profiles are smooth
and they are supported around points where the function $u$ is not
smooth.

\subsubsection{Two-dimensional case}
\label{subsubsec: 2d case}
The definition of the 2D viscosity operator follows similarly as the
one for the 1D case, with adequately modified versions of the
underlying functions and operators. In detail, in the present 2D case
we define $\widetilde R$ and $R$ ($\widetilde R[\eta]_{ij} = R(\eta_{ij})$) as
in the 1D case, but using the values $R(1) = 1.5$, $R(2)=1$,
$R(3) = 0.5$, and $R(4)=0$. The $7\times 7$ localization stencils
$L^{i, j}$ ($(i, j) \in \mathcal{I}$) are defined in terms of the 1D
localization stencils $L^i$ and $L^j$ via the relation
$L^{i, j}= L^i \times L^j$. For the 2D scalar Burgers equation, the
MWSB operator $S[\textit{\textbf{e}}]$ (resp. the discrete operator
$\widetilde S[\textit{\textbf{e}}_h]$) is taken to equal the maximum
characteristic speed, that is $S[\textit{\textbf{e}}] = |u|$ (resp.
$\widetilde S[\textit{\textbf{e}}_h]_{ij} = |u_{ij}|$ for
$(i, j) \in \mathcal{I}$). In the case of the 2D Euler problem, the
MWSB operator $S$ used in this paper assigns to $\textbf{\textit{e}}$
the upper bound $S(\textbf{\textit{e}}) = |u| + |v| + a$ on the speed
of propagation $\mathbf{u} \cdot \vec{\kappa} + a$ of the wave
corresponding to the largest eigenvalue of the 2D Flux-Jacobian
(which, in a direction supported by the unit vector $\vec{\kappa}$,
equals $\mathbf{u} \cdot \vec{\kappa} + a$~\cite[Sec. 16.3 and
16.5]{hirsch1990numerical}), so that the discrete operator $\widetilde S$
we propose is given on the grid by
\begin{equation}\label{Flux_der_Euler_2}
  \widetilde S[\textit{\textbf{e}}_h]_{ij} = |u_{ij}| + |v_{ij}|  + a_{ij}.
\end{equation}
Note that this MWSB operator, which equals $a$ plus the sum of the
absolute values of the components of the velocity vector $\mathbf{u}$,
differs slightly from the upper-bound selected
in~\cite{guermond2011entropy, kornelus2018flux}, where the
(equivalent) Euclidean-norm of $\mathbf{u}$ was used instead.  The
two-dimensional local-window operator, finally, is given by
\begin{equation} \label{eq:Lambda2D} \Lambda[b]_{ij} =
  \sum_{(k, \ell) \in \mathcal{I}}\check{W}^k_i  \check{W}^\ell_j b_{k\ell}
\end{equation}
---which, clearly, can be obtained in practice by applying
consecutively the 1D local window operator introduced in
Section~\ref{subsubsec:1d_case} in the horizontal and vertical
directions consecutively.  As in the 1D case, using these operators
and functions, we define the 2D artificial viscosity operator
\begin{equation} \label{eq: V2D} \tilde\mu[\textit{\textbf{e}}_h]_{ij} =
  \Lambda[\widetilde{R}(\tilde\tau_{xy}[\bm{\phi}])]_{ij} \cdot \max_{(k,
    \ell)\in L^{ij}}(S[\textit{\textbf{e}}_h]_{k \ell}) h.
\end{equation}

\subsection{\label{filtering}Spectral filtering}

Spectral methods regularly use filtering strategies in order to
control the error growth in the unresolved high frequency modes. One
such ``global'' filtering strategy is employed in the context of this
paper as well, in conjunction with
FC~\cite{albin2011spectral,amlani2016fc}, as detailed in
Section~\ref{global filtering}. Additionally, a new ``localized''
filtering strategy~\ref{localized filtering} is introduced in this
paper, in order to regularize discontinuous initial conditions, while
avoiding the over-smearing of smooth flow-features. Details regarding
the global and localized filtering strategies are provided in what
follows.

\subsubsection{Global filtering strategy \label{global filtering}}

As indicated above, the proposed algorithm employs spectral filters in
conjunction with the FC method to control the error growth in
unresolved high frequency modes~\cite{albin2011spectral,amlani2016fc}.
For a given Fourier Continuation expansion
\[
  \sum_{k = -M}^{M} \hat{F}^c_k \exp (2\pi i k x / \beta)
\]
the corresponding globally-filtered Fourier Continuation expansion is
given by
\begin{equation} \label{eq: filtering}
  \widetilde F_g = \sum_{k = -M}^{M}
  \hat{F}^c_k \sigma \Big( \frac{2k}{N + C} \Big) \exp (2\pi i k x /
  \beta)
\end{equation}
where
\begin{equation*} \label{eq: filter function} \sigma \Big( \frac{
    2k}{N+C} \Big) = \exp\Big( -\alpha_f \Big(\frac{2k}{N+C} \Big)^{p_f}
  \Big)
\end{equation*}
for adequately chosen values of the positive integer $p_f$ and the
real parameter $\alpha_f>0$. For applications involving
two-dimensional domains the spectral filter is applied sequentially,
one dimension at a time.

In the algorithm proposed in this paper, all the components of the
unknown solution vector $\textbf{\textit{e}}$ are filtered using this
procedure at every time step following the initial time, using the
parameter values $\alpha_f = 10$ and $p_f = 14$, as indicated in
Algorithm~\ref{alg_one}.

\subsubsection{Localized discontinuity-smearing for initial
  data\label{localized filtering}}

In order to avoid the introduction of spurious oscillations arising
from discontinuities in the initial condition, the spectral filter
considered in the previous section is additionally applied, in a
modified form, before the time-stepping process is initiated. A
stronger filter is used to treat the initial conditions, however,
since, unlike the flow field for positive times, the initial
conditions are not affected by artificial viscosity. In order to avoid
unduly degrading the representation of the smooth features of the
initial data, on the other hand, a localized discontinuity-smearing
method, based on use of filtering and windowing is used, that is
described in what follows.

We first present the discontinuity-smearing approach for a 1D
function $F:I\rightarrow \mathbb{R}$, defined on a one-dimensional
interval $I$, which is discontinuous at a single point $z \in I$. In
this case, the smeared-discontinuity function $\widetilde{F}_\mathrm{sm}$,
which combines the globally filtered function $\widetilde F_g$ in a
neighborhood of the discontinuity with the unfiltered function
elsewhere, is defined by
\begin{equation}\label{q}
  \widetilde{F}_\mathrm{sm}(x) =  q_{c, r}(x - z) \widetilde F_g(x) +  (1 - q_{c, r}(x-z)) F(x),
\end{equation}

\begin{remark}
  Throughout this paper, windows $q_{c, r}$ with $c=18$ and $r=9$ and
  globally filtered functions $\widetilde F_g$ with filter parameters
  $\alpha_f = 10$ and $p_f = 2$, which is depicted in
  Figure~\ref{fig:visc_filt_window} right, were used for the
  initial-condition filtering problem, except as noted below in cases
  resulting in window overlap.
\end{remark}
In case multiple discontinuities exist the procedure is repeated
around each discontinuity point. Should the support of two or more of
the associated windowing functions overlap, then each group of
overlapping windows is replaced by a single window which equals zero
outside the union of the supports of the windows in the group, and
which equals one except in the rise regions for the leftmost and
rightmost window functions in the group.

\begin{algorithm}
  \small
  \begin{algorithmic}[1]
    \State \textbackslash \textbackslash Initialization.
    \State Input the trained ANN weights and biases (Section~\ref{subsubsec:
nn archirecture and training}).
  \State  Initialize the unknown solution vector $\textit{\textbf{e}}_h$ (Section~\ref{laws}) to  the given initial-condition values over the given spatial grid.
  \State Initialize time: $t = 0$.
  \While {$t < T$}  
  \State Evaluate the proxy variable $\bm{\phi}$ corresponding to $\textit{\textbf{e}}_h$ at all spatial grid points.
 \State Obtain the smoothness classification values ($\tilde \tau[\bm{\phi}]$, equation~\eqref{tau}, in the 1D case, or $\tilde{\tau}_{xy}[\bm{\phi}]$, equation \eqref{tauxy}, in the 2D case) at all grid points by applying steps~\eqref{pt1} through~\eqref{pt6} in Section~\ref{subsec: smoothness classifier} as required in each case, 1D or 2D.
 \State Evaluate the MWSB operator $\widetilde S [\textit{\textbf{e}}_h]$ at all spatial grid points (Section~\ref{subsubsec:1d_case} in the 1D cases, and Section~\ref{subsubsec: 2d case} in the 2D cases).
 \State Determine the artificial viscosity assignments $\tilde \mu[\textit{\textbf{e}}_h]$ (Equation~\eqref{eq: V} in the 1D case or equation~\eqref{eq: V2D} in the 2D case) at all spatial grid points.
 \State (Case $t = 0$) Apply localized discontinuity-smearing (Section~\ref{localized filtering}) to the solution vector $\textit{\textbf{e}}_h$ and overwrite $\textit{\textbf{e}}_h$ with the resulting values.
    \State (Case $t > 0$) Apply global filtering (Section~\ref{global filtering}) to the solution vector $\textit{\textbf{e}}_h$ and overwrite $\textit{\textbf{e}}_h$ with the resulting values. 
  \State Evaluate the temporal step-size $\Delta t$ by substituting the discrete version $\widetilde S[\textit{\textbf{e}}_h]$ and $\tilde \mu[\textit{\textbf{e}}_h]$ of $S[\textit{\textbf{e}}]$ and $\mu[\textit{\textbf{e}}]$ in equation~\eqref{eq: CFL}. 
  \State Apply the SSPRK-4 time stepping scheme and FC-based spatial differentiation for given (Dirichlet or Neumann) boundary conditions (Sections~\ref{sec:fc} and~\ref{FC-tstep})  to the discrete version of the  viscous system of equations~\eqref{f_visc}-\eqref{eq: convection diffusion eqn} with $\mu[\textit{\textbf{e}}]$ substituted by  $\tilde \mu[\textit{\textbf{e}}_h]$.
  \State Update time: $t = t + \Delta t$
 \EndWhile   
  \caption{FC-SDNN algorithm}
\label{alg_one}
\end{algorithmic}
\end{algorithm}

In the 2D case the localized filtering strategy is first performed
along every horizontal line $y = y_{j:}$ for $0 \leq j \leq
N_2-1$. The resulting ``partially'' filtered function is then filtered
along every vertical line $x = x_i$ for $0 \leq i \leq N-1$ using the
same procedure.

\begin{remark}\label{filter_vector}
  For PDE problems involving vectorial unknowns, the localized
  discontinuity-smearing strategies for the initial condition is
  applied independently to each component.
\end{remark}

As indicated in Algorithm~\ref{alg_one}, all the components of the
initial ($t=0$) solution vector $\textbf{\textit{e}}$ are filtered
using the procedure described in this section.

\subsection{Algorithm pseudo-code\label{pseudo-code}}
A pseudo-code for the complete FC-SDNN numerical method for the
various equations considered in this paper, and for both 1D and 2D
cases, is presented in Algorithm~\ref{alg_one}.

\section{\label{sec:numerical results}Numerical results}

This section presents results of application of the FC-SDNN method to
a number of non-periodic test problems (with the exception of a
periodic linear-advection problem, in Section~\ref{dispersionless},
demonstrating the limited dispersion of the method), time-dependent
boundary conditions, shock waves impinging on physical boundaries
(including a corner point and a non rectangular domain), etc.  All of
the examples presented in this section resulted from runs on Matlab
implementations of the various methods used. Computing times are not
reported in this paper in view of the inefficiencies associated with
the interpreter computer language used but, for reference,
we note
from~\cite{shahbazi2011multi} that, for the types of equations
considered in this paper, the FC implementations can be quite
competitive, in terms of computing time, for a given accuracy. Our
experiments indicate, further, that the relative cost of application
of smooth-viscosity operators of the type used in this paper decreases
as the mesh is refined, and that for large enough discretizations the
relative cost of the neural network algorithm becomes
insignificant. We thus expect that, as in~\cite{shahbazi2011multi},
efficient implementations of the proposed FC-SDNN algorithm will prove
highly competitive for general configurations.

\subsection{\label{subsec:linear advection results}Linear advection}
The simple 1D linear-advection results presented in this section
demonstrate, in a simple context, two main benefits resulting from the
proposed approach, namely 1)~Effective handling of boundary conditions
(Section~\ref{bound_1D}); and 2)~Essentially dispersionless character
(Section~\ref{dispersionless}). For the examples in this section FC
expansions with $d = 5$ were used.
\subsubsection{Boundary conditions \label{bound_1D}}
\begin{figure}
     \centering
         \includegraphics[width=0.32\textwidth]{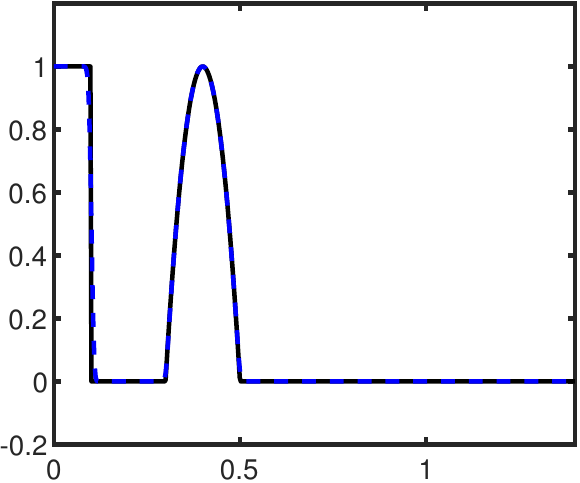}
         \includegraphics[width=0.32\textwidth]{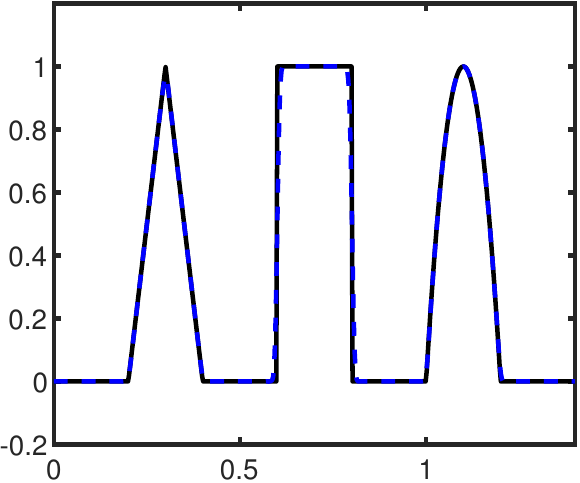}      
         \includegraphics[width=0.32\textwidth]{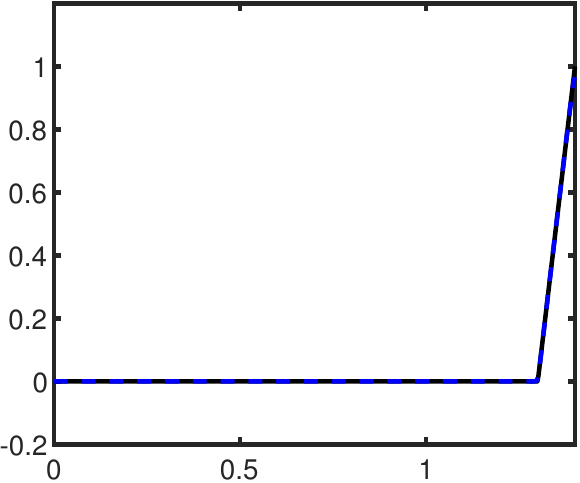}         
         \caption{Solution of the non-periodic one dimensional linear
           advection problem at three different points in time, using
           $N = 500$ discretization points. Exact (black solid
           line), FC-SDNN (blue dashed line). Left: $t = 0.5$. Center:
           $t = 1.3$. Right: $t = 2.3$.}
    \label{fig: LA non-periodic}
\end{figure}

Figure~\ref{fig: LA non-periodic} displays the FC-SDNN solution to the
linear advection problem~\eqref{eq: linear advection}, in which three
waves with various degrees of smoothness emanate from the left
boundary, on the interval $[0, 1.4]$ and with an initial condition
$u(x, 0) = 0$ and boundary condition at $x = 0$ given by
\begin{equation} \label{eq: linear advection BC}
    u(0, t)=
    \left\lbrace
        \begin{array}{cl}
            100 t (t-0.2) & \mbox{if }  0 < t < 0.2\\
            1 & \mbox{if } 0.2 < t < 0.4\\
            10(t - 0.8) & \mbox{if } 0.8 < t < 0.9\\
            1 - 10(t - 0.9) & \mbox{if } 0.9 < t < 1\\
            0 & \mbox{otherwise.}
        \end{array}\right.
\end{equation}
At the outflow boundary $x = 1.4$, the solution was evolved similarly
as for the interior points. The solution was computed up to time
$T = 2$, using an adaptive time step given by~\eqref{eq: CFL}, with
$\textrm{CFL} = 2$. As shown in Figure \ref{fig: LA non-periodic}, the
waves travel within the domain and match the position of the exact
solution, showcasing the dispersionless character of the solver. The
introduction of waves and irregularities at times $t = 0$, $t = 0.2$,
$t = 0.4$, $t = 0.6$, $t = 0.8$, $t = 0.9$, $t = 1$ through the left
boundary is automatically accompanied by the assignment of artificial
viscosity on a small area near that boundary. The FC-SDNN algorithm
then stops assigning viscosity after a short time, when the
irregularities are sufficiently smeared. This is a desirable property
of the original algorithm~\cite{schwander2021controlling} which is
preserved in the present context. Finally, the waves exit the domain
without producing undesired reflective artifacts around the boundary.

\subsubsection{\label{dispersionless}Limited dispersion}
To demonstrate the limited dispersion inherent the FC-SDNN algorithm
we consider a problem of cyclic advection of a ``bump'' solution over
a bounded 1D spatial domain---thus effectively simulating, in a
bounded domain, propagation over arbitrarily extended spatial
regions. To do this we utilize the smooth cut-off ``bump'' function
$\omega = \omega(x, q_1, q_2)$ ($q_1 \neq q_2$) defined by
\begin{equation}
  \omega(x, q_1, q_2) =
  \left\lbrace
    \begin{array}{ccc}
      \exp(2\frac{e^{-1/\xi}}{\xi - 1}) & \mbox{if} & q_1 \leq |x| \leq q_2 \\
      1 & \mbox{if} & |x| \leq q_1  \\
      0 & \mbox{if} & |x| \geq q_2,  \\
    \end{array}\right. \quad\mbox{where}\quad \xi = \frac{|x| - q_1}{q_2 - q_1}.
\end{equation}

For this example we solve the equation~\eqref{eq: linear advection}
with $a=1$, starting from a smooth initial condition given by
$u(x,0) = \omega(x - 0.5, 0, 0.2)$, over the domain
$\left[0, 1 \right]$ under periodic boundary conditions. In order to
enforce such periodic conditions, the FC differentiation scheme is
adapted, by using the same precomputed matrices $A_{\ell}$, $A_r$ and
$Q$ (see Section~\ref{sec:fc}) in conjunction with the "wrapping"
procedure described in~\cite[Sec. 3.3]{albin2011spectral}. Using the
FC-SDNN algorithm of order $d=5$ and SSPRK-4, the solution was evolved
for five hundred periods, up to $T = 500$; the resulting $t=500$
solution $u(x,500)$ is displayed in Figure~\ref{fig: LA
  dispersion}. For comparison this figure presents numerical results
obtained by means of a 6-th order central finite-difference scheme
(also with SSPRK-4 time stepping). The finite-difference method uses a
constant time-step $\Delta t = 0.0034$, while the FC-SDNN uses the
adaptive time step defined in \eqref{eq: CFL} for which, in the
present case, with $\mu = 0$ and $S=1$, we have $\Delta t =
0.0036$. Clearly, the FC-SDNN solution matches the exact solution
remarkably well even after very long times, showcasing the low
dissipation and dispersion afforded by the FC-based approach. In
comparison, the higher-order finite difference solution suffers
greatly from the dissipation and dispersion effects. Thus, significant
advantages arise from use the FC-based strategy as well as the
localized-viscosity approach that underlies the FC-SDNN method.

\begin{figure}[H]
     \centering
     \begin{subfigure}[b]{0.28\textwidth}
         \centering
         \includegraphics[width=\textwidth]{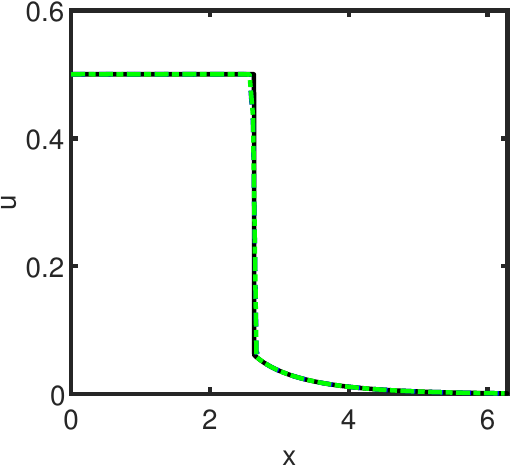}
         \caption{$t = 2 \pi$}
     \end{subfigure}%
     \begin{subfigure}[b]{0.28\textwidth}
         \centering
         \includegraphics[width=\textwidth]{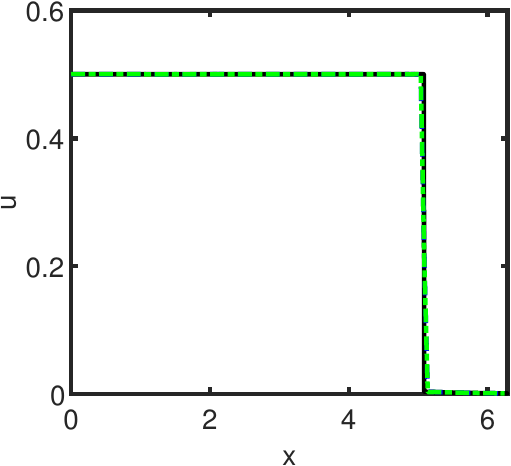}
         \caption{$t = 5 \pi$}
     \end{subfigure}%
     \begin{subfigure}[b]{0.28\textwidth}
         \centering
         \includegraphics[width=\textwidth]{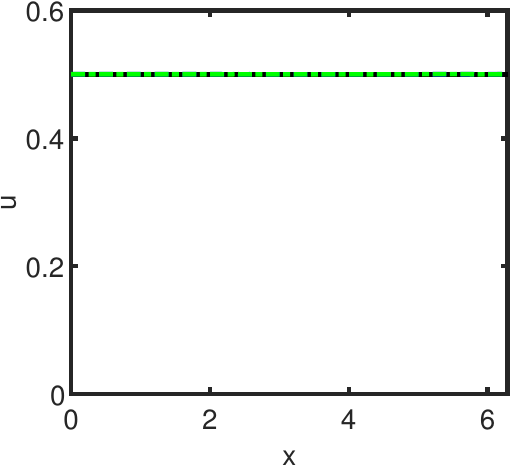}
         \caption{$t = 8 \pi$}
     \end{subfigure}  
     \begin{subfigure}[b]{0.5\textwidth}
         \centering
         \includegraphics[width = 0.7\textwidth]{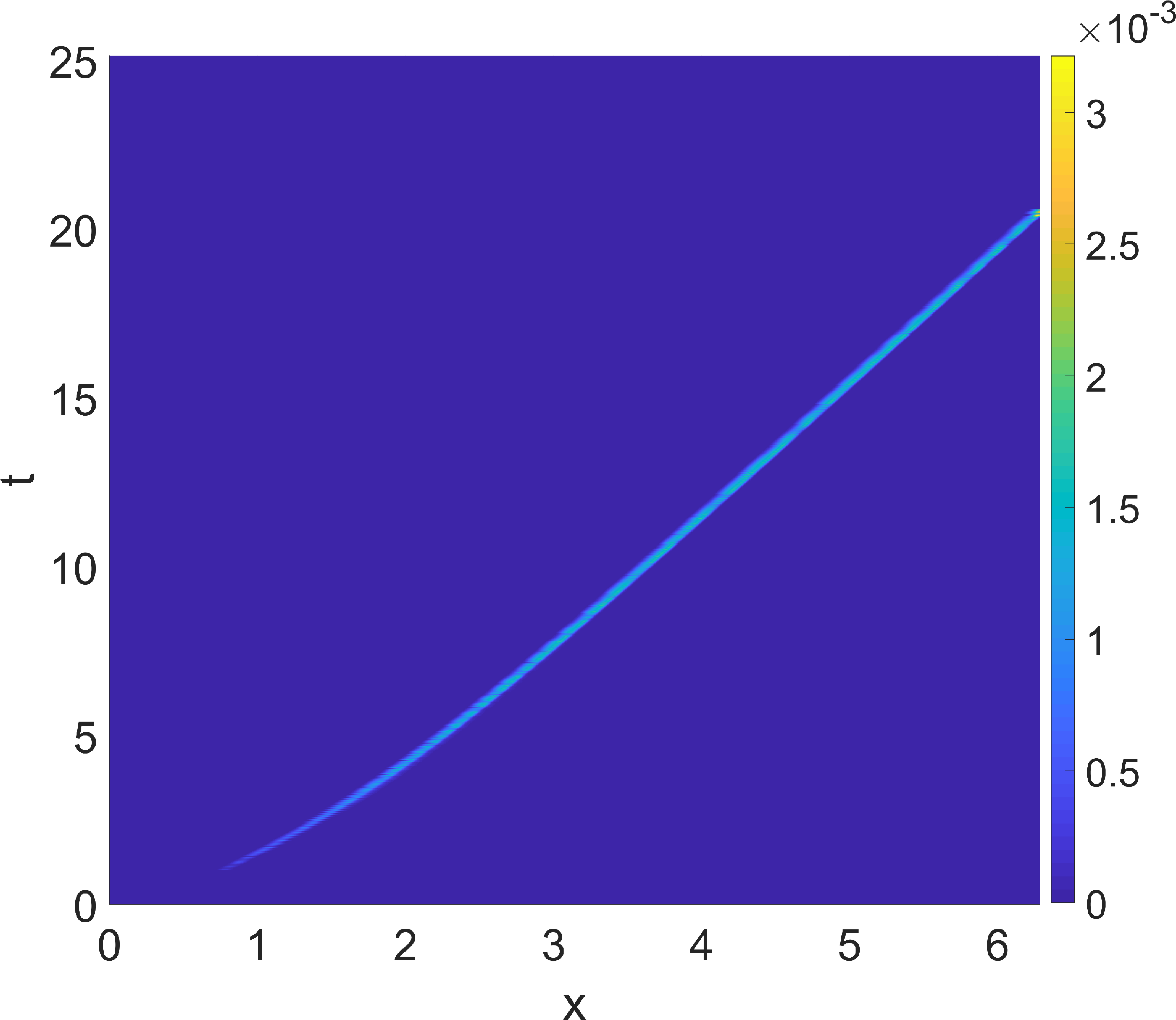}
         \caption{Time history of FC-SDNN artificial viscosity.}
     \end{subfigure}       
     \caption{Solutions to the non-periodic one dimensional Burgers
       equation produced by the FC-SDNN and FC-EV algorithms of order
       $d = 2$ at three different times $t$. Black solid line: finely
       resolved FC-SDNN ($N=10,000$, for reference). Blue dashed line:
       FC-SDNN with $N=500$. Green dot-dashed line: FC-EV with $N=500$.}
\label{fig: Burgers 1D}
\end{figure}

\subsection{\label{subsec:burgers results}Burgers equation}

The 1D and 2D Burgers equation tests presented in this section
demonstrate the FC-SDNN solver's performance for simple nonlinear
non-periodic problems. In particular, the 2D example showcases the
ability of the algorithm to handle multi-dimensional problems where
shocks intersect domain boundaries (a topic that is also considered in
Section~\ref{Euler 2D results} in the context of the Euler equations).

\subsubsection{\label{subsubsec:1D burgers results}1D Burgers equation} 

Figure~\ref{fig: Burgers 1D} displays solutions to the 1D Burgers
equation (\ref{eq: burgers 1d}), where a shock forms from the sharp
features in the initial condition
\begin{equation} \label{eq: 1d burgers IC}
    u_0(x)=\frac{1}{\exp(x - \frac{3}{20})[\tanh(10x - 3) + 1] - \tanh(10x-3) + 1}.
  \end{equation}
  Results of simulations produced by means of the FC-SDNN and FC-EV
  algorithms are presented in the figure. In both cases Dirichlet
  boundary conditions at the inflow boundary $x = 0$ were used, while
  the solution at the outflow boundary $x = 2\pi$ was evolved
  numerically in the same manner as the interior domain points. As
  shown in the figure, the shock is sharply resolved by both
  algorithms, with no visible oscillations. In both cases the shock
  eventually exits the physical domain without any undesired
  reflections or numerical artifacts. For this example FC expansions
  with $d = 2$ were used. The reference solution (black) was computed
  on a 10000-point domain, with the FC-SDNN method.

\subsubsection{\label{subsubsec:2D burgers results}2D Burgers equation}

\begin{figure}
\centering
  \begin{subfigure}[b]{0.35\textwidth}
    \includegraphics[width=1\textwidth]{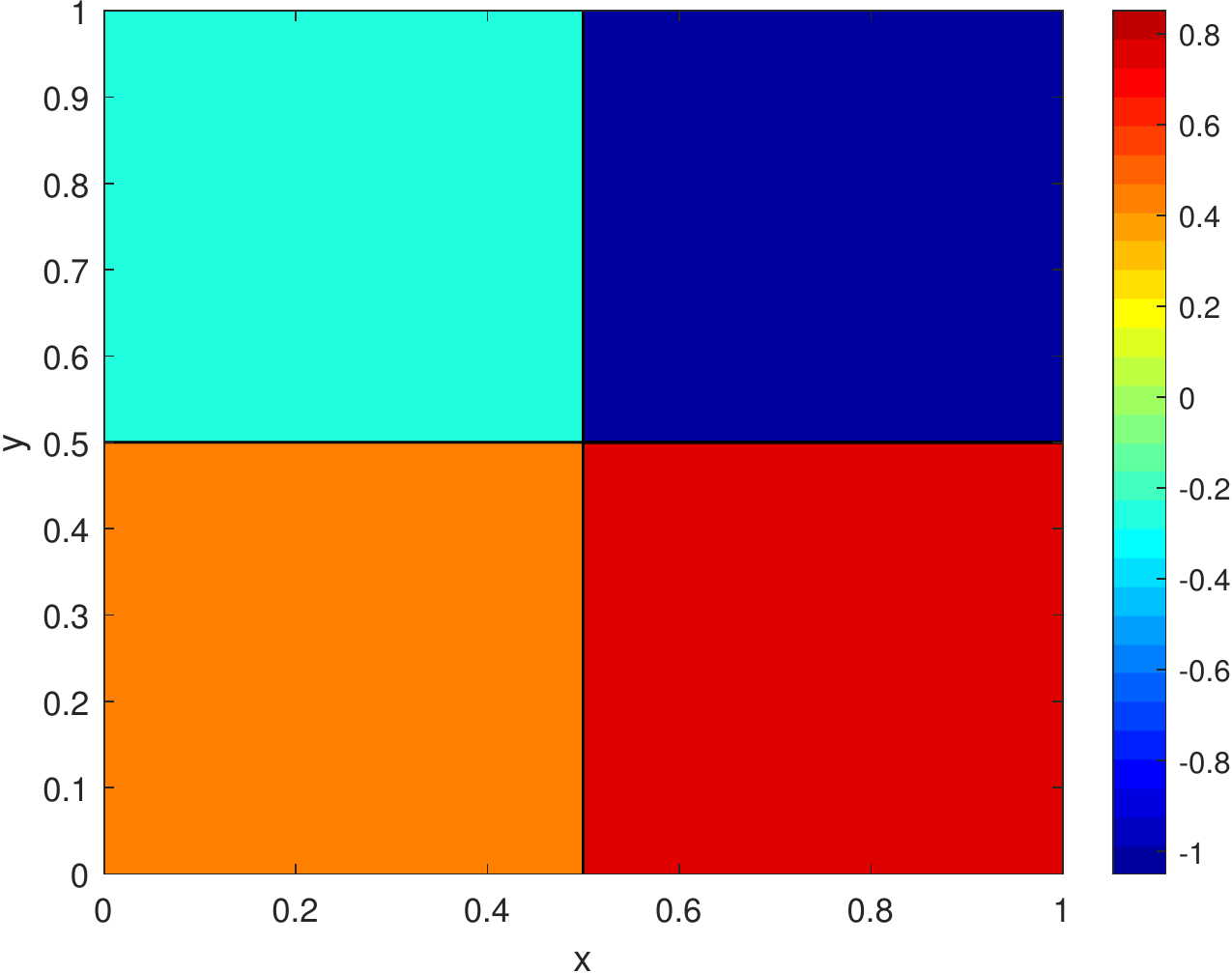}
    \caption{Initial condition}
    \label{fig:Burgers 2D initial}
  \end{subfigure}
  \begin{subfigure}[b]{0.35\textwidth}
    \includegraphics[width=\textwidth]{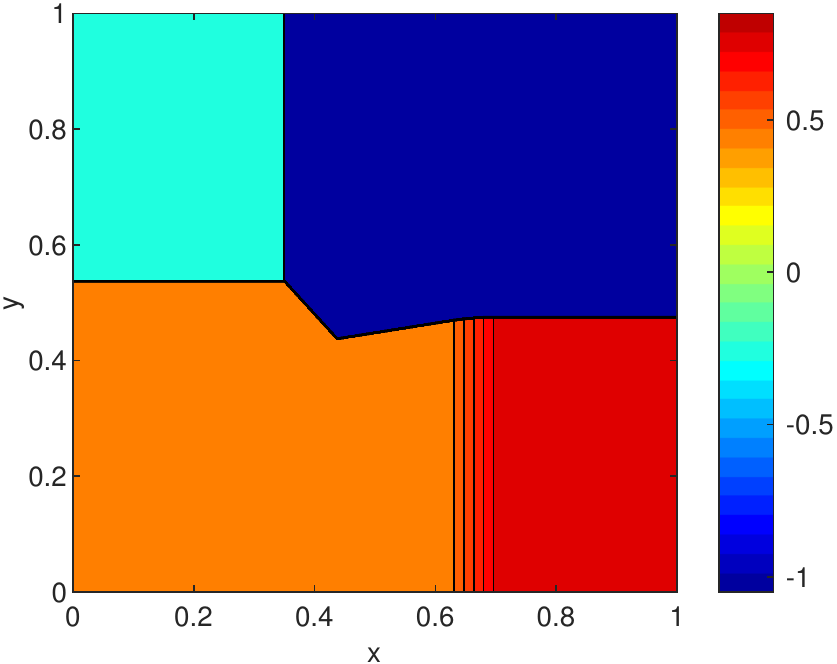}
    \caption{Exact solution}
    \label{fig:Burgers 2D exact}
  \end{subfigure}
  \begin{subfigure}[b]{0.35\textwidth}
    \includegraphics[width=1\textwidth]{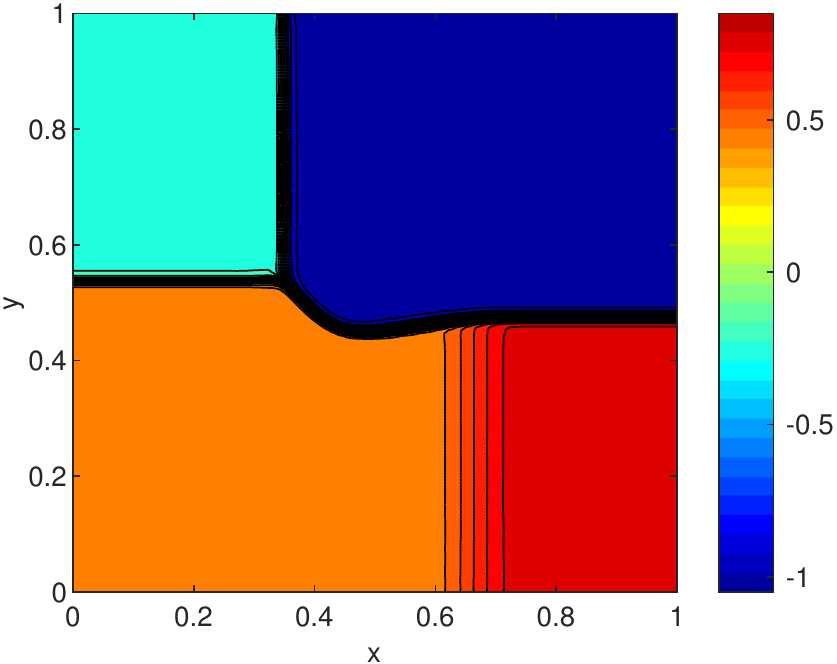}
    \caption{N = 200}
    \label{fig:Burgers 2D solution N=200}
  \end{subfigure}
  \begin{subfigure}[b]{0.35\textwidth}
    \includegraphics[width=\textwidth]{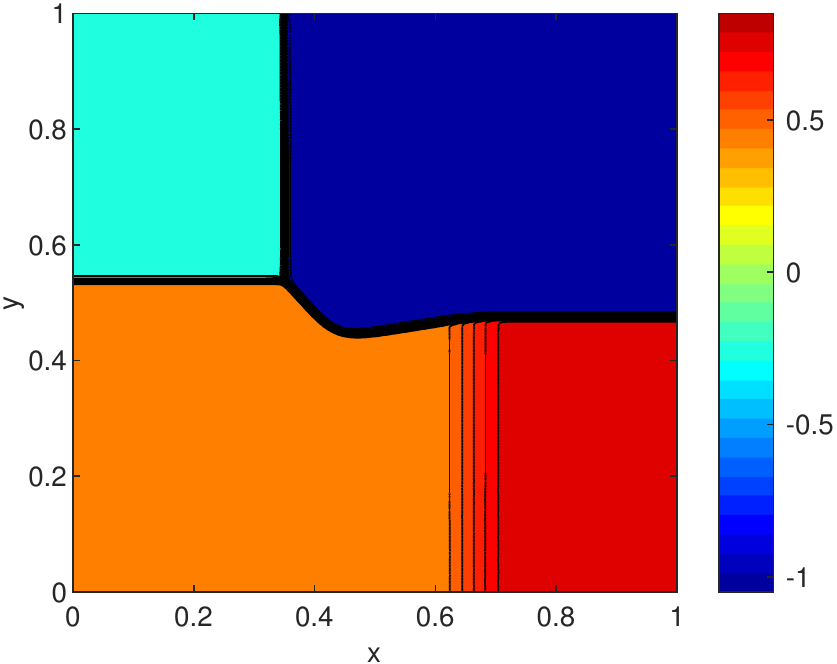}
    \caption{N = 400}
    \label{fig:Burgers 2D solution N=400}
  \end{subfigure}
  \begin{subfigure}[b]{0.35\textwidth}
    \includegraphics[width=\textwidth]{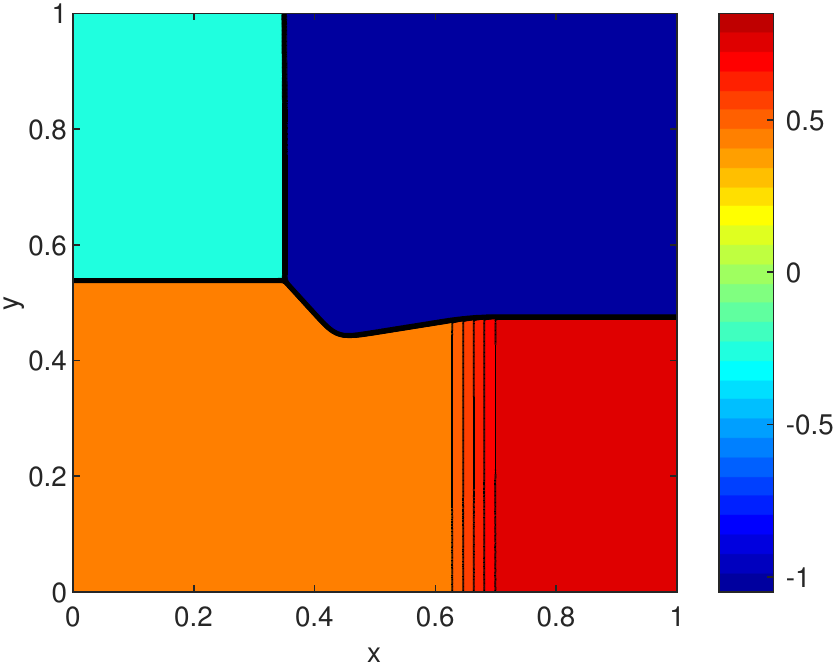}
    \caption{N = 1000}
    \label{fig:Burgers 2D solution N=1000}
  \end{subfigure}
  \begin{subfigure}[b]{0.35\textwidth}
    \includegraphics[width=\textwidth]{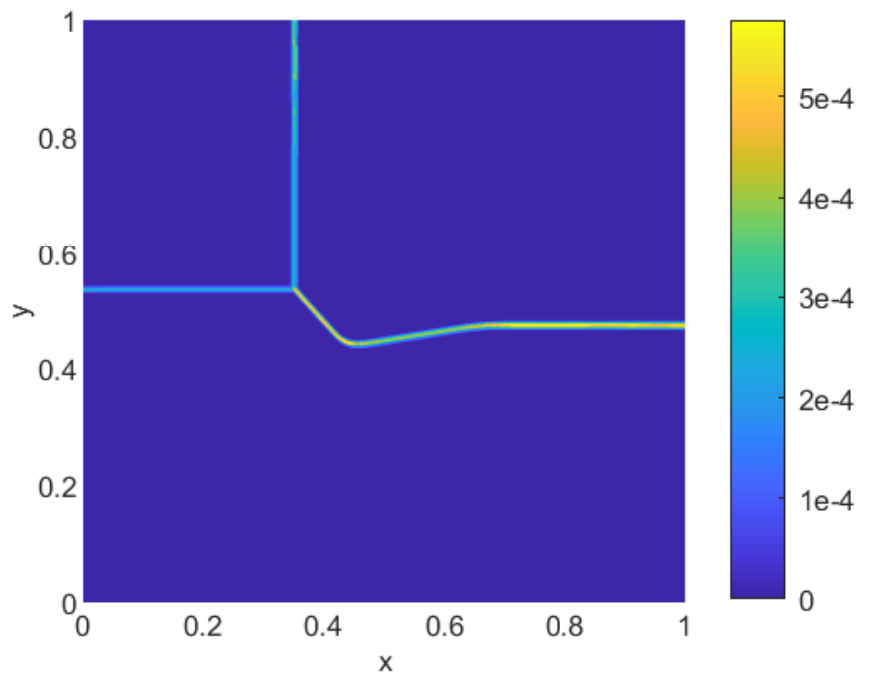}
    \caption{Art. visc. ($t = 0.25$, $N=1000$)}.
    \label{fig:Burgers 2D viscosity N=1000}
  \end{subfigure}  
  \caption{Fifth-order ($d = 5$) FC-SDNN numerical solution to the 2D
    Burgers equation with initial condition displayed on the upper-left
    panel, whose exact solution at $t = 0.25$ is displayed on the
    upper-right panel. The middle and left-lower panels display the
    FC-SDNN numerical solutions at $t = 0.25$ obtained by using
    $N\times N$ spatial grids with three different values of
    $N$, as indicated in each panel. The FC-SDNN numerical viscosity
    at $t = 0.25$ for the case $N=1000$ is presented in the
    lower-right panel. }
\label{fig: Burgers 2D solution}
\end{figure}

To demonstrate the solver's performance and correct handling of
shock-boundary interactions for 2D problems, we consider the 2D
Burgers scalar equation (\ref{eq: burgers 2d}) on the domain
$\mathcal{D} = \left[0, 1\right] \times \left[0, 1 \right]$, with an
initial condition given by the function
\begin{equation} \label{eq: smooth burgers 2d IC}
    u_0(x)=
    \left\lbrace
        \begin{array}{ccc}
            -1 & \mbox{if} & x \in  [0.5, 1]  \mbox{ and } y \in [0.5, 1]\\
            -0.2 & \mbox{if} & x \in  [0, 0.5]  \mbox{ and } y \in [0.5, 1]\\
            0.5 & \mbox{if} & x \in  [0, 0.5]  \mbox{ and } y \in [0, 0.5]\\
            0.8 & \mbox{if} & x \in  [0.5, 1]  \mbox{ and } y \in [0, 0.5],            
        \end{array}\right.
\end{equation}
and with vanishing normal derivatives at the boundary. This problem
admits an explicit solution (displayed in Figure~\ref{fig:Burgers 2D
  exact} at time $t = 0.25$) which includes three shock waves and a
rarefaction wave, all of which travel orthogonally to various straight
boundary segments.  Figures~\ref{fig:Burgers 2D solution
  N=200},~\ref{fig:Burgers 2D solution N=400}, and~\ref{fig:Burgers 2D
  solution N=1000} present the corresponding numerical solutions
produced by the FC-SDNN algorithm at time $t = 0.25$ resulting from
use of various spatial discretizations and with adaptive time step
given by~\eqref{eq: CFL} with $\textrm{CFL} = 2$. Sharply resolved
shock waves are clearly visible for the finer discretizations, as is
the rarefaction wave in the lower part of the figure. The viscosity
assignments, which are sharply concentrated near shock positions as
the mesh is refined, suffice to avert the appearance of spurious
oscillations. It is interesting to note that non-vanishing viscosity
values are only assigned around shock discontinuities. (The SDNN
algorithm assigns zero viscosity to the rarefaction wave for all
time---as a result of the discontinuity-smearing introduced by the
algorithm on the initial condition for the velocity, described in
Section~\ref{localized filtering}, which the FC-SDNN then preserves for
all times in regions near the rarefaction wave on account of the
resulting smoothness of the numerical solution in such regions). For
this example FC expansions with $d = 5$ were used.

\subsection{\label{subsec:euler results}1D and 2D Euler systems}

This section presents a range of 1D and 2D test cases for the Euler
system demonstrating the FC-SDNN algorithm's performance in the
context of a nonlinear systems of equations. The test cases include
well-known 2D arrangements, including the shock-vortex interaction
example~\cite{shu1999high}, 2D Riemann problem
flow~\cite{lax1998solution}, Mach 3 forward facing
step~\cite{woodward1984numerical}, and Double Mach
reflection~\cite{woodward1984numerical}. In particular the results
illustrate the algorithm's ability to handle contact discontinuities
and shock-shock interactions as well as shock reflection and
propagation along physical and computational boundaries.\looseness = -1

\subsubsection{\label{Euler 1D results}1D Euler problems}

The 1D shock-tube tests considered in this section demonstrate the
solver's ability to capture not only shock-wave discontinuities (that
also occur in the Burgers test examples considered in
Section~\ref{subsec:burgers results}) but also contact
discontinuities. Fortunately, in view of the localized spectral
filtering strategy used for the initial-data (Section~\ref{localized
  filtering}), the algorithm completely avoids the use of artificial
viscosity around contact discontinuities, and thus leads to excellent
resolution of these important flow features. This is demonstrated in a
variety of well known test cases, including the
Sod~\cite{sod1978survey}, Lax~\cite{lax1954weak},
Shu-Osher~\cite{shu1989efficient} and Blast
Wave~\cite{ray2018artificial} problems, with flows going from left to
right---so that the left boundary point (resp. right boundary point)
is the inflow (resp. outflow) boundary. In all cases an adaptive time
step given by~\eqref{eq: CFL} with $\textrm{CFL} = 2$ was used, and,
following~\cite[Sec. 19]{hirsch1990numerical}, inflow (resp. outflow)
boundary condition were enforced at the inflow boundaries
(resp. outflow boundaries) by setting $\rho$ and $u$ (resp $p$)
identically equal, for all time, to the corresponding boundary values
of these quantities at the initial time $t=0$. For the examples in this section FC expansions with $d = 5$ were used.

\begin{description}

\item[Sod problem.]  We consider a Sod shock-tube problem for the 1D
  Euler equations~\eqref{eq: euler 1d equation} on the interval
  $[-4, 5]$ with initial conditions
    \begin{equation*}
        (\rho, u, p)=
        \left\lbrace
            \begin{array}{ccc}
              (1, 0, 1) & \mbox{if} & x < 0.5\\
              (0.125, 0, 0.1) & \mbox{if} & x > 0.5,\\
            \end{array}\right.
    \end{equation*}
    a setup that gives rise (from right to left) to a right-moving
    shock wave, a contact discontinuity and a rarefaction wave (upper
    and middle left images in Figure~\ref{fig: Sod solution}). The
    solution was computed up to time $T = 2$. The results presented in
    Figure~\ref{fig: Sod solution} show well resolved shocks (upper
    and middle right) and contact discontinuities (upper and middle
    center), with no visible
    Gibbs oscillations in any case. The FC-SDNN and FC-EV solvers
    demonstrate a similar resolution in a vicinity of the shock, but
    the FC-SDNN method provides a much sharper resolution of the
    contact-discontinuity. As shown in Figure~(\ref{fig: Sod
      viscosity}), after a short time the FC-SDNN algorithm does not
    assign artificial viscosity in a vicinity of the contact
    discontinuity, leading to the significantly more accurate
    resolution observed for this flow feature.

\begin{figure}
  \begin{subfigure}[t]{1\linewidth}
    \centering
    \includegraphics[width=.33\linewidth]{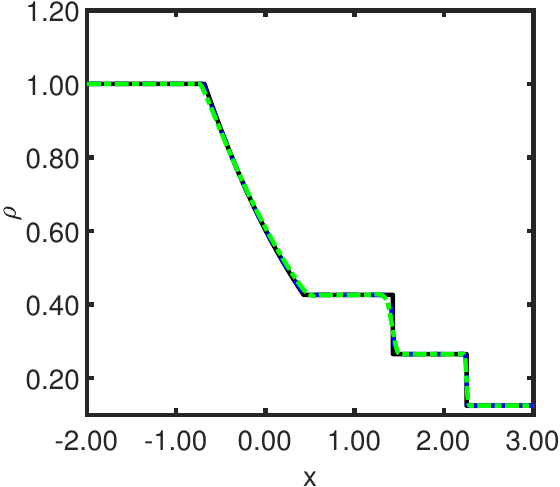}%
    \includegraphics[width=.33\linewidth]{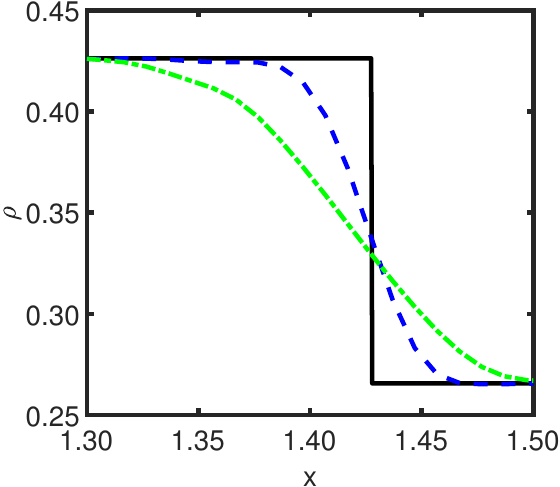}%
    \includegraphics[width=.33\linewidth]{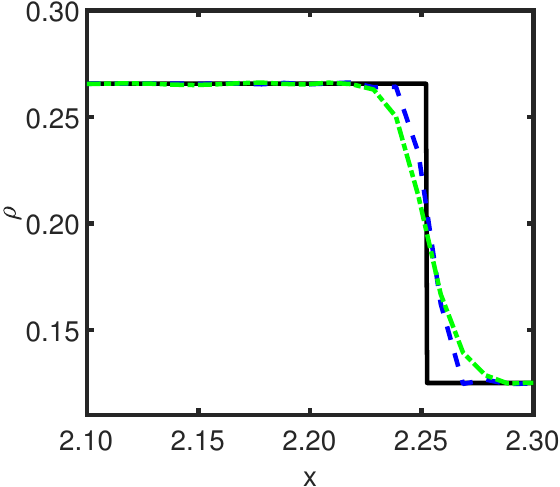}
    \caption{N = 500}
  \end{subfigure}
  \begin{subfigure}[t]{1\linewidth}
    \centering
    \includegraphics[width=.33\linewidth]{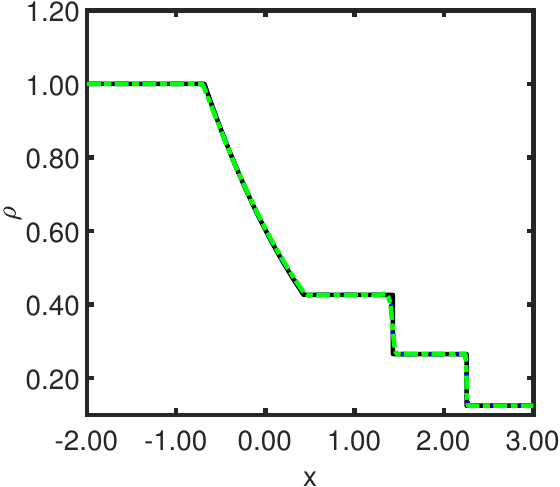}%
    \includegraphics[width=.33\linewidth]{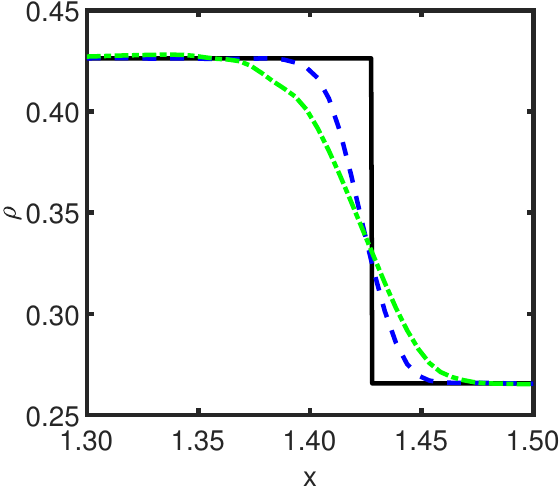}%
    \includegraphics[width=.33\linewidth]{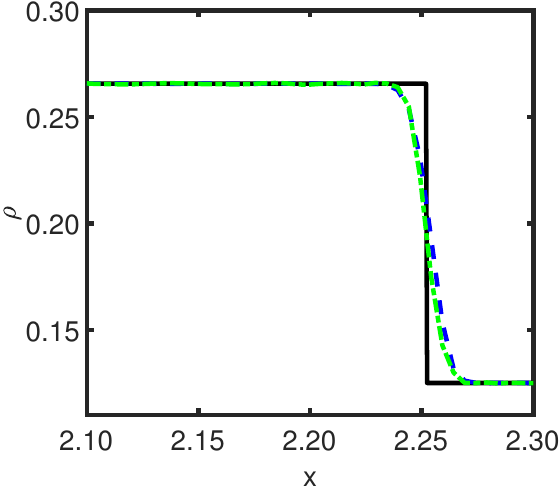}
    \caption{N = 1000}
  \end{subfigure}  
  \begin{subfigure}[t]{1\linewidth}
    \centering
    \includegraphics{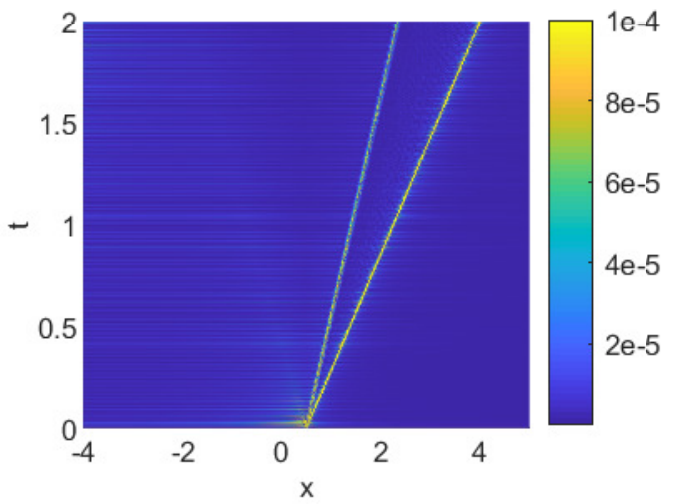}
    \includegraphics{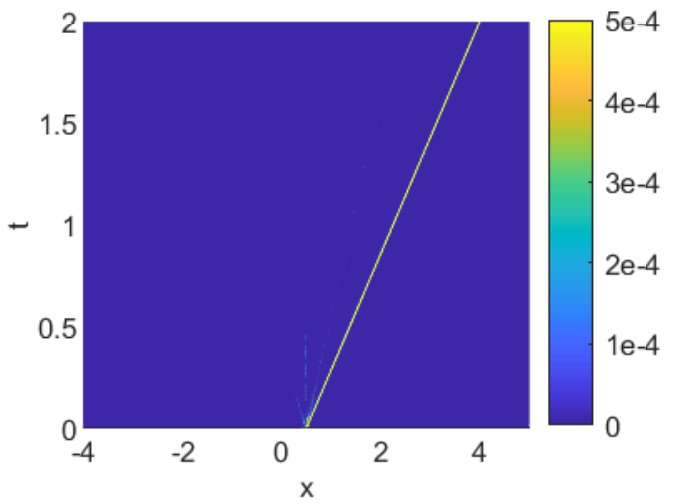}
    \caption{Time history of assigned artificial viscosity for N = 1000, for FC-EV (left) and FC-SDNN (right).}
    \label{fig: Sod viscosity}
  \end{subfigure}        
  \caption{Solutions to the Sod problem produced by the FC-SDNN and
    FC-EV algorithms of order $d = 5$ at $t= 0.2$. Exact solution:
    solid black line. FC-SDNN solution: Blue dashed-line. FC-EV
    solution: Green dot-dashed line. Numbers of discretization points:
    $N=500$ in the upper panels and $N=1000$ in the middle panels.
    Bottom panels: artificial viscosity assignments.}
\label{fig: Sod solution}
\end{figure}

\item[Lax problem.]

  We consider a Lax problem on the interval $[-5, 5]$, with initial
  condition
  \begin{equation*}
      (\rho, u, p)=
      \left\lbrace
            \begin{array}{ccc}
                (0.445, 0.698, 3.528) & \mbox{if} & x < 0\\
                (0.5, 0, 0.571) & \mbox{if} & x > 0\\
            \end{array}\right.
        \end{equation*}
        which results in a combination (from right to left) of a shock
        wave, a contact discontinuity and a rarefaction wave (upper
        and middle left images in Figure~\ref{fig: Lax solution}). The
        solution was computed up to time $T = 1.3$.  The results are
        presented in Figure~\ref{fig: Lax solution}, which shows
        well resolved shocks (upper and middle right) without
        detectable Gibbs oscillations. The viscosity time history
        displayed in Figure~\ref{fig: Lax viscosity} shows that the
        FC-SDNN method only assigns artificial viscosity in the
        vicinity of the shock discontinuity but, as discussed above,
        not around the contact discontinuity, leading to a sharper
        resolution by the FC-SDNN method in this region (upper and
        middle center). The shock resolution is similar for the EV and
        SDNN algorithms, but the latter approach is significantly more
        accurate around the contact discontinuity.

\begin{figure}
  \begin{subfigure}[t]{1\linewidth}
    \centering
    \includegraphics[width=.33\linewidth]{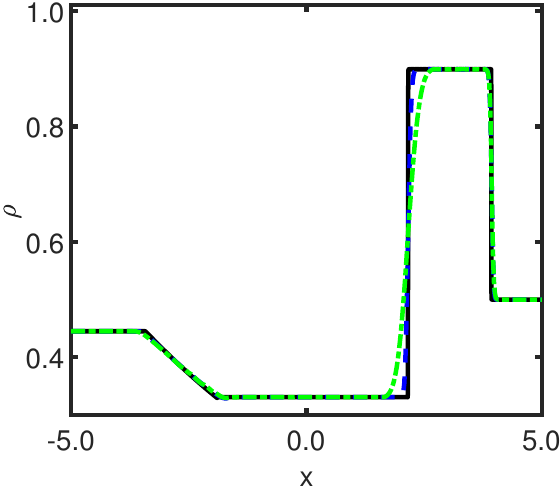}%
    \includegraphics[width=.33\linewidth]{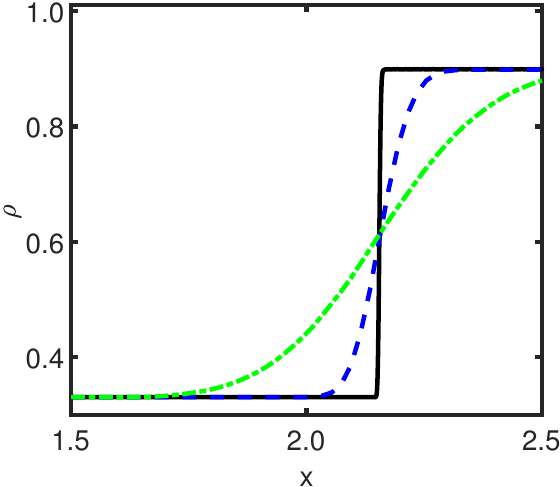}%
    \includegraphics[width=.33\linewidth]{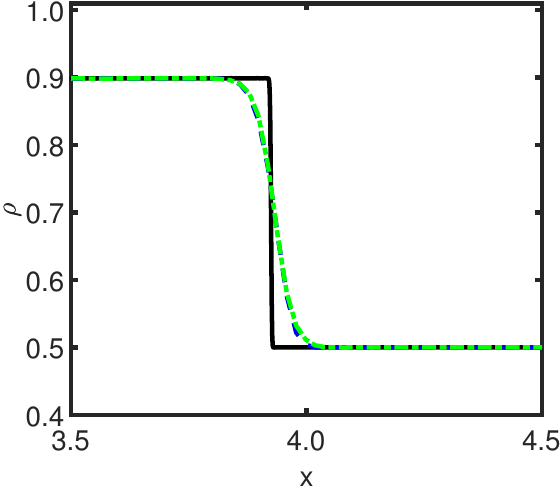}
    \caption{N = 500}
  \end{subfigure}
  \begin{subfigure}[t]{1\linewidth}
    \centering
    \includegraphics[width=.33\linewidth]{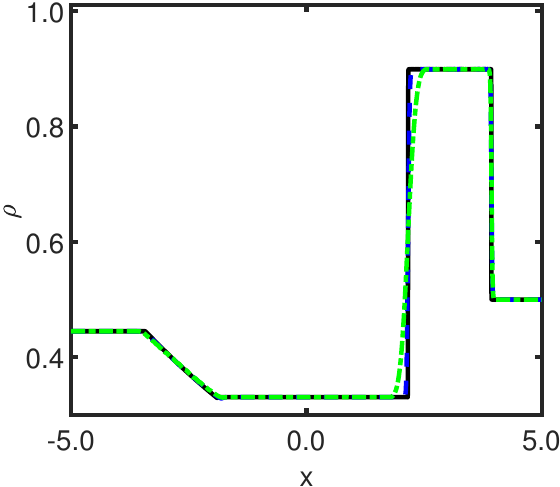}%
    \includegraphics[width=.33\linewidth]{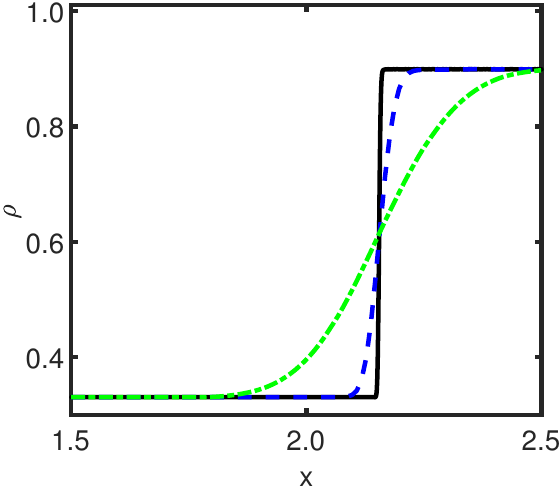}%
    \includegraphics[width=.33\linewidth]{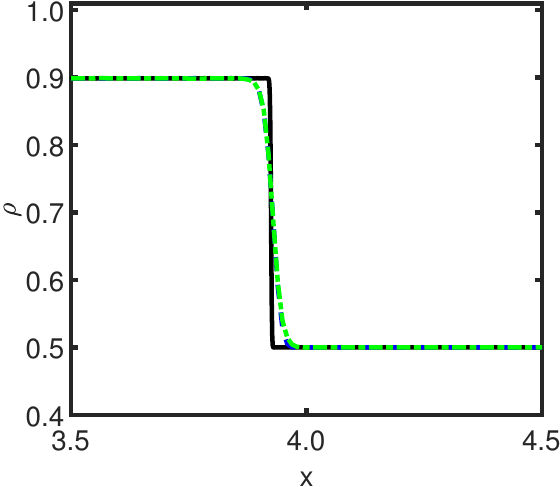}
    \caption{N = 1000}
  \end{subfigure}  
  \begin{subfigure}[t]{1\linewidth}
    \centering
    \includegraphics{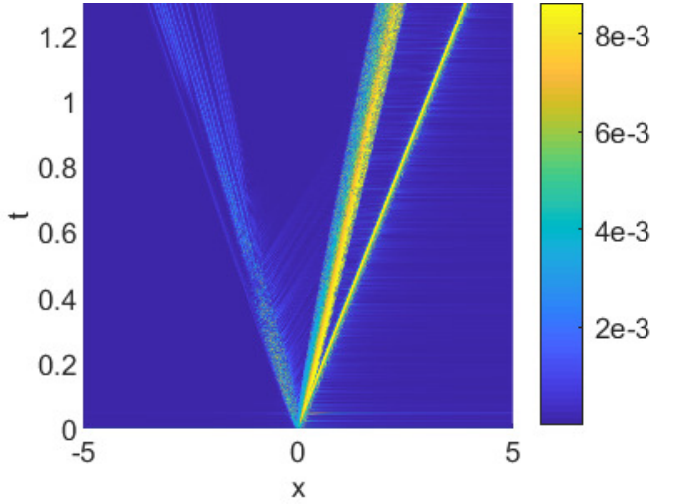}
    \includegraphics{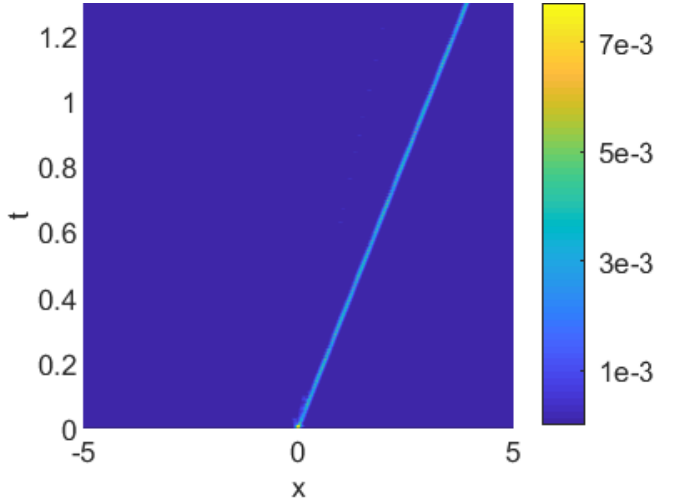}
    \caption{Time history of assigned artificial viscosity for N = 1000, for FC-EV (left) and FC-SDNN (right).}
    \label{fig: Lax viscosity}
  \end{subfigure}        
  \caption{Solutions to the Lax problem produced by the FC-SDNN and
    FC-EV algorithms of order $d = 5$ at $t= 1.3$. Exact solution:
    solid black line. FC-SDNN solution: Blue dashed-line. FC-EV
    solution: Green dot-dashed line. Numbers of discretization points:
    $N=500$ in the upper panels and $N=1000$ in the middle panels.
    Bottom panels: artificial viscosity assignments.}
\label{fig: Lax solution}
\end{figure}

\item[Shu-Osher problem.] 
    We consider the Shu-Osher shock-entropy problem on the interval $[-5, 5]$, with initial condition given by
    \begin{equation} \label{eq: shu osher IC}
        (\rho, u, p)=
        \left\lbrace
            \begin{array}{ccc}
                (3.857143, 2.6929369, 10.33333) & \mbox{if} & x < -4\\
                (1 + 0.2\text{sin}(5x), 0, 1) & \mbox{if} & x > -4.\\
            \end{array}\right.
        \end{equation}
        The solution is computed up to time $T = 1.8$.  In this
        problem, a shock wave encounters an oscillatory smooth
        wavetrain. This test highlights the FC-SDNN solver's low
        dissipation, as artificial viscosity is only assigned in the
        vicinity of the right-traveling shock as long as the waves
        remain smooth, allowing for an accurate representation of the
        smooth features. In particular, Figure~(\ref{SO viscosity})
        shows that the support of the FC-SDNN artificial viscosity is
        much more narrowly confined around the shock position than the
        artificial viscosity resulting in the FC-EV approach. As a
        result, and as illustrated in Figure~(\ref{fig: SO solution})
        (upper and middle right images), the FC-SDNN method provides a
        more accurate resolution in the acoustic region (behind the
        main, rightmost, shock).

\begin{figure}
  \begin{subfigure}[t]{1\linewidth}
    \centering
    \includegraphics[width=.33\linewidth]{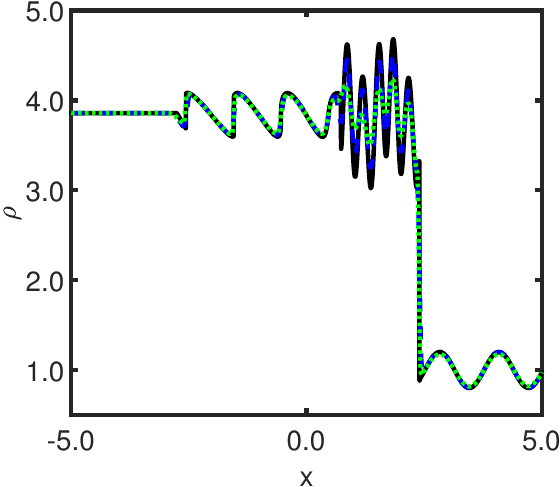}%
    \includegraphics[width=.33\linewidth]{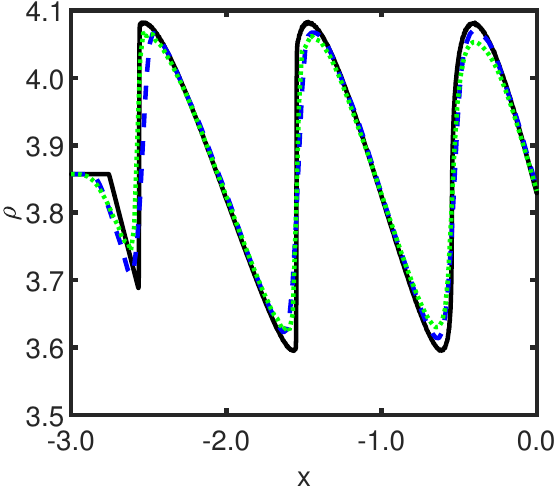}%
    \includegraphics[width=.33\linewidth]{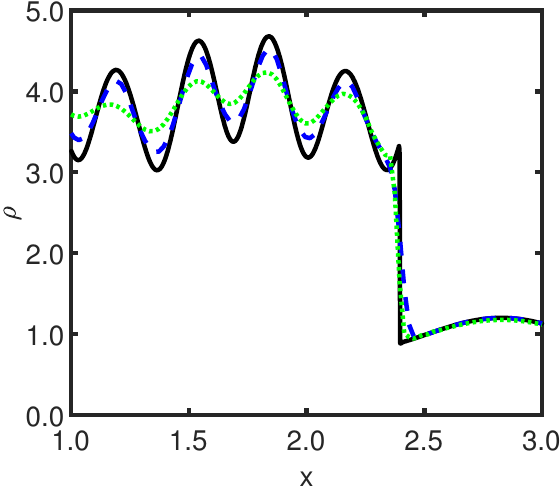}
    \caption{N = 500}
  \end{subfigure}
  \begin{subfigure}[t]{1\linewidth}
    \centering
    \includegraphics[width=.33\linewidth]{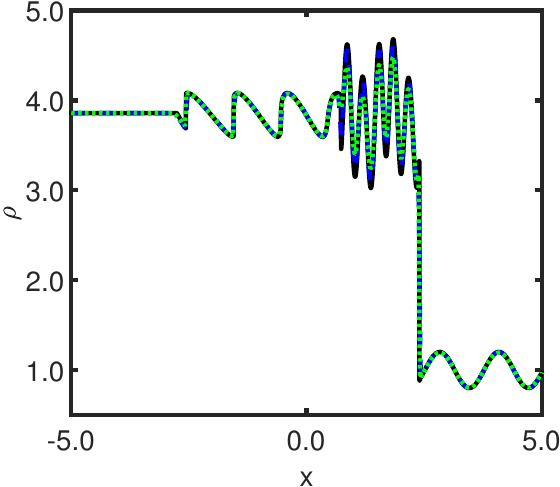}%
    \includegraphics[width=.33\linewidth]{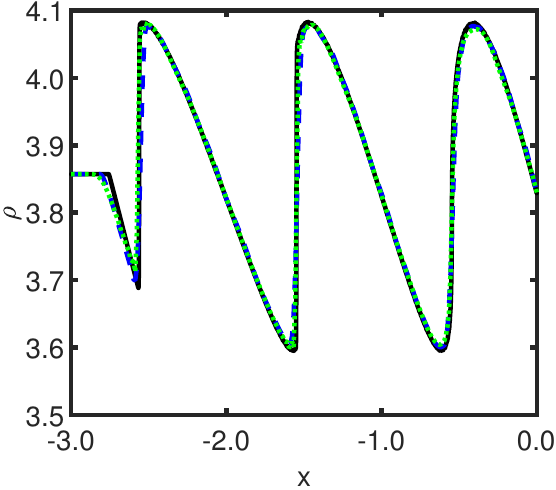}%
    \includegraphics[width=.33\linewidth]{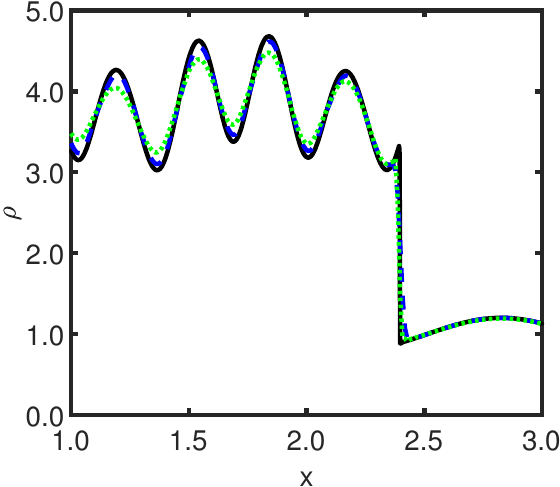}
    \caption{N = 1000}
  \end{subfigure}  
  \begin{subfigure}[t]{1\linewidth}
    \centering
    \includegraphics{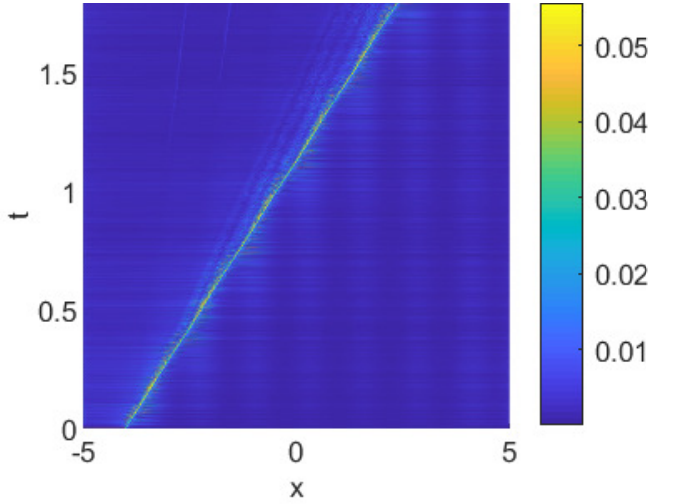}
    \includegraphics{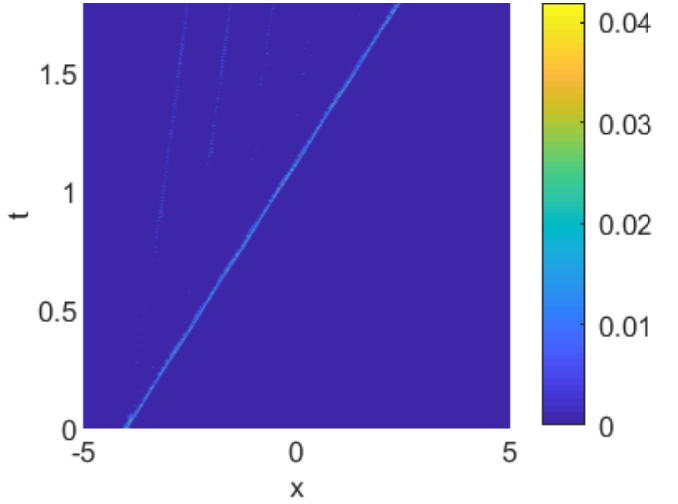}
    \caption{Time history of assigned artificial viscosity for N = 1000, for FC-EV (left) and FC-SDNN (right).}
    \label{SO viscosity}
  \end{subfigure}        
  \caption{Solutions to the Shu-Osher problem produced by the FC-SDNN
    and FC-EV algorithms of order $d = 5$ at $t= 1.8$. Solid black
    line: finely resolved FC-SDNN ($N = 10,000$, for reference). Blue
    dashed line: FC-SDNN with $N=500$ and $N=1000$ (middle
    panels). Green dot-dashed line: FC-EV with $N=500$ (upper panels)
    and $N=1000$ (middle panels).  Bottom panels: artificial viscosity
    assignments.}
\label{fig: SO solution}
\end{figure}

\item[Blast Wave problem.] Finally, we consider the Blast Wave problem
  as presented in~\cite{ray2018artificial}, on the interval $[0, 1]$,
  with initial conditions given by
    \begin{equation} \label{eq:blast wave IC}
        (\rho, u, p)=
        \left\lbrace
            \begin{array}{ccc}
                (1, 0, 1000) & \mbox{if} & x < 0.5\\
                (1, 0, 0.01) & \mbox{if} & x > 0.5,\\
            \end{array}\right.
    \end{equation}
    up to time $T = 0.012$. This setup is similar to the one
    considered in the Sod problem, but with a much stronger
    right-moving shock.  In order to avoid unphysical oscillations at
    the boundaries, which could result from the presence of the strong
    shock, the value of the operator $\tau$ is set to $1$ on the
    leftmost and rightmost nine points in the domain, thus effectively
    assigning a small amount of viscosity at the boundaries at every
    time step of the simulation. As shown in Figure~\ref{fig: BW}, the
    shock is sharply resolved by both the FC-SDNN and FC-EV
    approaches. As the mesh is refined the contact discontinuity is
    resolved more sharply by the former method which, as in the
    previous examples, does not assign viscosity around such features.

    \begin{figure}
     \centering
     \begin{subfigure}[b]{1\textwidth}
         \centering
         \includegraphics[width=0.33\textwidth]{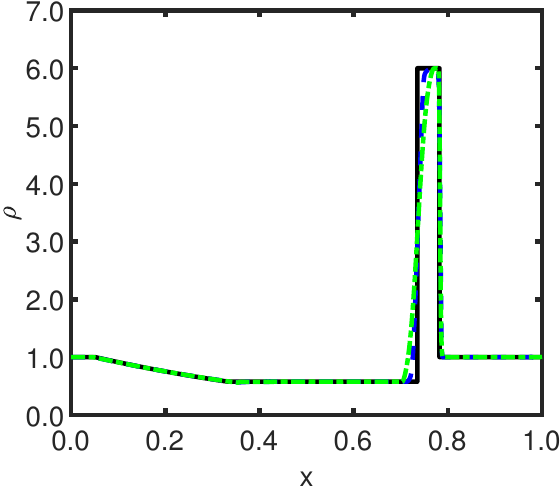}%
         \includegraphics[width=0.33\textwidth]{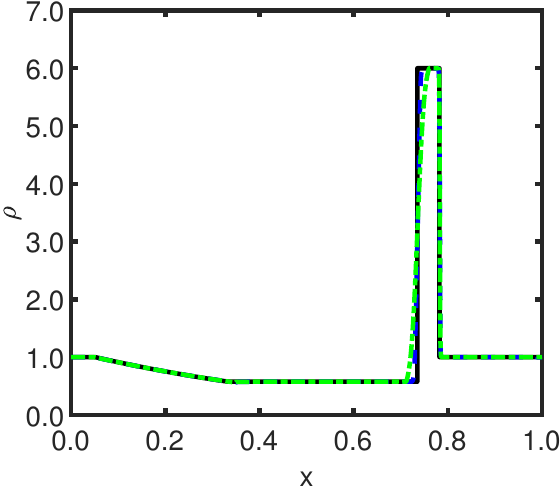}%
         \includegraphics[width=0.33\textwidth]{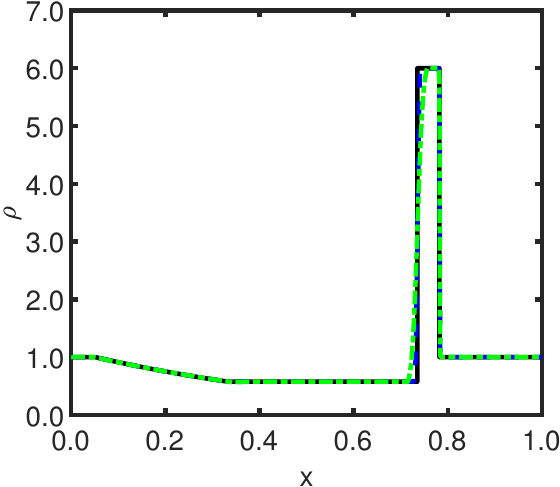}
         \caption{"wide"}
         \label{fig: N = 1000}
     \end{subfigure}%
     \begin{subfigure}[b]{1\textwidth}

     \end{subfigure}%
 
     \begin{subfigure}[b]{1\textwidth}
          \centering
         \includegraphics[width=0.33\textwidth]{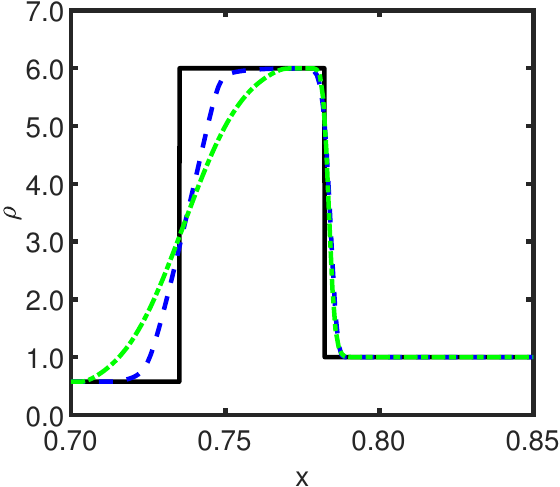}%
         \includegraphics[width=0.33\textwidth]{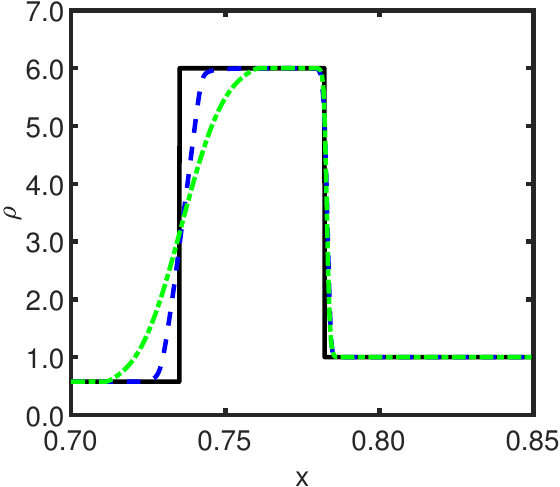}%
         \includegraphics[width=0.33\textwidth]{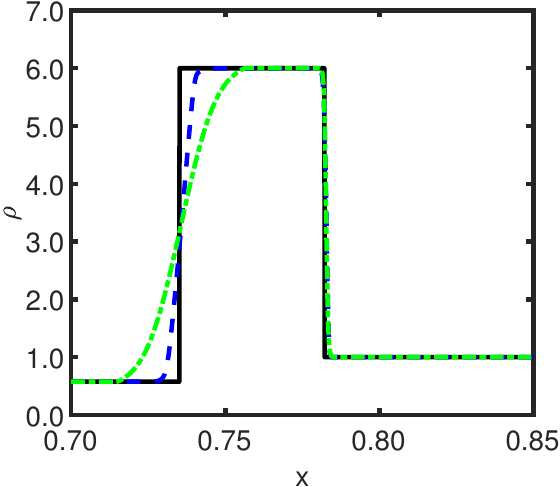}
         \caption{"zoom"}
         \label{fig: N = 2000}   
     \end{subfigure}
    \begin{subfigure}[b]{1\textwidth}
         \centering
         \includegraphics{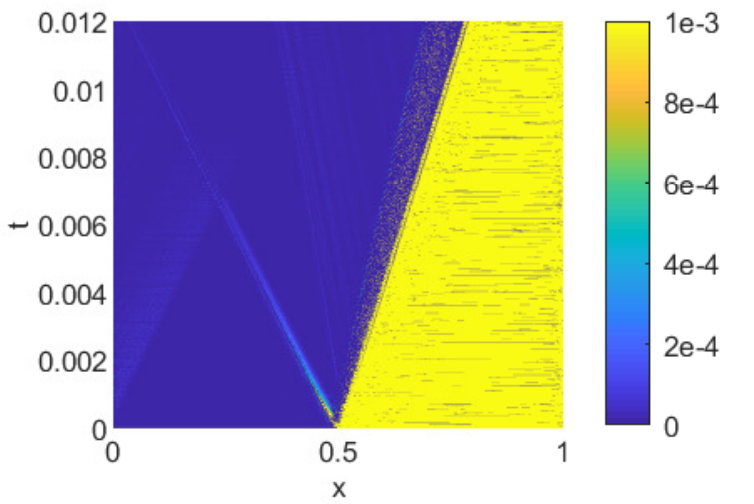}%
         \includegraphics{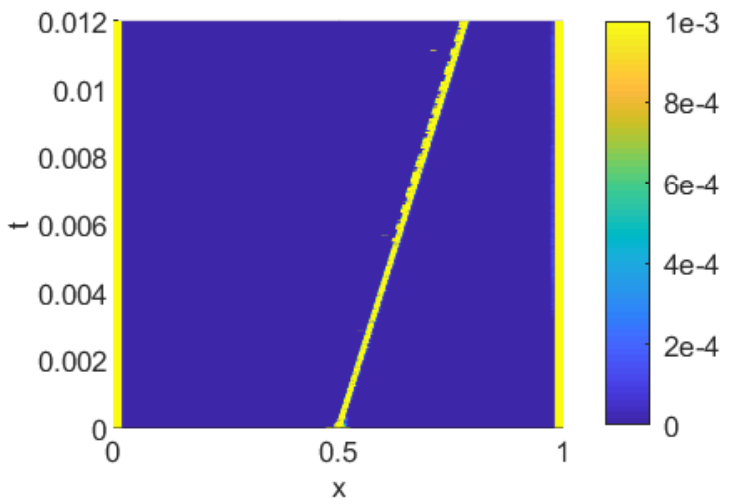}      
         \caption{Time history of assigned artificial viscosity for $N = 1000$, for FC-EV (left) and FC-SDNN (right).}
         \label{fig: BW viscosity} 
    \end{subfigure}     
    \caption{Solutions to the Blast wave problem produced by the
      FC-SDNN and FC-EV algorithms of order $d = 5$ at $t=
      0.012$. Solid black line: exact solution. Blue dashed line:
      FC-SDNN with $N=1000$ (upper- and middle-left panels), $N=2000$
      (upper- and middle-center panels), and $N = 3000$ (upper- and
      middle-right panels). Green dot-dashed line: FC-EV with $N=1000$
      (left panels), $N=2000$ (center panels), and $N = 3000$ (right
      panels). Bottom panels: artificial viscosity assignments.}
     \label{fig: BW} 
    \end{figure}

\end{description}
\subsubsection{\label{Euler 2D results} 2D Euler problems}

The test cases considered in this section showcase the FC-SDNN
method's ability to handle complex shock-shock, as well as
shock-boundary interactions, including shocks moving orthogonally to
the boundaries as in the Riemann 2D and the Shock vortex problems, or
moving obliquely to the boundary of the domain, after reflecting on a
solid wedge, in the Double Mach reflection problem, or reflecting
multiple times on the solid walls of a wind tunnel with a step, in the
Mach 3 forward facing step problem. For the examples in this section
FC expansions with $d = 2$ were used. As a result the FC method can provide
significantly improved accuracy over other approaches of the same or
even higher accuracy orders. In all cases the solutions obtained are
in agreement with solutions obtained previously by various
methods~\cite{woodward1984numerical, lax1998solution,
  guermond2011entropy, mazaheri2019bounded}. As indicated in the
introduction, the proposed approach leads to smooth flows away from
shocks, as evidenced by correspondingly smooth level set lines for the
various flow quantities, in contrast with corresponding results
provided by previous methods.

\begin{description}

\item[Riemann problem (Configuration 4 in~\cite{lax1998solution}).]  We consider a
  Riemann problem configuration on the domain
  $[0, 1.2]\times[0, 1.2]$, with initial conditions given by
    \begin{equation} \label{eq: R4 IC}
        (\rho, u, v, p)=
        \left\lbrace
            \begin{array}{ccc}
                (1.1, 0, 0, 1.1) & \mbox{if} & x \in  [0.6, 1.2]  \mbox{ and } y \in [0.6, 1.2]\\
                (0.5065, 0.8939, 0, 0.35) & \mbox{if} & x \in  [0, 0.6)  \mbox{ and } y \in [0.6, 1.2]\\
                (1.1, 0.8939, 0, 0.35) & \mbox{if} & x \in  [0, 0.6)  \mbox{ and } y \in [0, 0.6)\\
                (0.5065, 0, 0.8939, 0.35) & \mbox{if} & x \in  [0.6, 1.2]  \mbox{ and } y \in [0, 0.6) ,           
            \end{array}\right.
    \end{equation}  
    and with vanishing normal derivatives for all variables on the
    boundary, at all times. The solution is computed up to time
    $T = 0.25$ by means of the FC-SDNN approach. The initial setting
    induces four interacting shock waves, all of which travel
    orthogonally to the straight segments of the domain boundary. The
    results, presented in Figure~\ref{fig: R4 solutions}, show a
    sharpening of the shocks as the mesh is refined and an absence of
    spurious oscillations in all cases. As shown on the left image in
    Figure~\ref{fig: Euler 2D viscosities}, the FC-SDNN viscosity
    assignments are sharply concentrated near the shock positions.
    
    \begin{figure}
    \begin{subfigure}[t]{0.33\linewidth}
        \centering
        \includegraphics[width=1\linewidth]{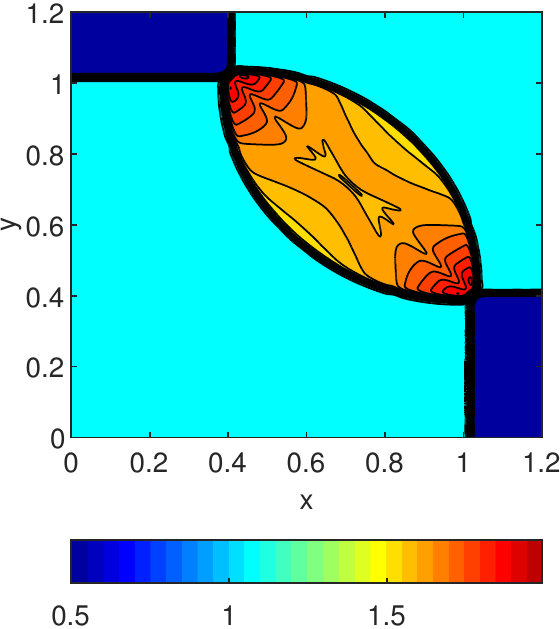}
        \caption{N = 400}
      \end{subfigure}%
      \begin{subfigure}[t]{0.33\linewidth}
        \centering  
        \includegraphics[width=1\linewidth]{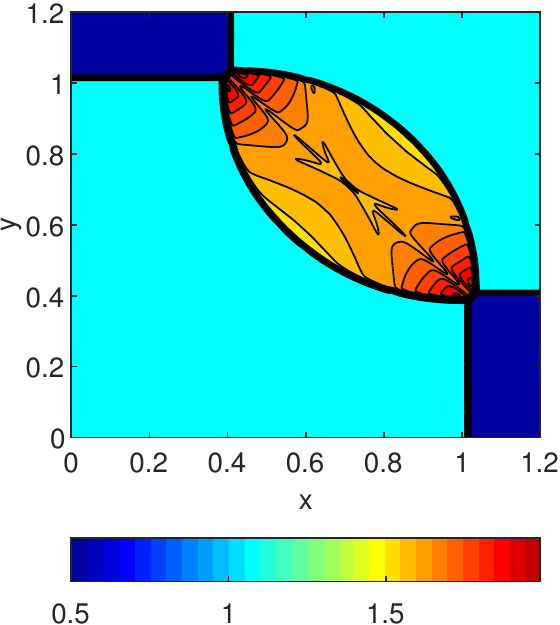}
        \caption{N = 600}
      \end{subfigure}%
      \begin{subfigure}[t]{0.33\linewidth}
        \centering  
        \includegraphics[width=1\linewidth]{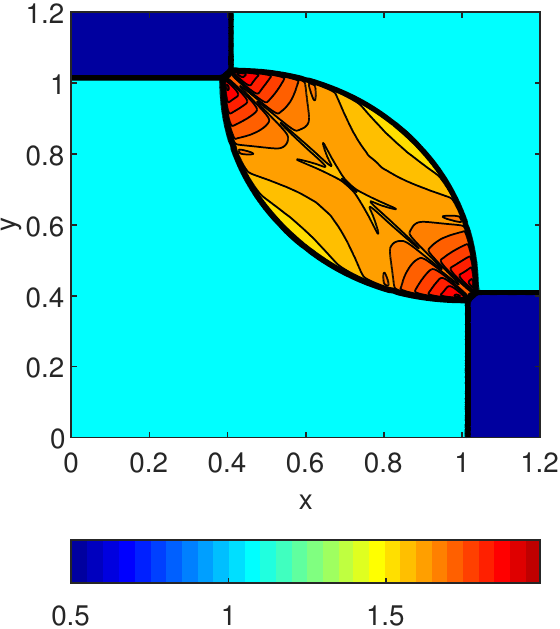}
        \caption{N = 800}
      \end{subfigure} 
      \caption{Second-order ($d = 2$) FC-SDNN numerical solution to
        the Euler 2D Riemann problem considered in Section~\ref{Euler
          2D results}, at $t = 0.25$, obtained by using a spatial
        discretization containing $N \times N$ grid points with three
        different values of $N$, as indicated in each panel. For each
        discretization, the solution is represented using thirty
        equispaced contours between $\rho = 0.5$ and $\rho = 1.99$.}
    \label{fig: R4 solutions}
    \end{figure}

\begin{figure}
    \centering
    \includegraphics{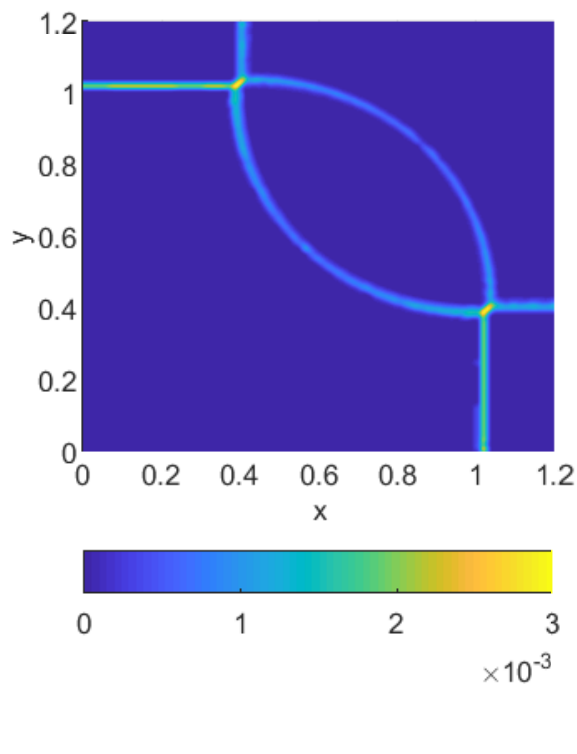}
    \includegraphics{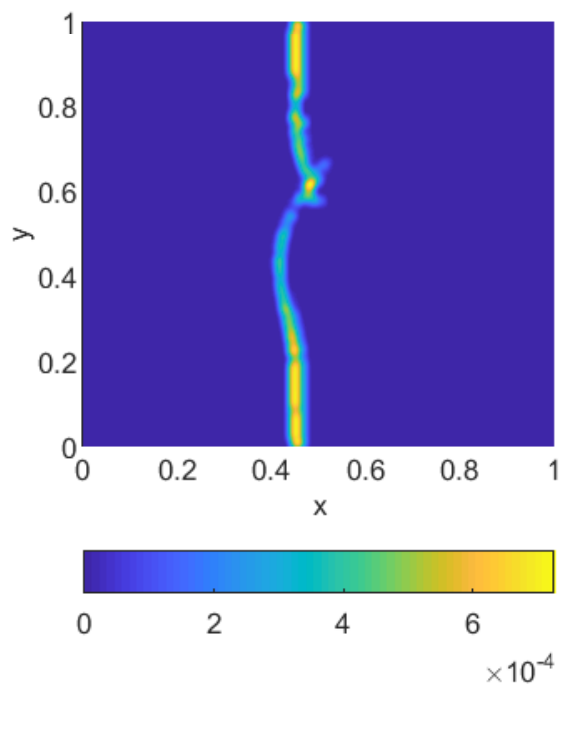}
    \caption{Artificial viscosity profiles for the Euler 2D Riemann
      problem considered in Figure~\ref{fig: R4 solutions} at
      $t = 0.25$ (left) and for the Shock vortex problem
      considered in  Figure~\ref{fig: SV solutions} at $t = 0.35$ (right).}
    \label{fig: Euler 2D viscosities}
\end{figure}

\item[Shock vortex problem~\cite{shu1999high}.]  We next consider a
  ``shock-vortex'' problem in the domain $[0, 1]\times[0,1]$, in which
  a shock wave collides with an isentropic vortex. The initial
  conditions are given by
    \begin{equation} \label{eq: SV IC}
      (\rho, u, v, p)=
        \left\lbrace
            \begin{array}{ccc}
                 (\rho_L + \tilde \rho, u_L + \tilde u, v_L + \tilde v, p_L + \tilde p) & \mbox{if} & x \in  [0, 0.5)  \\
                (\rho_R, u_R, v_R, p_R) & \mbox{if} & x \in  [0.5, 1]          
            \end{array}\right.
    \end{equation}
    where the left state equals the combination of the unperturbed
    left-state $(\rho_L, u_L, v_L, p_L) = (1, \sqrt{\gamma}, 0, 1)$ in
    the shock wave with the isentropic vortex
    \begin{equation}\label{vortex}
    \tilde u = \frac{x - x_c}{r_c}\Phi(r), \quad \tilde v = -\frac{x - x_c}{r_c}\Phi(r), \quad  \frac{\gamma - 1}{4 \zeta \gamma} \Phi(r)^2 =\frac{p_L}{\rho_L} - \frac{p_L + \tilde p}{\rho_L + \tilde \rho}, \quad p_L + \tilde p = (\rho_L + \tilde
    \rho)^{\gamma},
  \end{equation}
  centered at $(x_c, y_c) = (0.25, 0.5)$ (where
  $r = \sqrt{(x - x_c)^2 + (y - y_c)^2}$,
  $\Phi(r) = \epsilon e^{\zeta(1 - (r/r_c)^2)}$).  As in previous
  references for this example we use the vortex parameter values
  $r_c = 0.05$, $\zeta = 0.204$, and $\epsilon = 0.3$. The initial
  right state is given by
    \[
        \rho_R = \rho_L\frac{(\gamma + 1)p_R + \gamma - 1}{(\gamma - 1)p_R + \gamma + 1}, \quad u_R = \sqrt{\gamma} + \sqrt{2}\frac{1 - p_R}{\sqrt{\gamma - 1 + p_R(\gamma - 1)}}, \quad v_R = 0, \quad p_R = 1.3.
      \]
      Vanishing normal derivatives for all variables were imposed on
      the domain boundary at all times.

      The solution was obtained up to time $T = 0.35$, for which
      the vortex has completely crossed the shock.  The solutions
      displayed in Figure~\ref{fig: SV solutions} demonstrate the
      convergence of the method as the spatial and temporal
      discretizations are refined: shocks become sharper with each
      mesh refinement, while the vortex features remain smooth after
      the collision with the shock wave---a property that other
      solvers do not enjoy, and which provides an indicator of the
      quality of the solution. The right image in Figure~\ref{fig:
        Euler 2D viscosities} shows that, as in the previous examples,
      the support of the artificial viscosity imposed by the SDNN
      algorithm is narrowly focused in a vicinity of the shock.

    \begin{figure}
      \begin{subfigure}[b]{0.33\linewidth}
        \centering
        \includegraphics[width=1\linewidth]{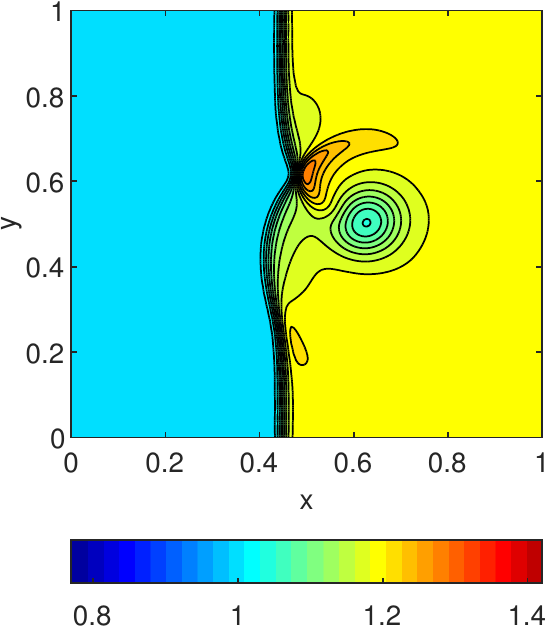}
        \caption{N = 200}
      \end{subfigure}%
      \begin{subfigure}[b]{0.33\linewidth}
        \centering  
        \includegraphics[width=1\linewidth]{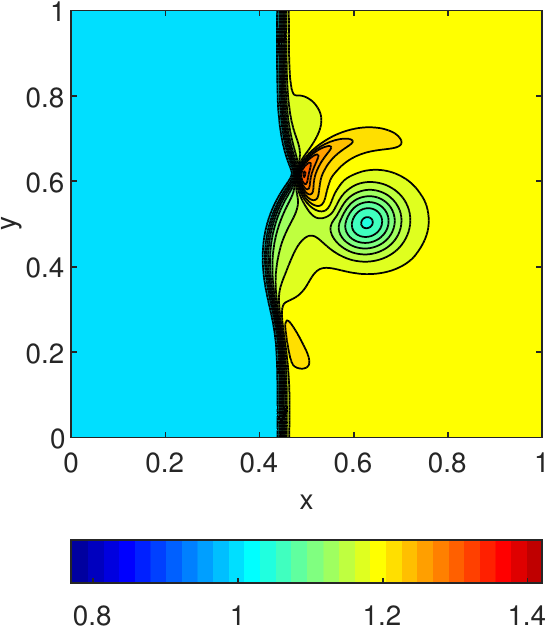}
        \caption{N = 300}
      \end{subfigure}%
      \begin{subfigure}[b]{0.33\linewidth}
        \centering  
        \includegraphics[width=1\linewidth]{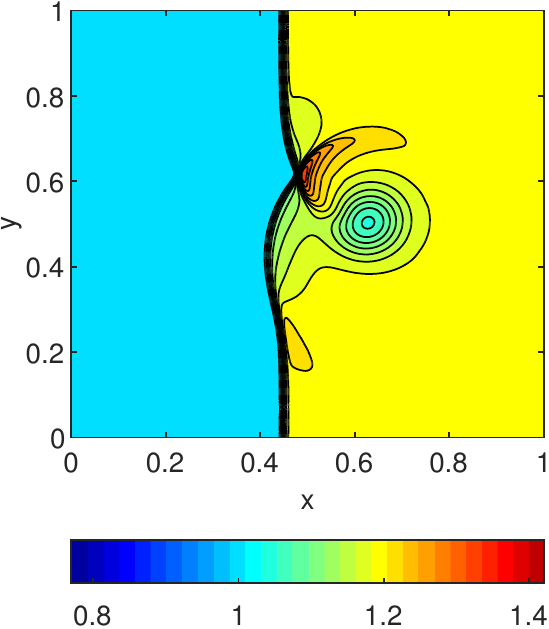}
        \caption{N = 400}
      \end{subfigure} 
      \caption{Second-order ($d = 2$) FC-SDNN numerical solution to
        the Euler 2D Shock-vortex problem considered in
        Section~\ref{Euler 2D results}, at $t = 0.35$, obtained by
        using a spatial discretization containing $N \times N$ grid
        points with three different values of $N$, as indicated in
        each panel. For each discretization, the solution is
        represented using thirty equispaced contours between
        $\rho = 0.77$ and $\rho = 1.42$.}
    \label{fig: SV solutions}
    \end{figure}

  \item[Mach 3 forward facing step~\cite{woodward1984numerical}.] We now consider a ``Mach 3
    forward facing step problem'' on the domain
    $([0, 0.6]\times[0, 1]) \cup ([0.6, 3]\times[0.2, 1])$, in which a
    uniform Mach 3 flow streams through a wind tunnel with a forward
    facing step, of 0.2 units in height, located at $x = 0.6$. The
    initial condition is given by
\begin{equation} \label{eq: M3 IC}
(\rho, u, v, p)=
\left\lbrace
    \begin{array}{ccc}
        (1.4, 3, 0, 1) & \mbox{if} & (x, y) \in  ([0, 0.6)\times[0, 1]) \cup ([0.6, 3]\times[0.2, 1])   \\
        (1.4, 0, 0, 1) & \mbox{if} & (x, y) \in  \{0.6\} \times [0, 0.2],
    \end{array}\right.
\end{equation}
and the solution is computed up to time $T = 4$. An inflow condition
is imposed at the left boundary at all times which coincides with the
initial values on that boundary, while the equations are evolved at
the outflow right boundary. Reflecting boundary conditions (zero
normal velocity) are applied at all the other
boundaries. (Following~\cite{woodward1984numerical,
  guermond2011entropy}, no boundary condition is enforced at the node
located at the step corner.) The simulation was performed using the
adaptive time step~(\ref{eq: CFL}) with $\textrm{CFL} = 1$. The
density solution is displayed in Figure~\ref{M3 solutions}.

\item[Double Mach reflection~\cite{woodward1984numerical}]
  We finally consider the ``Double Mach reflection problem'' on the
  domain $[0, 4] \times [0, 1]$. This problem contains a reflective
  wall located on the $x \geq x_r$ part of the bottom boundary $y = 0$
  (here we take $x_r = 1/6$), upon which there impinges an incoming
  shock wave forming a $\theta = \pi/3$ angle with the positive
  $x$-axis. The initial condition is given by
\begin{equation} \label{eq: DM IC}
    \mathbf{e}(x, y, 0) = (\rho, \rho u, \rho v, E)=
    \left\lbrace
        \begin{array}{ccc}
             (8, 57.1597, -33.0012 , 563.544) & \mbox{if} &  0 \leq x \leq  x_r + \frac{y}{\sin(\theta)}  \\
            (1.4, 0, 0, 2.5) & \mbox{if} & x_r + \frac{y}{\sin(\theta)} \leq x \leq  4          
        \end{array}\right.
\end{equation}
and the solution is computed up to time $T = 0.2$.  This setup,
introduced in~\cite{woodward1984numerical}, gives rise to the
reflection of a strong oblique shock wave on a wall. (Equivalently,
upon counter-clockwise rotation by $30^\circ$, this setup can be
interpreted as a vertical shock impinging on a $30^\circ$ ramp.) As a
result of the shock reflection a number of flow features arise,
including, notably, two Mach stems and two contact discontinuities
(slip lines), as further discussed below in this section.

The initial and boundary values used in this context do not exactly
coincide with the ones utilized in~\cite{woodward1984numerical}. Indeed, on one hand, in order to avoid density
oscillations near the intersection of the shock and the top
computational boundary, we utilize ``oblique'' Neumann boundary
conditions on all flow variables
$\mathbf{e} = (\rho, \rho u,\rho v,E)$. More precisely, we enforce
zero values on the derivative of $\mathbf{e}$ with respect to the
direction parallel to the shock, along both the complete upper
computational boundary and the region $0\leq x\leq x_r$ (that is, left
of the ramp) on the lower computational boundary. This method can be
considered as a further development of the approach proposed
in~\cite{vevek2019alternative}, wherein an extended domain in the oblique direction
was utilized in conjunction with Neumann conditions along the normal
direction to the oblique boundary. In order to incorporate the oblique
Neumann boundary condition we utilized the method described
in~\cite{amlani2016fc}, which, in the present application, proceeds by
obtaining relations between oblique, normal and tangential derivatives
of the flow variables $\mathbf{e}$. Inflow boundary conditions were
used which prescribe time-independent values of $\rho$, $u$ and $v$ on
the left boundary; outflow conditions on the right boundary, in turn,
enforce a time independent value of $p$.\looseness = -1
\begin{figure}
  \begin{subfigure}[t]{1\linewidth}
    \centering
    \includegraphics[width=0.7\linewidth]{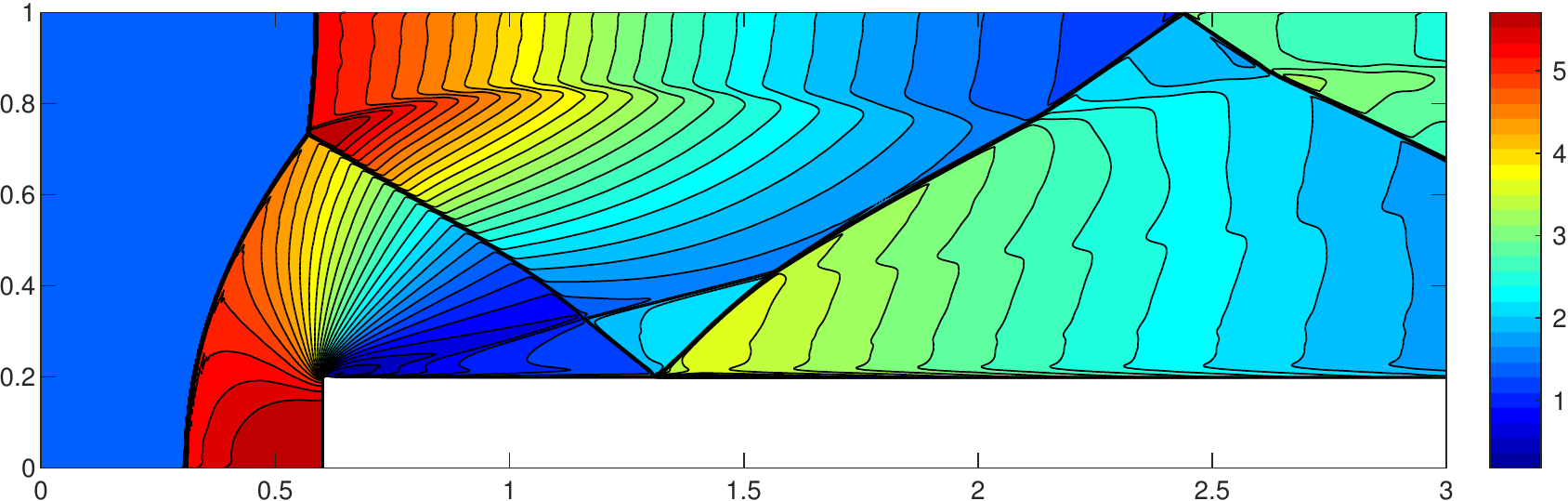}
    \caption{Thirty equispaced $\rho$ contours between 0.18 and 5.71\label{M3 solutions_a}}
  \end{subfigure} 
  \begin{subfigure}[t]{1\linewidth}
    \centering  
    \includegraphics[width=0.7\linewidth]{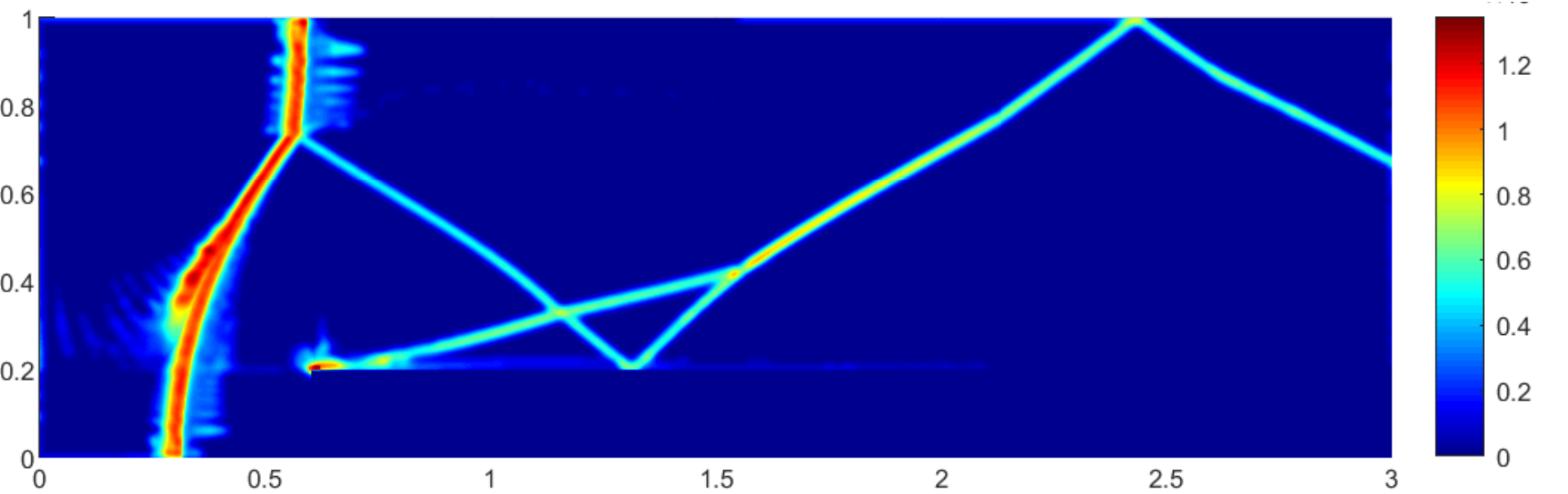}
    \caption{Artificial viscosity at time t=4.}
  \end{subfigure}

  \caption{Top panel: Second-order ($d = 2$) FC-SDNN numerical
    solution to the Mach 3 forward facing-step problem considered in
    Section~\ref{Euler 2D results}, at $t = 4$, obtained by using
    $1200 \times 400$ spatial grid. Bottom panel: Viscosity
    assignment.}
\label{M3 solutions}
\end{figure}

Additionally, following~\cite{vevek2019alternative}, we utilize a numerical viscous
incident shock as an initial condition in order to avoid well-known
post-shock oscillations, as noted in~\cite{johnsen2008numerical}, that result from
the use of a sharp initial profile. Such a smeared shock is obtained
in our context by applying the FC-SDNN solver to the propagation of an
oblique flat shock on all of space (without the ramp) up to time
$T = 0.2$, including imposition of oblique Neumann boundary conditions
throughout the top and bottom boundary. The solution
$\hat{\mathbf{e}}_h$ amounts to a smeared shock profile on the $(x,y)$
plane which, at time $t$, is centered on the straight line
$x = x_s(y,t)$ where
$x_s(y,t) = x_r + \frac{U_s}{\sin (\theta)}t + \frac{y}{\tan
  (\theta)}$. This shock is followed by some back-trailing
oscillations, as has been observed in~\cite{johnsen2008numerical, vevek2019alternative}. In order to
eliminate these artifacts from the initial condition, a diagonal strip
of $N_d$ points centered around the shock location is selected. Then,
on a strip of $N_d$ points surrounding the initial position of the
shock, the initial condition $\mathbf{e}_h$ to the Double Mach
reflection problem is defined as
\[
    \mathbf{e}_h(x, y, 0) = q_{c, r}(x_s(y, 0) - x) \hat{\mathbf{e}}_h(x+x_s(y, T), y, T) +  (1 - q_{c, r}(x_s(y, 0)-x)) \mathbf{e}(x, y, 0)
\]
where $q_{c, r}$ is the window function defined in~(\ref{eq:
  qf}). $N_d$ is taken to be small enough to exclude the back-trailing
oscillations from the strip, and large enough as to provide a
relatively smooth profile for the initial condition. For the
double-Mach solution depicted in Figure~\ref{DM solutions}, which was
obtained on the basis of a $3200\times 800$-point grid, the
incident-shock parameter values $N_d = 125$, $c = 25$ and $r = 50$
were used.
\begin{figure}
  \begin{subfigure}[t]{1\linewidth}
    \centering
    \includegraphics[width=0.7\linewidth]{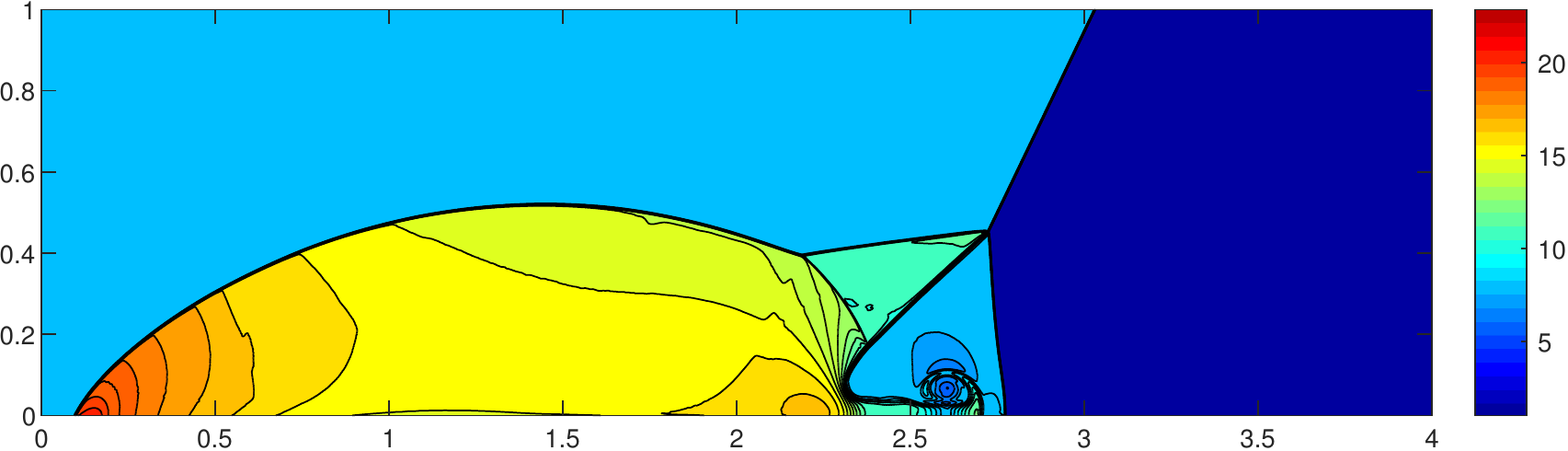}
    \caption{Thirty equispaced $\rho$ contours between 1.0 and 22.8}
  \end{subfigure} 
  \begin{subfigure}[t]{1\linewidth}
    \centering  
    \includegraphics[width=0.7\linewidth]{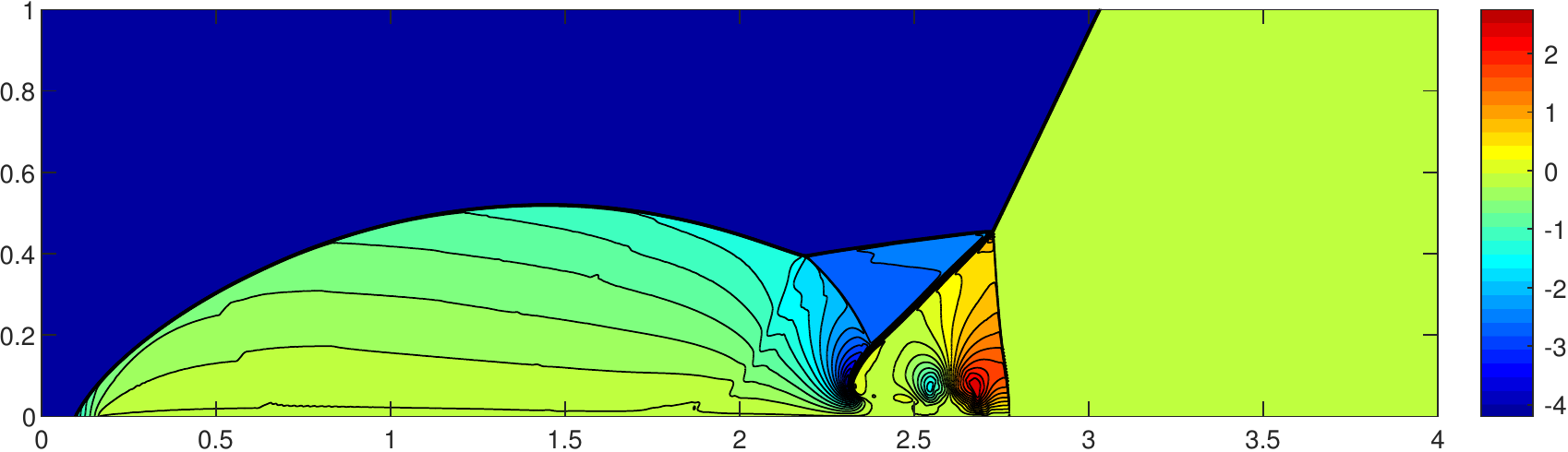}
    \caption{Thirty equispaced $v$ contours between -4.2 and 2.7}
  \end{subfigure}
  \begin{subfigure}[t]{1\linewidth}
    \centering  
    \includegraphics[width=0.7\linewidth]{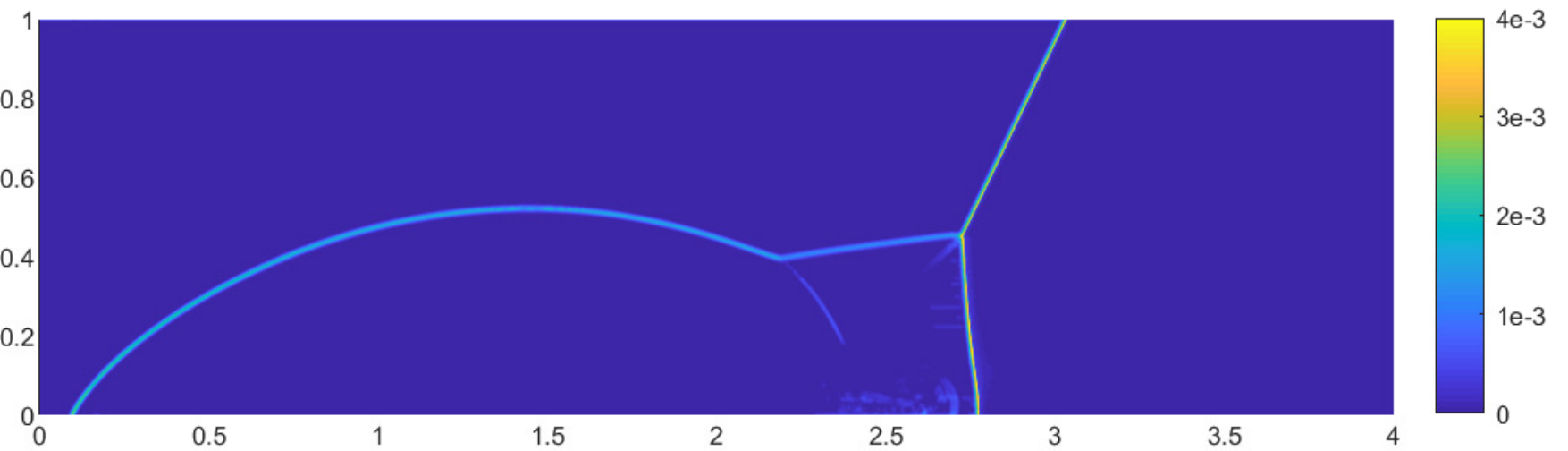}
    \caption{Artificial viscosity at time $t = 0.2$}
  \end{subfigure}

  \caption{Top and middle panels: Second-order ($d = 2$) FC-SDNN
    numerical solution to the Double Mach reflection problem
    considered in Section~\ref{Euler 2D results}, at $t = 0.2$,
    obtained by using $3200 \times 800$ spatial grid. Bottom panel:
    Viscosity assignment.}
\label{DM solutions}
\end{figure}

Among several notable features in the $t = 0.2$ solution we mention
the density-contour roll-ups of the primary slip line that are clearly
visible in the density component of the solution presented in
Figure~\ref{DM solutions}a near the reflective wall around $x=2.6$, as
well as the extremely weak secondary slip line, that is most easily
noticeable in the variable $v$ (Figure~\ref{DM solutions}b) as a dip
in the contour lines within the rectangle $[1.9,2.2]\times [0,0.4]$,
but which is also visible in this area as a blip in the density
contour lines (Figure~\ref{DM solutions}a). The accurate simulation of
these two features has typically been found 
challenging~\cite{woodward1984numerical,kemm2016proper,
  vevek2019alternative}.  The elimination of back-trailing shock
oscillations, which was accomplished, as discussed above, by resorting
to use of a numerically viscous incident shock, is instrumental in the
resolution of the secondary slip line---which might otherwise be
polluted by such oscillations to the point of being unrecognizable.

\end{description}

\section{\label{sec:conclusion}Conclusions}

This paper introduced the FC-SDNN method, a neural network-based
artificial viscosity method for evaluation of shock dynamics in
non-periodic domains. In smooth flow areas the method enjoys the
essentially dispersionless character inherent in the FC method. The
smooth but localized viscosity assignments allow for a sharp
resolution of shocks and contact discontinuities, while yielding
smooth flow profiles away from jump discontinuities. An efficient
implementation for general 2D and 3D domains, which could be pursued
on the basis of an overlapping-patch setup~\cite{albin2011spectral,
  brown1997overture}, is left for future work.\looseness = -1

\vspace{0.3cm}
\noindent {\bf \large Acknowledgments.} OB and DL gratefully acknowledges
support from NSF under contracts DMS-1714169 and DMS-2109831, from
AFOSR under contract FA9550-21-1-0373, and from the NSSEFF Vannevar
Bush Fellowship under contract number N00014-16-1-2808. Initial
conversations with J. Paul are gratefully acknowledged.

\small

\bibliographystyle{abbrv}
\bibliography{bibliography}

\begin{thebibliography}{10}

\bibitem{albin2011spectral}
N.~Albin and O.~P. Bruno.
\newblock A spectral {FC} solver for the compressible {N}avier--{S}tokes
  equations in general domains {I}: {E}xplicit time-stepping.
\newblock {\em Journal of Computational Physics}, 230:6248--6270, 2011.

\bibitem{amlani2016fc}
F.~Amlani and O.~P. Bruno.
\newblock An {FC}-based spectral solver for elastodynamic problems in general
  three-dimensional domains.
\newblock {\em Journal of Computational Physics}, 307:333--354, 2016.

\bibitem{brown1997overture}
D.~L. Brown, W.~D. Henshaw, and D.~J. Quinlan.
\newblock Overture: An object-oriented framework for solving partial
  differential equations.
\newblock In {\em International Conference on Computing in Object-Oriented
  Parallel Environments}, pages 177--184. Springer, 1997.

\bibitem{bruno2019higher}
O.~P. Bruno, M.~Cubillos, and E.~Jimenez.
\newblock Higher-order implicit-explicit multi-domain compressible
  {N}avier-{S}tokes solvers.
\newblock {\em Journal of Computational Physics}, 391:322--346, 2019.

\bibitem{bruno2010high}
O.~P. Bruno and M.~Lyon.
\newblock High-order unconditionally stable {FC-AD} solvers for general smooth
  domains {I}. {B}asic elements.
\newblock {\em Journal of Computational Physics}, 229(6):2009--2033, 2010.

\bibitem{bruno2020two}
O.~P. Bruno and J.~Paul.
\newblock Two-dimensional {F}ourier {C}ontinuation and applications.
\newblock {\em arXiv:2010.03901}, 2020.

\bibitem{carpenter1995theoretical}
M.~H. Carpenter, D.~Gottlieb, S.~Abarbanel, and W.-S. Don.
\newblock The theoretical accuracy of {R}unge--{K}utta time discretizations for
  the initial boundary value problem: a study of the boundary error.
\newblock {\em SIAM Journal on Scientific Computing}, 16(6):1241--1252, 1995.

\bibitem{discacciati2020controlling}
N.~Discacciati, J.~S. Hesthaven, and D.~Ray.
\newblock Controlling oscillations in high-order {D}iscontinuous {G}alerkin
  schemes using artificial viscosity tuned by neural networks.
\newblock {\em Journal of Computational Physics}, 409:109304, 2020.

\bibitem{gentry1966eulerian}
R.~A. Gentry, R.~E. Martin, and B.~J. Daly.
\newblock An {E}ulerian differencing method for unsteady compressible flow
  problems.
\newblock {\em Journal of Computational Physics}, 1(1):87--118, 1966.

\bibitem{glorot2010understanding}
X.~Glorot and Y.~Bengio.
\newblock Understanding the difficulty of training deep feedforward neural
  networks.
\newblock In {\em JMLR Workshop and Conference Proceedings}, pages 249--256,
  2010.

\bibitem{goodfellow2016}
I.~Goodfellow, Y.~Bengio, and A.~Courville.
\newblock {\em Deep {L}earning}.
\newblock MIT Press, 2016.

\bibitem{gottlieb2005high}
S.~Gottlieb.
\newblock On high order strong stability preserving {R}unge-{K}utta and multi
  step time discretizations.
\newblock {\em Journal of Scientific Computing}, 25(1):105--128, 2005.

\bibitem{guermond2011entropy}
J.-L. Guermond, R.~Pasquetti, and B.~Popov.
\newblock Entropy viscosity method for nonlinear conservation laws.
\newblock {\em Journal of Computational Physics}, 230(11):4248--4267, 2011.

\bibitem{harten1987uniformly}
A.~Harten, B.~Engquist, S.~Osher, and S.~R. Chakravarthy.
\newblock Uniformly high order accurate essentially non-oscillatory schemes,
  {III}.
\newblock In {\em Upwind and high-resolution schemes}, pages 218--290.
  Springer, 1987.

\bibitem{hirsch1990numerical}
C.~Hirsch.
\newblock Numerical computation of internal and external flows. {V}ol.
  2-{C}omputational methods for inviscid and viscous flows.
\newblock {\em John Wiley \& Sons, 1990, 708}, 1990.

\bibitem{jiang1996efficient}
G.-S. Jiang and C.-W. Shu.
\newblock Efficient implementation of weighted {ENO} schemes.
\newblock {\em Journal of Computational Physics}, 126(1):202--228, 1996.

\bibitem{johnsen2008numerical}
E.~Johnsen and S.~Lele.
\newblock Numerical errors generated in simulations of slowly moving shocks.
\newblock {\em Center for Turbulence Research, Annual Research Briefs}, pages
  1--12, 2008.

\bibitem{kemm2016proper}
F.~Kemm.
\newblock On the proper setup of the double {M}ach reflection as a test case
  for the resolution of gas dynamics codes.
\newblock {\em Computers \& Fluids}, 132:72--75, 2016.

\bibitem{kopriva2009implementing}
D.~A. Kopriva.
\newblock {\em Implementing spectral methods for partial differential
  equations: {A}lgorithms for scientists and engineers}.
\newblock Springer Science \& Business Media, 2009.

\bibitem{kornelus2018flux}
A.~Kornelus and D.~Appel{\"o}.
\newblock Flux-conservative {H}ermite methods for simulation of nonlinear
  conservation laws.
\newblock {\em Journal of Scientific Computing}, 76(1):24--47, 2018.

\bibitem{lapidus1967detached}
A.~Lapidus.
\newblock A detached shock calculation by second-order finite differences.
\newblock {\em Journal of Computational Physics}, 2(2):154--177, 1967.

\bibitem{lax1959systems}
P.~Lax.
\newblock Systems of conservation laws.
\newblock Technical report, LOS ALAMOS NATIONAL LAB NM, 1959.

\bibitem{lax1954weak}
P.~D. Lax.
\newblock Weak solutions of nonlinear hyperbolic equations and their numerical
  computation.
\newblock {\em Communications on pure and applied mathematics}, 7(1):159--193,
  1954.

\bibitem{lax1998solution}
P.~D. Lax and X.-D. Liu.
\newblock Solution of two-dimensional {R}iemann problems of gas dynamics by
  positive schemes.
\newblock {\em SIAM Journal on Scientific Computing}, 19(2):319--340, 1998.

\bibitem{leveque1992numerical}
R.~J. LeVeque.
\newblock {\em Numerical methods for conservation laws}, volume 132.
\newblock Springer, 1992.

\bibitem{leveque2002finite}
R.~J. LeVeque et~al.
\newblock {\em Finite volume methods for hyperbolic problems}, volume~31.
\newblock 2002.

\bibitem{liu1994weighted}
X.-D. Liu, S.~Osher, and T.~Chan.
\newblock Weighted essentially non-oscillatory schemes.
\newblock {\em Journal of Computational Physics}, 115(1):200--212, 1994.

\bibitem{lyon2010high}
M.~Lyon and O.~P. Bruno.
\newblock High-order unconditionally stable {FC-AD} solvers for general smooth
  domains {II}. {E}lliptic, parabolic and hyperbolic {PDE}s; theoretical
  considerations.
\newblock {\em Journal of Computational Physics}, 229(9):3358--3381, 2010.

\bibitem{mazaheri2019bounded}
A.~Mazaheri, C.-W. Shu, and V.~Perrier.
\newblock Bounded and compact weighted essentially nonoscillatory limiters for
  discontinuous {G}alerkin schemes: {T}riangular elements.
\newblock {\em Journal of Computational Physics}, 395:461--488, 2019.

\bibitem{pathria1997correct}
D.~Pathria.
\newblock The {C}orrect {F}ormulation of {I}ntermediate {B}oundary {C}onditions
  for {R}unge--{K}utta {T}ime {I}ntegration of {I}nitial {B}oundary {V}alue
  {P}roblems.
\newblock {\em SIAM Journal on Scientific Computing}, 18(5):1255--1266, 1997.

\bibitem{persson2006sub}
P.-O. Persson and J.~Peraire.
\newblock Sub-cell shock capturing for discontinuous {G}alerkin methods.
\newblock In {\em 44th AIAA Aerospace Sciences Meeting and Exhibit}, page 112,
  2006.

\bibitem{ramani2019space1}
R.~Ramani, J.~Reisner, and S.~Shkoller.
\newblock A space-time smooth artificial viscosity method with wavelet noise
  indicator and shock collision scheme, {P}art 1: {T}he {1-D} case.
\newblock {\em Journal of Computational Physics}, 387:81--116, 2019.

\bibitem{ramani2019space2}
R.~Ramani, J.~Reisner, and S.~Shkoller.
\newblock A space-time smooth artificial viscosity method with wavelet noise
  indicator and shock collision scheme, {P}art 2: the {2-D} case.
\newblock {\em Journal of Computational Physics}, 387:45--80, 2019.

\bibitem{ray2018artificial}
D.~Ray and J.~S. Hesthaven.
\newblock An artificial neural network as a troubled-cell indicator.
\newblock {\em Journal of Computational Physics}, 367:166--191, 2018.

\bibitem{reisner2013space}
J.~Reisner, J.~Serencsa, and S.~Shkoller.
\newblock A space--time smooth artificial viscosity method for nonlinear
  conservation laws.
\newblock {\em Journal of Computational Physics}, 235:912--933, 2013.

\bibitem{richtmyer1948proposed}
R.~Richtmyer.
\newblock Proposed numerical method for calculation of shocks.
\newblock {\em Los Alamos Report}, 671, 1948.

\bibitem{schwander2021controlling}
L.~Schwander, D.~Ray, and J.~S. Hesthaven.
\newblock Controlling oscillations in spectral methods by local artificial
  viscosity governed by neural networks.
\newblock {\em Journal of Computational Physics}, 431:110144, 2021.

\bibitem{shahbazi2011multi}
K.~Shahbazi, N.~Albin, O.~P. Bruno, and J.~S. Hesthaven.
\newblock Multi-domain {F}ourier-continuation/{WENO} hybrid solver for
  conservation laws.
\newblock {\em Journal of Computational Physics}, 230(24):8779--8796, 2011.

\bibitem{shu1999high}
C.-W. Shu.
\newblock High order {ENO} and {WENO} schemes for computational fluid dynamics.
\newblock In {\em High-order methods for computational physics}, pages
  439--582. Springer, 1999.

\bibitem{shu1989efficient}
C.-W. Shu and S.~Osher.
\newblock Efficient implementation of essentially non-oscillatory
  shock-capturing schemes, {II}.
\newblock In {\em Upwind and High-Resolution Schemes}, pages 328--374.
  Springer, 1989.

\bibitem{sod1978survey}
G.~A. Sod.
\newblock A survey of several finite difference methods for systems of
  nonlinear hyperbolic conservation laws.
\newblock {\em Journal of Computational physics}, 27(1):1--31, 1978.

\bibitem{stevens2020enhancement}
B.~Stevens and T.~Colonius.
\newblock Enhancement of shock-capturing methods via machine learning.
\newblock {\em Theoretical and Computational Fluid Dynamics}, 34:483--496,
  2020.

\bibitem{vevek2019alternative}
U.~Vevek, B.~Zang, and T.~H. New.
\newblock On alternative setups of the double {M}ach reflection problem.
\newblock {\em Journal of Scientific Computing}, 78(2):1291--1303, 2019.

\bibitem{vonneumann1950method}
J.~VonNeumann and R.~D. Richtmyer.
\newblock A method for the numerical calculation of hydrodynamic shocks.
\newblock {\em Journal of Applied Physics}, 21(3):232--237, 1950.

\bibitem{woodward1984numerical}
P.~Woodward and P.~Colella.
\newblock The numerical simulation of two-dimensional fluid flow with strong
  shocks.
\newblock {\em Journal of Computational Physics}, 54(1):115--173, 1984.

\end{thebibliography}

\end{document}